\documentclass[11pt,reqno]{amsart}

\usepackage{tikz}
\usetikzlibrary{calc}
\usepackage{wasysym,hyperref}
\usepackage{accents,esint,amssymb}
\usepackage{tableaucolors}

\newtheorem{theorem}{Theorem}[section]
\newtheorem{lemma}[theorem]{Lemma}
\newtheorem{corollary}[theorem]{Corollary}

\def\t{\theta}

\def\T{\Theta}
\def\sg{\sigma}
\def\vs{\varsigma}

\def\a{\alpha}
\def\d{\delta}
\def\w{\omega}
\def\g{\gamma}
\def\k{\kappa}
\def\G{\Gamma}
\def\c{{\rm c}}
\def\s{{\rm s}}
\def\tf{{\rm t}}

\def\del{\partial}
\def\dot{\accentset{\mbox{\Large\bfseries .}}}

\def\ol{\overline}
\def\ul{\underline}
\def\vp{\varphi}
\def\B#1{\mathbb#1}
\def\mc#1{\mathcal#1}
\def\com#1{\quad\text{#1}\quad}
\def\TS{\textstyle}
\def\wh{\widehat}
\def\wt{\widetilde}
\def\whc#1{\wh{\mc{#1}}}

\def\IR#1{\scalebox{1.1}{$\frac{\mc I#1\mc R}2$}}
\def\pb{{\TS\frac\pi2}}
\def\piot{{\TS\frac{2\pi}T}}
\def\piotk{{\TS\frac{2\pi}{T_k}}}
\def\x#1{$#1\times #1$}
\def\({\left(\begin{array}{cccccc}}
\def\){\end{array}\right)}
\def\dot{\accentset{\mbox{\Large\bfseries .}}}
\def\ddot{\accentset{\mbox{\Large\bfseries .\kern-1.75pt.}}}
\def\bpr{\backprime}
\def\half{{\TS\frac12}}

\numberwithin{equation}{section}

\begin{document}

\title[Time-periodic solutions of compressible Euler]{Time-periodic
  Solutions of the Compressible Euler Equations and the Nonlinear
  Theory of Sound}

\author{Blake Temple}
\address{Department of Mathematics,
  University of California, Davis, CA 95616}
\email{temple@math.ucdavis.edu}

\author{Robin Young}
\address{Department of Mathematics and Statistics,
  University of Massachusetts, Amherst, MA 01003}
\email{rcy@umass.edu}

\date{\today\ -- \jobname}

\begin{abstract}
  We prove the existence of ``pure tone'' nonlinear sound waves of all
  frequencies.  These are smooth, time periodic, oscillatory solutions
  of the $3\times3$ compressible Euler equations satisfying periodic
  or acoustic boundary conditions in one space dimension.  This
  resolves a centuries old problem in the theory of Acoustics, by
  establishing that the pure modes of the linearized equations are the
  small amplitude limits of solutions of the nonlinear equations.
  Riemann's celebrated 1860 proof that compressions always form shocks
  is known to hold for isentropic and barotropic flows, but our proof
  shows that for generic entropy profiles, shock-free periodic
  solutions containing nontrivial compressions and rarefactions exist
  for every wavenumber $k$.
\end{abstract}

\maketitle

\section{Introduction}

The theory of Shock Waves was initiated in 1848 by Stokes in his paper
\emph{``On a Difficulty in the Theory of Sound''} \cite{Stokes}.
Stokes confirmed an observation of Challis that oscillatory solutions
of the nonlinear compressible Euler equations,
\begin{equation}
  \label{euler}
  \del_t\rho + \nabla\cdot\big(\rho\,u\big) = 0, \qquad
  \del_t\big(\rho\,u\big) +
  \nabla\cdot\big(\rho\,u\,u^\dag + p\,I\big) = 0,
\end{equation}
with pressure $p=p(\rho)$, exhibit gradient blowup.  It follows that
the nonlinear equations generate a qualitatively different theory of
sound compared to the sinusoidal oscillations associated with the
linear wave equation
\[
  \varrho_{tt}-c^2\Delta \varrho =0, \qquad c^2 = p'(\ol\rho),
\]
which is obtained by linearizing \eqref{euler} around a constant state
with density $\rho=\ol\rho$ and velocity $u=0$.  In this paper we
prove that there is no fundamental inconsistency between the linear
and nonlinear theories of sound, by proving the existence and
genericity of pure modes of oscillation when the entropy is
non-constant.  We prove that when conservation of energy is accounted
for, and entropy varies generically, \emph{all} pure modes of the
linearized wave equation perturb to nearby shock-free ``pure tone''
oscillatory solutions of the nonlinear equations.  That is, for
generic entropy profiles, we prove that the fundamental oscillations
of the \emph{nonlinear equations} are well approximated by the
corresponding linearized sound wave solutions.

Stokes proposed the use of surfaces across which solutions jump
discontinuously, now familiar as shock waves.  Because solutions
containing shocks are discontinuous, they cannot be understood as
classical solutions of the PDEs.  The fundamental ``difficulty''
referred to by Stokes, then, is that when the pressure depends only on
the density $p=p(\rho)$, nonlinear oscillations evolve away from the
linear oscillations, ultimately forming shock waves, which are far
from linear.  This was later clarified by Riemann~\cite{Riemann}, who
proved that \emph{all} compressive plane wave (that is
one-dimensional) solutions of \eqref{euler} must \emph{necessarily}
form shocks in finite time, provided that $p''(\rho)\ne0$.  This was
made definitive in the celebrated work of Glimm and Lax in 1970, in
which they proved that all space periodic solutions of isentropic or
barotropic genuinely nonlinear systems necessarily form shocks and
decay to average at rate $1/t$~\cite{GL}.  At the time, the prevailing
sentiment in the field was that this result would extend to
non-isentropic solutions as well.

Since then there have been several attempts to resolve the question of
shock formation in oscillatory solutions of \x3 Euler, including
refinements of Riemann's shock formation result by the analysis of
characteristics~\cite{Lax64, J, TPLdec, TPLblowup, Chen1, CPZ1, CY1},
construction of periodic solutions in approximations based on weakly
nonlinear geometric optics~\cite{MRS, HMR, P}, and the theoretical and
analytic study of simplified models~\cite{Yex, Yper, Ysus}.  Direct
numerical simulations suggested the existence
of an attractor consisting of nonconstant long time limits of
spatially periodic solutions \cite{Shef, Vayn}.

The steepening of gradients and formation of discontinuities in
solutions of the nonlinear equations, which are presumed to model the
continuum more exactly, diverge both qualitatively and quantitatively
from the sinusoidal oscillations of the \emph{linearized} wave
equation, the classical model for sound waves.  This raises the
apparent paradox that the linearized equations do not provide a good
model for the nonlinear equations.  On the one hand, the nonlinear
equations should provide a more accurate description of the physics of
sound, while on the other hand, the linear wave equation, which in
principle should only be an approximation, provides an apparently more
accurate and time-tested way of representing actual sound waves.
Indeed, only a small nonlinear steepening of a sinusoidal sound wave
results in an audibly noticeable degradation of sound quality long
before shock formation, which is at odds with everyday
experience~\cite{Y-sound-demo}.  On the other hand, compressions
\emph{necessarily} become shocks in the one-dimensional isentropic
system, even for very small amplitudes.  This conflict raises the
important question:

\medskip
\begin{center}
  \begin{minipage}[h]{0.8\linewidth}
    \emph{In what way can the use of the linearized equations be
    mathematically justified when the linear and nonlinear equations
    produce phenomenologically different solutions?}
  \end{minipage}
\end{center}
\par
\medskip

Our results are summarized in the following theorem, which is
stated more precisely in Theorems \ref{thm:linear}, \ref{thm:nl} and
\ref{thm:top}, below.  We consider the problem of perturbing a
\emph{``quiet state''} to a solution of the fully nonlinear
compressible Euler equations.  A \emph{quiet state} is defined to be a
stationary solution of the form
\begin{equation}
  \label{quiet}
  p(x,t) = \ol p, \qquad
  u(x,t) = 0, \qquad
  s(x,t) = s(x), \qquad 0\le x \le \ell,
\end{equation}
where the specific entropy $s(\cdot)$ is piecewise $C^1$, either with
or without discontinuities.  The density is \emph{not} constant in a
quiet state but is determined by the equations through the consitutive
relation $\rho = \rho(p,s)$, and varies spatially with the entropy.
Nevertheless, quiet states, and the perturbed solutions constructed
here with discontinuous entropy, are classical shock-free solutions in
the sense that $u$ and $p$ remain continuous and the equations are
satisfied almost everywhere.  This class includes time-reversible
contact discontinuities, across which $s$ and $\rho$ jump, while $u$
and $p$ remain continuous.  Shock waves, across which $u$ and $p$ are
discontinuous, are not time-reversible, and hence are inconsistent
with time-periodic evolution.

We consider equations of state which are \emph{separable}, in that
the density $\rho = \rho(p,s)$ has the form
\begin{equation}
  \label{sepEoS}
  \frac{\del\rho}{\del p} = \rho^2\,a(p)^2\,\vs(s)^2,
  \com{with}
  \frac{da}{dp} \ne 0, \quad \frac{d\vs}{ds} \ne 0,
\end{equation}
uniformly in $(p,s)$.  In particular, this class includes all
ideal polytropic $\gamma$-law gases, $\gamma>1$, for which
\[
  p = \rho^\gamma\,e^{s/c_v}, \com{or}
  \rho = p^{1/\gamma}\,e^{s/c_p},
\]
but excludes the isothermal case $\gamma=1$.

The starting point of our theory is to identify a self-adjoint
boundary condition in \eqref{gpbc} below, which generalizes the
classical acoustic boundary condition $u=0$ at $x=0,\ \ell$.  The main
purpose of this paper is then to prove that for each quiet state
\eqref{quiet}, nonresonant pure modes of the linearized equations
\eqref{lin}, \eqref{gpbc}, perturb to nonlinear pure tone solutions of
the compressible Euler equations, according to the following theorem.

\begin{theorem}
  \label{thm:main}
  Solutions which satisfy the self-adjoint boundary condition
  \eqref{gpbc} extend, by a reflection principle, to global space and
  time periodic solutions.  This holds for the equations obtained by
  linearizing the compressible Euler equations about any quiet state,
  as well as the fully nonlinear equations themselves, both expressed
  in material coordinates.  Moreover,
  \begin{enumerate}
  \item Pure-tone $k$-mode solutions of the boundary value problem for the
    linearized equations exist, by Sturm-Liouville theory, for a
    sequence of distinct frequencies $\omega_k$, $k=1,2,3,...$, which
    depend on the quiet state.  General solutions of the linearized
    equations can be expressed as infinite sums of pure-tone solutions.
  \item For separable equation of state, each \emph{non-resonant} pure
    $k$-mode solution of the linearized equations perturbs to a one
    parameter family of ``pure tone'' $k$-mode solutions of the
    nonlinear compressible Euler equations, meeting the same
    self-adjoint boundary conditions at $x=0$, $x=\ell$.  A $k$-mode
    is non-resonant if $\omega_k$ is not a rational multiple of any
    other $\omega_j,$ $j\neq k$, of the linearized problem.
    Non-resonance forces the small divisors inherent in the problem to
    be non-zero.
  \item Entropy profiles are generically ``fully non-resonant''.  In a
    fully nonresonant profile, \emph{all} $k$-modes are non-resonant
    and perturb to nonlinear solutions.
  \end{enumerate}
  That is, for generic quiet states, every $k$-mode solution of the
  linearized equations perturbs to a one parameter family of ``pure
  tone'' nonlinear $k$-mode solutions of the compressible Euler
  equations satisfying boundary conditions \eqref{gpbc}, and thus
  extend to shock-free space and time periodic solutions of
  compressible Euler.
\end{theorem}

For the proof we identify a restricted class of symmetries which are
met by both linear and nonlinear solutions of the boundary value
problem.  Assuming these symmetries {\it a priori}, we deduce the
existence of a new nonlinear functional which imposes our boundary
conditions by projection onto symmetry, and at the same time has the
property that the small divisors inherent in the problem factor
uniformly out of the operator when the equation of state is separable.
This produces a nonlinear functional to which the Implicit Function
Theorem can be directly applied when a pure mode is
\emph{non-resonant}, a condition guaranteeing the small divisors are
non-zero.  Remarkably, our method does not require decay estimates on
the small divisors, which are \emph{not} uniformly non-zero with
respect to the entropy profile~\cite{TYdiff2}.  By this we overcome
difficult technical obstacles such as diophantine estimates, and the
expunging of resonances, which are inherent in KAM and Nash-Moser
problems, and which have been required for constructing periodic
solutions in physically simpler
settings~\cite{CW,BW,Rabinowitz,Moser}.

It has long been an open problem as to whether space and time periodic
solutions of compressible Euler exist.  The existence of such
solutions establishes definitively that shock formation can be avoided
in globally bounded classical solutions which exhibit sustained
nonlinear interactions, by which we mean solutions which \emph{do not}
decay time-asymptotically to non-interacting wave patterns.  Our
results also provide the first global existence result for classical
solutions of the acoustic boundary value problem.

Our theorem shows that for nonresonant entropy profiles and separable
equations of state (EoS), the oscillatory solutions of the linearized
system do indeed approximate pure tone solutions of the nonlinear
compressible Euler equations, so that when energy conservation is
taken into account, the linear acoustics equations are in fact a good
approximation of the nonlinear system.

\begin{corollary}
  \label{cor:acoustics}
  For separable EoS and nonresonant entropy profiles, which are
  generic, all pure mode solutions of the linearized wave equation are
  the small amplitude limits of (scaled) solutions of the nonlinear
  compressible Euler system; that is, for the solutions $(p,u)$ of
  \eqref{NLsol}, the limits
  \[
    \begin{aligned}
      P(x,t)
      &:= \lim_{\a\to0}\frac{p(x,t)-\ol p}{\a} = \cos(k\,t)\,\vp_k(x),\\
      U(x,t)
      &:= \lim_{\a\to0}\frac{u(x,t)}{\a} = \sin(k\,t)\,\psi_k(x),
    \end{aligned}
  \]
  exist globally in $x$ and $t$, and satisfy the linear wave equation
  \eqref{lin}, or equivalently \eqref{linwv}.  Here $\vp_k$ and
  $\psi_k$ are solutions of the associated Sturm-Liouville problem.
\end{corollary}

Most interesting to the authors, these solutions point to a new
physical mechanism that attenuates shock wave formation in oscillatory
solutions of compressible Euler when the entropy is non-constant.
Namely, non-constant entropy breaks the coincidence of wave speeds
between characteristic and phase velocities in oscillatory solutions.
The physical mechanism for this is repeated nonlinear interactions, or
``echoes'', with the background entropy profile.  That is, the
reflection of waves by the entropy gradient deflects the speed of
characteristics away from the speed of constant Riemann invariants,
setting up a consistent periodic pattern which cycles characteristics
through balancing regions of compression and rarefaction.  The effect
is that when entropy is non-constant, linear and nonlinear
characteristics generically move ergodically through the periods.  By
this, nonlinear solutions explore balancing regions of rarefaction and
compression, and this has the effect of strongly attenuating shock
wave formation in oscillatory solutions when the signals are not too
strong.  Our work here demonstrates that this phenomenon prevents
shock formation \emph{entirely} in pure mode oscillations.  The
authors believe that this is an important effect in the physical
propagation of nonlinear signals, which to our knowledge has not
previously been definitively expressed or understood.

The results of this paper are the culmination of authors' earlier
ideas, developed in~\cite{TYperStr, TYperLin,TYperBif,TYperG,TYperEV,
  TYperNM, TYdiff1, TYdiff2}.  The results provide a new class of
solutions of the nonlinear system \eqref{1Dsys} and show that
\emph{Riemann's general shock formation result holds only for the
  isentropic system~\eqref{psys}, but not for the entropic
  system~\eqref{1Dsys}}.

\subsection{Statement of Results}

The Euler equations \eqref{euler} represent conservation of mass and
momentum, respectively.  It was realized after the work of Stokes and
Riemann, that they should be augmented with the conservation of
energy, expressed as
\begin{equation}
  \label{energy}
  \del_t\big(\TS\frac12\rho\,u^2 + \rho\,e\big)
  + \nabla\cdot\big(\TS\frac12\rho\,u^3 + \rho\,u\,e + u\,p\big) = 0,
\end{equation}
while the system should be closed using the Second Law of
Thermodynamics,
\begin{equation}
  \label{TD}
  de = \theta\,ds-p\,dv.
\end{equation}
Here the variables are velocity $u$, density $\rho$, pressure $p$,
specific internal energy $e$, temperature $\theta$, specific volume
$v:=1/\rho$ and specific entropy $s$, and using \eqref{TD}, any two of
the thermodynamic variables determine all others, and we use $p$ and
$s$, so the system is closed by taking
\begin{equation}
  \label{evps}
  e=e(p,s) \com{and} v=v(p,s)
\end{equation}
as the constitutive relations.

For \emph{classical} solutions, the energy equation \eqref{energy} is
equivalent to the conservation of entropy,
\begin{equation}
  \label{entropy}
    \del_t\big(\rho\,s\big) + \nabla\cdot\big(\rho\,s\,u\big) = 0,
\end{equation}
which is structurally simpler, especially when $s$ is used as one of
the thermodynamic variables~\cite{CF,S}.  On the other hand, if
discontinuities are present in the solution, \eqref{entropy} is to be
interpreted as the \emph{entropy inequality}, and this serves to
select admissible discontinuities, which are shock waves.  In this
paper, the solutions we obtain are classical in the sense that $p$ and
$u$ are everywhere continuous, with possibly discontinuous $\rho$ and
$s$, and equations \eqref{euler}, \eqref{entropy} are satisfied almost
everywhere.

For us it is preferable to work in the \emph{Lagrangian frame}, which
refers to the \emph{material} coordinate ${\bf x}$, in which the (full)
compressible Euler equations take on the simpler form
\begin{equation}
  \label{ce}
  \begin{gathered}
  \del_t v - \nabla_{\bf x}\cdot u = 0, \\
  \del_t u + \nabla_{\bf x}\; p = 0, \\
  \del_t\big(\TS\frac12 u^2+e\big) + \nabla_{\bf x}\cdot\big(u\,p\big) = 0,
  \end{gathered}
\end{equation}
again augmented by \eqref{TD}, and \eqref{evps} closes the system.
The two systems \eqref{euler}, \eqref{energy} and \eqref{ce} are
equivalent for weak and classical solutions~\cite{CF,Wag}, and for
classical solutions, the energy equation is equivalent to the
(Lagrangian) entropy equation
\begin{equation}
  \label{slagr}
  \del_t s = 0, \com{so that} s({\bf x},t) = s^0({\bf x}),
\end{equation}
so that the entropy can be regarded as a stationary wave in the material
variable ${\bf x}$.

It follows that for \emph{classical solutions without shocks}, the
full system of the compressible Euler equations can be written as
\begin{equation}
  \label{sys}
  \del_t v - \nabla_{\bf x}\cdot u = 0, \qquad
  \del_t u + \nabla_{\bf x}\;p = 0,
\end{equation}
in which the constitutive equation can be taken to be either
\begin{equation}
  \label{stand}
  v = v\big(p,s^0({\bf x})\big) \com{or}
  p = p\big(v,s^0({\bf x})\big),
\end{equation}
and this in turn can be expressed as a single \emph{nonlinear wave
  equation}, namely
\begin{equation}
  \label{wveq}
  \del_t^2 v + \Delta_{\bf x}\,p = 0,
\end{equation}
together with \eqref{stand}.  We make the standard assumption that
\[
  \frac{\del p}{\del v}<0, \com{so also}
  \frac{\del v}{\del p}<0,
\]
which implies that the system is strictly hyperbolic~\cite{CF}, as
well as the separability condition \eqref{sepEoS}.  It is now evident
that linearization of \eqref{wveq} around the constant pressure
$\ol p$ yields the linear wave equation
\begin{equation}
  \label{linwv}
  \del_t^2 P - c^2\,\Delta_{\bf x}\,P = 0,
\end{equation}
but now the wavespeed $c = 1/\sqrt{-\frac{\del v}{\del p}}$ depends
also on the entropy $s = s^0({\bf x})$.

We now restrict to one space dimension, so we are studying plane waves
in direction ${\bf k}$, say, with one-dimensional variable
$x:={\bf k}\cdot {\bf x}$.  In this case, provided no shocks are
present, the entropy can be regarded as a standing wave, $s=s(x)$, and
\eqref{sys} becomes the \x2 system
\begin{equation}
  \label{1Dsys}
  \del_x p + \del_t u = 0. \qquad
  \del_x u - \del_t v = 0, \qquad
  v = v\big(p,s(x)\big).
\end{equation}
The Lagrangian formulation was known to Euler, and the mapping between
spatial (Eulerian) coordinate and material (Lagrangian) coordinate
${\bf x}$ is solution dependent, and in one dimension reduces to the
map $X\to x$, given by
\[
x(X,t) := x_0(t) + \int_0^X \rho(\chi,t)\;d\chi,
\]
where $x_0(t)$ is a given particle path\cite{CF,LaxBook,S}.

We treat the system as evolving in material coordinate $x$, with
inhomogeneous nonlinear term $v\big(p,s(x)\big)$, in which the entropy
profile is regarded as fixed and known.  If the entropy is constant
(the isentropic case) or if specific volume $v$ is a function of
pressure alone (the barotropic case), system \eqref{1Dsys} reduces to
the well-known $p$-system,
\begin{equation}
  \label{psys}
  \del_x p + \del_t u = 0. \qquad
  \del_x u - \del_t v(p) = 0,
\end{equation}
and this is the system for which Riemann showed that \emph{all
  compressions eventually lead to shocks}~\cite{Riemann,Lax64,S}, and
Glimm and Lax showed that periodic solutions decay to average by shock
wave dissipation at rate $1/t$~\cite{GL}.

Because sound waves are pressure variations, we treat solutions as
perturbations of constant ``quiet states'' $(\ol p,\ol u)$.  Without
loss of generality, a Galilean transformation in the spatial variable
allows us to take $\ol u=0$.  For a given entropy profile, a quiet
state is fully described by a single constant $\ol p$, while the
specific volume $v\big(\ol p,s(x)\big)$ varies with the entropy.  The
linearization of \eqref{1Dsys} around the constant quiet state
$(\ol p,0)$, with stationary entropy profile $s(x)$, is the
inhomogeneous linear system
\begin{equation}
  \label{lin}
  \del_x P + \del_t U = 0. \qquad
  \del_x U + \sg^2\,\del_t\;P = 0,
\end{equation}
with varying linear wavespeed $\sg = 1/c(\ol p,s)$ given by
\begin{equation}
  \label{sigma}
  \sg = \sg(x) := \sqrt{-\frac{\del v}{\del p}\big(\ol p,s(x)\big)}.
\end{equation}
Because the inhomogeneous linearized wave speed is not time dependent
in a Lagrangian frame, one of the fundamental simplifications in our
point view here is to evolve in the material coordinate $x$, rather
than the time $t$, as is standard in most treatments.  This has the
effect of putting the inhomogeneity $s(x)$ into a Sturm-Liouville
framework.

Our procedure is to first solve the linear system \eqref{lin} by
separation of variables, and then prove that generically, all pure
mode solutions of \eqref{lin} perturb to time and space periodic
solutions of the fully nonlinear system~\eqref{1Dsys}.  We accomplish
this by solving a \emph{periodic boundary value problem} (PBVP) on the
region $(x,t)\in[0,\ell]\times[0,T]$, which we refer to as a
\emph{tile}.

To begin, we introduce the \emph{acoustic boundary value problem}
(ABVP) for time-periodic solutions, by imposing a reflection
boundary condition, namely
\begin{equation}
 \label{abvp}
  \begin{gathered}
    u(x,0) = u(x,T), \qquad
    p(x,0) = p(x,T), \qquad
    x\in[0,\ell],\\
    u(0,t) = 0, \qquad
    u(\ell,t) = 0, \qquad
    t\in[0,T].
  \end{gathered}
\end{equation}
This models pure reflections off a solid wall in acoustics.  Note that
we can impose time-periodicity because the inhomogeneity is isolated
in the spatial variations of $s(x)$.

We will prove the existence of pure tone solutions of this ABVP, for
which $p$ is even and $u$ is odd as functions of $x$.  These in turn
generate space and time periodic solutions with time period $T$ and
spatial period $2\ell$ by a reflection principle.  This accounts for
the periodic solutions for \emph{even} modes, but to generate the
\emph{odd} modes we must generalize the boundary condition.  By this
we find that the odd modes generate $4\ell$ spatially periodic
solutions from the same tile under the same even/odd conditions, when
$p$ and $u$ satisfy the more general \emph{periodic} boundary
conditions $u(0,t) = 0$, with
\begin{equation}
  \label{pbc}
  \IR-\,\mc T^{-T/4}\,p(\ell,t) = 
  \IR-\,\mc T^{-T/4}\,u(\ell,t) = 0,
\end{equation}
for $t\in[0,T]$.  Here $\mc R$ and $\mc T$ are linear reflection and
translation operators, defined respectively by
\[
  \mc R[f](t) := f(-t), \com{and}
  \mc T^\tau[f](t) := f(t-\tau).
\]
The presence of the term $\IR-$, which is a projection onto odd modes,
only imposes conditions on half of the modes, which is why the two
conditions in \eqref{pbc} are consistent.  Moreover, it is convenient
to now introduce a new variable $y$ and write \eqref{pbc} as the single
condition
\begin{equation}
  \label{ydef}
  \IR-\,\mc T^{-T/4}\,y(\ell,t) = 0, \com{where}
  y := p+u.
\end{equation}
The quantity $y$ is dimensionally inconsistent, because it adds $p$
and $u$.  However, since $p$ is even and $u$ is odd, we can formally
add these orthogonal pieces, and we use $y$ as a convenient notational
device to track the separate evolutions of $p$ and $u$ simultaneously.
We then recover $p$ and $u$ and their derivatives from $y$ and its
derivatives by simply projecting $y$ onto even or odd modes as
appropriate.  Formally, $y$ solves a nonlocal nonlinear scalar balance
law, which in turn simplifies the analysis of solutions by the method
of characteristics, because it reduces to a single characteristic
family.

To incorporate all of the above conditions into a single problem for
the generating tile, it suffices to consider the following general
boundary conditions:
\begin{equation}
  \label{gpbc}
  \begin{aligned}
    u(x,0) = u(x,T), \qquad
    p(x,0) = p(x,T), \qquad
    &x\in[0,\ell],\\
    u(0,t) = 0,\qquad
    \IR-\,\mc T^{-T/4}\big(p(\ell,t)+u(\ell,t)\big) = 0,
    \qquad
    &t\in[0,T].
  \end{aligned}
\end{equation}
These boundary conditions are consistent with and suggested by the
\emph{critical observation} that a combination of reflections and
translations extend the solution to a tiling of the whole plane in
which both $p(x,t)$ and $u(x,t)$ are jointly periodic, namely
\[
  (p,u)(x+4\,\ell,t) = (p,u)(x,t) \com{and}
  (p,u)(x,t+T) = (p,u)(x,t),
\]
for all $x$ and $t$, for both the linear and nonlinear systems.
Moreover, these boundary conditions imply that the Sturm-Liouville
operator obtained by separating variables in \eqref{lin} is
self-adjoint.

The acoustic boundary condition \eqref{abvp} is physically realizable
for arbitrary entropy profiles on $[0,\ell]$.  On the other hand, the
general boundary condition \eqref{gpbc} together with our reflections
imply that the entropy is $2\ell$ periodic while the solution is
$4\ell$ periodic in $x$.  Periodic entropy profiles in all of space
should thus be regarded as finely tuned, but the boundary value
problem that generates them is robust in both the linear and nonlinear
problems.  To the authors' knowledge, this is the first instance of a
proof of the existence of globally defined, bounded, \emph{classical}
solutions of the acoustic boundary value problem \eqref{abvp} for the
nonlinear system of compressible Euler equations.  Although Nishida
and Smoller obtained globally bounded solutions in $BV$ for a general
boundary value problem in the special case of an \emph{isothermal}
equation of state~\cite{NS2}, any non-constant solutions \emph{must}
contain shock waves and necessarily decay to a constant state because
the equation of state is barotropic.
 
For our purposes it suffices to restrict to the class $\mc P$ of
piecewise $H^b$ entropy profiles,
\begin{equation}
  \label{P}
  \mc P := \Big\{s:[0,\ell]\to \B R\;\Big|\;
  s(x) \text{ piecewise }H^b\;\Big\},
\end{equation}
for Sobolev exponent $b>3/2$, which includes both $C^1$ and piecewise
constant profiles, and is dense in $L^1$.  It follows that $\sg$
defined by \eqref{sigma} is also piecewise $H^b$, $\sg\in\mc P$.  We
adopt the standard definition of piecewise $H^b$ as meaning a bounded
function with finitely many jumps in the interior of the interval
$[0,\ell]$.  For brevity we abbreviate the trigonometric functions
$\c:=\cos$ and $\s:=\sin$.  Within this framework, we state our main
results in the following three theorems which together precisely
restate Theorem~\ref{thm:main}:

\begin{theorem}
  \label{thm:linear}
  For every entropy profile $s\in\mc P$, and every~$k\in\B N$, the
  linearized system \eqref{lin}, \eqref{gpbc} admits a $k$-mode
  solution
  \begin{equation}
    \label{linsol}
    P(x,t) = \c(\w_k\,t)\,\vp_k(x), \qquad
    U(x,t) = \s(\w_k\,t)\,\psi_k(x),
  \end{equation}
  where $(\vp_k,\psi_k)$ is a solution of the Sturm-Liouville
  system~\eqref{SL} with corresponding \emph{eigenfrequency} $\w_k$.
  When combined with the boundary conditions \eqref{gpbc}, this extends
  to a jointly periodic pure ``$k$-mode'' solution with periods $4\,\ell$
  in $x$ and time period $T = 2\,\pi\,k/\omega_k$ of the linearized
  system \eqref{lin}.
\end{theorem}

For a given entropy profile, we define the $k$-mode to be
\emph{nonresonant} provided that the
equation
\[
  k\,\w_l - j\,\w_k = 0, \qquad l,\ j\in\B N,
\]
has no solutions.  We then define the \emph{set of
  resonant profiles} by
\[
  \mc Z := \Big\{ s(\cdot)\in\mc P\;\Big|\;  k\,\w_l = j\,\w_k
  \text{ for some } (k,j,l)\in\B N^3\Big\},
\]
so that any $s\in \mc Z$ has at least one resonant mode.  In
particular, a sufficient condition for nonresonance is that the ratios
$\w_l/\w_k$ are irrational for each $l\ne k$.  Resonance of a mode is
fully determined at the level of the linearized equations.  We note
that in the linear homogeneous wave equation, obtained when the
entropy is constant, the frequencies all take the form
$\omega_k = k\,\ol\omega$, for some constant $\ol\omega$.  In this
case, all frequencies are determined by the single number $\ol\omega$,
and in particular constant entropy profiles are \emph{completely
  resonant}.  In the next theorem, we show that any nonresonant
$k$-mode solution of Theorem~\ref{thm:linear} perturbs to a
one-parameter family of solutions of the nonlinear compressible Euler
equations.

\begin{theorem}
  \label{thm:nl}
  For a separable EoS satisfying \eqref{sepEoS}, let any entropy
  profile $s(\cdot)\in\mc P$ and Sobolev index $b>5/2$ be given.  If
  the $k$-mode is nonresonant, $\w_j/\w_k\notin\B Q$ for all $j\ne k$,
  then there is some $\ol{\a}_k>0$, such that the nonlinear system
  \eqref{1Dsys}, \eqref{gpbc} admits a 1-parameter family of
  continuous and piecewise $H^b$ solutions of the form
  \begin{equation}
  \label{NLsol}
    \begin{aligned}
      p(x,t) &= \ol p + \a\,\c(\w_k\,t)\,\vp_k(x) + O(\a^2),\\
      u(x,t) &= \a\,\s(\w_k\,t)\,\psi_k(x) + O(\a^2),
    \end{aligned}
  \end{equation}
  for $|\a|<\ol{\a}_k$.  These in turn extend by reflection to jointly
  periodic \emph{pure tone} solutions of the nonlinear compressible
  Euler equations with periods $4\,\ell$ in $x$, and time period
  $T = 2\,\pi\,k/\omega_k$.
\end{theorem}

Finally, for the third theorem, we regard $\mc P$ as a dense subspace
of $L^1(0,\ell)$.  We give a sense in which the \emph{fully
  nonresonant} entropy profiles are generic.

\begin{theorem}
  \label{thm:top}
  Given any piecewise constant approximation to an $L^1$ entropy
  profile, the probability of that approximation having any resonant
  linearized modes is zero, within the class of all piecewise constant
  profiles with the same number of jumps.  In this sense, for
  \emph{generic} piecewise $H^b$ entropy profiles, \emph{every}
  linearized $k$-mode perturbs to a one-parameter family of nonlinear
  pure tone solutions of the form \eqref{NLsol} of the compressible
  Euler equations.
\end{theorem}

Note that the perturbation Theorem \ref{thm:nl} is the only one that
requires the separability condition \eqref{sepEoS}, the other theorems
remaining true for general equations of state. 

\subsection{Discussion of Results}

The linearized profiles $(\vp_k,\psi_k)$ in \eqref{NLsol} are the
eigenfunctions, corresponding to eigenvalue $\w_k^2$, of the
Sturm-Liouville (SL) system
\begin{equation}
  \label{SL}
  \dot\vp_k + \w_k\,\psi_k = 0, \qquad
  \dot\psi_k - \sg^2\,\w_k\,\vp_k = 0,
\end{equation}
which is obtained by separating variables using the ansatz
\eqref{linsol} in the linearization \eqref{lin}.  It is known that the
SL system \eqref{SL} has continuous solutions provided $\sg$ is
piecewise continuous, as long as continuity of the solution
$(\vp_k,\psi_k)$ is imposed at any jump in the SL coefficient
$\sg$~\cite{Pryce,Zettl,BR}.  In the linear problem continuity of
$(\vp_k,\psi_k)$ is equivalent to continuity of $p$ and $u$.  In the
nonlinear problem, imposing continuity in $p$ and $u$ forces a
time-reversible contact discontinuity in $\rho=1/v$~\cite{S,CF}.  The
linearized evolution through a piecewise $H^b$ entropy field can then
be realized as a series of non-commuting compositions of evolutions on
$H^b$ fields.  The concatenation of these smooth evolutions together
with continuity of $p$ and $u$ at entropy jumps yields a bounded
evolution operator in Sobolev spaces $H^b$ for arbitrary $b>5/2$.  In
particular, if the entropy profile is piecewise constant, the full
evolution is a non-commuting composition of \emph{isentropic}
evolutions and entropy jumps.

For the isentropic system, linearization around a constant state gives
the linear wave equation~\eqref{linwv}, so all frequencies have the
form $k\,\ol\w$.  In this case, all pairs of modes are resonant, so
that no modes can be perturbed by our methods.  Thus Corollary~1
\emph{fails} for isentropic profiles, consistent with the observations
of Stokes and Riemann.  Since the seminal work of Glimm and
Lax~\cite{GL}, experts have generally thought that shocks would always
form and solutions would decay to average in space periodic solutions
of the full compressible Euler equations.  In the context of
\cite{GL}, our results are thus surprising and provide the first
definitive proof that shock formation and decay to average
\emph{fails} when a full system of Riemann invariants is not
available.

In this paper we prove that the fundamental modes of the linear
theory, from which all solutions of \eqref{linwv} are obtained by
superposition, perturb to nonlinear solutions.  Although we do not
have a nonlinear superposition principle, and nonlinear solutions may
shock, stability of the fundamental modes strongly indicates that
generic solutions with or without shocks will not decay to average.
Thus, if there exists a Glimm-Lax decay theory for \x3 Euler, we
cannot expect decay to average as in the \x2 case.  Indeed, solutions
constructed here have the property that nonlinear rarefaction and
compression are in balance on average along every characteristic,
because characteristics move ergodically through the entire entropy
profile as a consequence of the nonresonance condition.  This
identifies a new mechanism that attenuates shock wave formation in
oscillatory solutions of compressible Euler when the entropy is
non-constant, and suggests that long-time limits are much richer than
constant states.  This is consistent with the numerical experiments
of~\cite{Shef, Vayn}, and with the predicted results of weakly
nonlinear geometric optics~\cite{MRS}.

Regarding stability, by continuity, perturbations of our nonlinear
periodic solutions will balance rarefaction and compression in the
limit of small perturbation.  This has the effect of sending the time
of shock formation to infinity as a perturbation tends to zero, an
aspect of stability.  A complete stability analysis would require
knowledge of all solutions which similarly balance rarefaction and
compression on average along characteristics and the effect of this
averaging process under perturbation.  Authors are exploring the
possibility that pure tone nonlinear solutions, like their linearized
pure mode counterparts, will perturb to nearby \emph{almost periodic}
classical solutions which balance compression and rarefaction
sufficiently to avoid shock formation.  If strong perturbations of
these solutions do form shocks, they will decay by entropy
dissipation, but we cannot expect the long-time limits of such
decaying solutions to be restricted to just ``quiet'' constant states;
rather, they would decay to any solution which balances rarefaction
and compression like the pure tone nonlinear solutions constructed
here.  The simulations of~\cite{Vayn, Shef} can be interpreted as a
numerical indication of stability.  Prior to the results of this
paper, the question of stability of such solutions could not be
formulated mathematically.

The authors began this project in~\cite{TYperStr, TYperLin, TYperG,
  TYperBif} by identifying the physical phenomenon of ``echoing''
which can prevent a single compression from forming a shock in finite
time by interaction with the entropy profile.  We conjectured that
space and time periodic solutions to Euler exist, and clarified the
difficult technical issues that arose~\cite{TYperNM, TYperEV},
including issues of small divisors and KAM-type
behavior~\cite{TYdiff1, TYdiff2}.  A valid Nash-Moser procedure was
set out in a scalar setting in \cite{TYdiff2}, but the difficult
technical issues described there were bypassed in this paper by the
discovery and use of symmetries.  With the results of this paper, our
initial conjectures, and much more, are confirmed.

In our previous treatments, we framed the problem by imposing
periodicity by ``periodic return'', namely
\begin{equation}
  \label{F1}
  \mc F_1(U) = \mc E^{L}(U) - U = 0.
\end{equation}
In this framework, and other Nash-Moser treatments, some expunging of
resonant parameters is required because there are no uniform bounds on
the decay rates of small divisors \cite{CW,BW}.

The key breakthrough in this paper is to identify fully, and
leverage, the \emph{symmetries} that are respected by \emph{both the
  linear and nonlinear} evolutions of the systems \eqref{lin} and
\eqref{1Dsys}.  These symmetries, identified only in material
coordinates, are $p$ even, $u$ odd, in both $x$ and $t$ separately.
That is, if the data is preserved by the transformations
\begin{equation}
  \label{symm}
  \begin{aligned}
    \big(p(x,t),u(x,t)\big) &\mapsto \big(p(-x,t),-u(-x,t)\big),\\
    \big(p(x,t),u(x,t)\big) &\mapsto \big(p(x,-t),-u(x,-t)\big),
\end{aligned}
\end{equation}
then solutions of systems \eqref{lin} and \eqref{1Dsys} also respect
\eqref{symm}.  This in turn means that the \emph{periodicity
  constraint} \eqref{F1} can be replaced by a \emph{projection
  constraint} through the boundary conditions \eqref{gpbc} on a smaller
tile, while simultaneously reducing the domain of the nonlinear
operators to symmetric data.  This has the effect of essentially
halving the (infinite) number of ``degrees of freedom'' while
retaining all essential nonlinear effects.

Utilizing the symmetries, we construct a new nonlinear functional
$\mc F$ which imposes periodicity by projection onto symmetry, such
that the nonlinear operator can be factored into two parts,
$\mc F = \mc A\,\mc N$.  The nonlinear piece $\mc N$ is regular and
has ``base linearization'' $D\mc N(\ol p) = \mc I$, so we can apply
the implicit function theorem directly to $\mc N$ without need of
Nash-Moser.  The prefactor $\mc A$ is a \emph{fixed}, \emph{diagonal},
\emph{constant} linear operator which is bounded, and whose inverse
$\mc A^{-1}$, although unbounded $H^b\to H^b$, exists and is closed
for each $\ol p$ so long as the entropy profile is nonresonant.  This
operator $\mc A$ encodes all the small divisors, and acts as a fixed,
constant preconditioner for the nonlinear problem, which can then be
inverted by the implicit function theorem.  To get uniform estimates,
it suffices to adjust the $H^b$ norm on the target space via scaling
by the fixed small divisors, and this has the effect of turning
$\mc A$ into an isometry.  We then invert $\mc A$ explicitly, so that
the nonlinear problem becomes essentially: find $y^0$ such that
$\mc N(y^0)=0$, with $D\mc N(\ol p)= \mc I$, where $\mc N$ is regarded
as a map $\mc N:H^b\to H^b$.  Once this factorization is in place, the
problem can be treated as a standard bifurcation problem.

In summary, the symmetries \eqref{symm} enabled us to discover a
nonlinear operator $\mc F$ which both encodes periodicity and
conveniently factors.  This allows us to bypass Nash-Moser and the
inherent problem of controlling the small divisors in ``nearby
linearized operators''.  To obtain periodic solutions it suffices to
invert the nonlinear operator $\mc N$ using the standard implicit
function theorem.  It is remarkable that for nonresonant profiles, we
get a fixed decay rate for small divisors of $\mc A=\mc A(\ol p)$ at
each $\ol p$, but this decay rate is not absolutely continuous in
$\ol p$.  In the bifurcation argument, it is essential that the
average value of $p$ be perturbed, because that is precisely our
bifurcation parameter.  In this way our perturbed nonlinear solutions
have a \emph{different} background pressure from the fixed base
pressure $\ol p$ which determines the small divisors.

The paper is laid out as follows.  In Section~\ref{sec:str}, we
introduce the symmetry and reflection principles which generate
periodic solutions, and derive the self-adjoint boundary conditions
\eqref{gpbc} and the nonlinear operator $\mc F$.  In
Section~\ref{sec:exist}, we provide a detailed description of the local
$H^b$ existence theorem that we use, including a uniform $L^\infty$
estimate for derivatives of solutions.  In Section~\ref{sec:lin}, we
linearize the problem and derive the corresponding Sturm-Liouville
system.  We then solve this SL system and develop the resonance
conditions.  In Section~\ref{sec:bif}, we develop the bifurcation
argument and prove existence of nonlinear pure tone solutions of the
boundary value problem, assuming technical lemmas which are proved in
later sections.  This involves rescaling of the target Hilbert space
by the small divisors, and produces periodic solutions using the
reflection principle.  In Section \ref{sec:SL} we carry out a detailed
analysis of the SL problem and estimate the growth rate of the
eigenfrequencies.  In Section \ref{sec:D2E} we show that the evolution
operator is twice differentiable and detail the role of genuine
nonlinearity.  In Section~\ref{sec:pwc}, we restrict to piecewise
constant entropy profiles, for which explicit formulae can be
developed, and use these to establish density of resonant profiles in
$L^1$, and genericity of non-resonant profiles.  Finally in
Section~\ref{sec:con}, we provide historical, scientific and
mathematical context for our results.

\section{Structure of Periodic Solutions}
\label{sec:str}

Our starting point is the one-dimensional compressible Euler equations
in a material frame \eqref{ce}, which is the $3\times3$ system in
conservative form,
\begin{equation}
  \label{lagr}
  v_t - u_x = 0, \quad u_t + p_x = 0, \quad
  \big(\TS{\frac12}\,u^2 + e\big)_t + (u\,p)_x = 0,
\end{equation}
see~\cite{CF,S}.  Here $x$ is the material coordinate and $u$ is the
Eulerian velocity, and the thermodynamic variables are specific volume
$v$, pressure $p$ and specific internal energy $e$.  The system is
closed by a \emph{constitutive relation} which satisfies the Second
Law of Thermodynamics,
\[
  de = \theta\,ds - p\,dv,
\]
so that
\[
  \theta = e_s(s,v) \com{and} p = -e_v(s,v),
\]
where $\theta$ is the temperature and $s$ the specific entropy.  It
follows that for classical smooth solutions, the third (energy)
equation can be replaced by the simpler \emph{entropy equation},
\[
  e_t + p\,v_t = 0,  \com{or}  s_t = 0,
\]
see~\cite{CF,S}.  Throughout this paper, we consider classical
solutions, so we regard the entropy $s(x)$ as a fixed, known,
stationary field.  Unlike solutions with shocks, solutions with
contact discontinuities, across which the pressure and velocity are
continuous, while entropy and density are discontinuous, can be
regarded as classical solutions because they are the uniform limit of
smooth solutions.  Indeed, in this paper, we explicitly incorporate
piecewise smooth entropy profiles.

This means that our general constitutive law becomes the inhomogeneous
condition
\[
  p = p\big(v,s(x)\big) \com{or} v = v\big(p,s(x)\big),
\]
and the entropy $s$ is no longer regarded as a dynamical variable.
In this point of view, we treat the pressure $p$ and velocity $u$ as
dependent variables, with general constitutive law $v = v(p,s)$.

For classical solutions, in a Lagrangian frame, we
write the system \eqref{lagr} as
\begin{equation}
  \label{genls}
  p_x + u_t = 0, \qquad
  u_x - v(p,s)_t = 0,
\end{equation}
or, in quasilinear form,
\[
  p_x + u_t = 0, \qquad
  u_x - v_p(p,s)\,p_t = 0,
\]
where we are again evolving in $x$, and since $s_t=0$, we regard the
entropy as a given fixed function $s(x)$.  Because we consider
\emph{classical} solutions, with possibly discontinuous entropy, we
\emph{require} that $p$ and $u$ be continuous throughout the
evolution.  This is consistent with the Rankine-Hugoniot jump
conditions for contact discontinuities, which are stationary in the
material frame~\cite{CF,S}.

\subsection{Symmetries and Boundary Conditions}

We construct periodic solutions within the restricted class of
solutions which admit the following essential symmetries respected by
the nonlinear system \eqref{genls}.  This represents a discovery more
than a loss of generality because the linearized sound wave solutions
which we perturb respect the same symmetries.

The first of the critical symmetries is \emph{even/odd symmetry in
time}, namely that the transformation
\begin{equation}
  \label{tsymm}
  p(x,t) \to p(x,-t), \quad u(x,t) \to -u(x,-t),
\end{equation}
is respected by the system, so that any data
$\big(p(x_0,\cdot),u(x_0,\cdot)\big)$ which satisfies \eqref{tsymm}
evolves in $x$ to a solution which also satisfies \eqref{tsymm}.

Solutions $(p,u)$ in this symmetry class can be described by the
equivalent scalar quantity \eqref{ydef}, namely
\[
  y(x,t) := p(x,t) + u(x,t),
\]
which satisfies the nonlinear, nonlocal, scalar  conservation law
\begin{equation}
  \label{ycl}
  y_x + g(y)_t = 0,
\end{equation}
in which the nonlocal scalar flux is given by
\begin{equation}
  \label{gy}
  g(y) := u + v(p,s) = \IR-y + v\big(\IR+y,s\big).
\end{equation}
Here $\mc R$ denotes the linear time reflection operator,
\begin{equation}
  [\mc R\,v](t) := v(-t),
  \label{refl}
\end{equation}
so that $\IR+$ and $\IR-$ are the projections onto even and odd parts
of a function, respectively.  Although the system \eqref{genls} and
scalar equation \eqref{ycl} are equivalent in the symmetry
class~\eqref{tsymm}, solutions of \eqref{ycl} are easier to analyze
because it has a single characteristic family.

The second symmetry condition is \emph{even/odd symmetry in material
  coordinate} $x$, namely that the transformation
\begin{equation}
  \label{xsymm}
  p(x,t) \to p(-x,t), \quad u(x,t) \to -u(-x,t),
\end{equation}
is also respected by system~\eqref{genls}, provided the entropy $s(x)$
is also even in $x$.  Because $x$ is the evolution variable, we impose
\eqref{xsymm} at the boundaries and use reflections to generate a
periodic tile.  Doing so imposes constraints on the data $y^0(\cdot)$
and on the solution $y(\ell,\cdot)$ which lead to our reformulation of
the nonlinear problem.

Recall that we are looking for continuous, classical solutions of the
system \eqref{genls}.  Imposition of the material symmetry
\eqref{xsymm} together with continuity requires $u(0,\cdot)=0$, or
equivalently $y(0,\cdot)$ even.  We thus \emph{restrict the domain} of
data $y^0$ for the scalar equation \eqref{ycl} to the ``half-space''
of even functions.  We evolve this even data forward in material space
to $x=\ell$.  This in turn yields the solution on the entire entropy
level $-\ell\le x\le\ell$ by reflection through
\begin{equation}
  \label{yrefl}
  y(-x,t) := \mc R\,y(x,t) = y(x,-t),
\end{equation}
which yields
\[
  \IR+y(-x,t) = \IR+y(x,t) \com{and}
  \IR-y(-x,t) = -\IR+y(x,t),
\]
which in turn is \eqref{xsymm}.  Note that our boundary condition is
the natural physical ``acoustic boundary reflection condition'' $u=0$
at $x=0$, and if this condition is imposed, the solution necessarily
satisfies \eqref{yrefl} by uniqueness.

It remains to impose a periodicity condition at the end of the
evolution.  For this we use the following principle: an even or odd
periodic function admits \emph{two} axes of symmetry, at both the
endpoint $0$ and midpoint $T/2$ of the interval, so that even or odd
periodic functions remain even or odd after a half-period shift,
respectively.  That is, if we set $\wt f(s):=f(s-T/2)$, then
\begin{equation}
  \label{reflper}
  \wt f(-s) = f(-s-T/2) = \pm f(s+T/2) = \pm f(s-T/2) = \pm\wt f(s).
\end{equation}
For periodicity, it thus suffices to impose the additional
\emph{shifted} reflection symmetry at the end of the evolution
$x=\ell$, which becomes the shifted symmetry axis.  This shifted
reflection is
\begin{equation}
  \begin{aligned}
    p(\ell+x',t) &= p(\ell-x',t-T/2), \\
    u(\ell+x',t) &= -u(\ell-x',t-T/2),
  \end{aligned}
\label{ulsymm}
\end{equation}
and setting $x'=0$, continuity of the solution at $\ell$ requires
\begin{equation}
  p(\ell,t-T/2)=p(\ell,t), \qquad
  u(\ell,t-T/2)=-u(\ell,t).
  \label{odddat}
\end{equation}
In other words, if $p(x,\cdot)$ and $u(x,\cdot)$ are defined for say
$0<x\le\ell$, and satisfy the symmetry condition \eqref{odddat} at
$x=\ell$, then use of \eqref{ulsymm} allows us to extend the solution
to the interval $0<x<2\ell$.  Assuming these symmetries, it follows
that we can generate a $4\,\ell$-periodic solution from a single tile
defined over $0\le x\le\ell$ via the reflections \eqref{yrefl} around
$x=0$ and \eqref{ulsymm} around $x=\ell$, respectively..

In order to fully express our nonlinear equation as a nonlinear
functional in the scalar field $y$, we must write condition
\eqref{odddat} in terms of $y$.  This is easiest done by introducing
the constant shift operator
\begin{equation}
  \label{shdef}
  \mc T^\tau[f](t) := f(t-\tau).
\end{equation}
Using a \emph{quarter-period} shift, we write \eqref{odddat} as
\[
  \big(\mc T^{T/4}\big)^2\, p(\ell,\cdot) = p(\ell,\cdot), \quad
  \big(\mc T^{T/4}\big)^2\, u(\ell,\cdot) = -u(\ell,\cdot).
\]
Since $p$ is even in $t$, $u$ is odd in $t$, and $y=p+u$, this is
\[
  \big(\mc T^{T/4}\big)^2\, y(\ell,\cdot) = \mc R\,y(\ell,\cdot),
\]
and since $\mc T^{-T/4}\,\mc R = \mc R\,\mc T^{T/4}$, this results
in the \emph{shifted reflection} condition
\begin{equation}
  \label{ySR}
  \IR-\,\mc T^{T/4}\,y(\ell,\cdot) = 0.
\end{equation}
We thus impose \eqref{ySR} as the nonlinear tiling condition that we
must solve: any even data $y^0$ whose evolution through $x=\ell$
satisfies \eqref{ySR} generates a periodic tile.

\begin{figure}[htb]
\includegraphics[width=\linewidth]{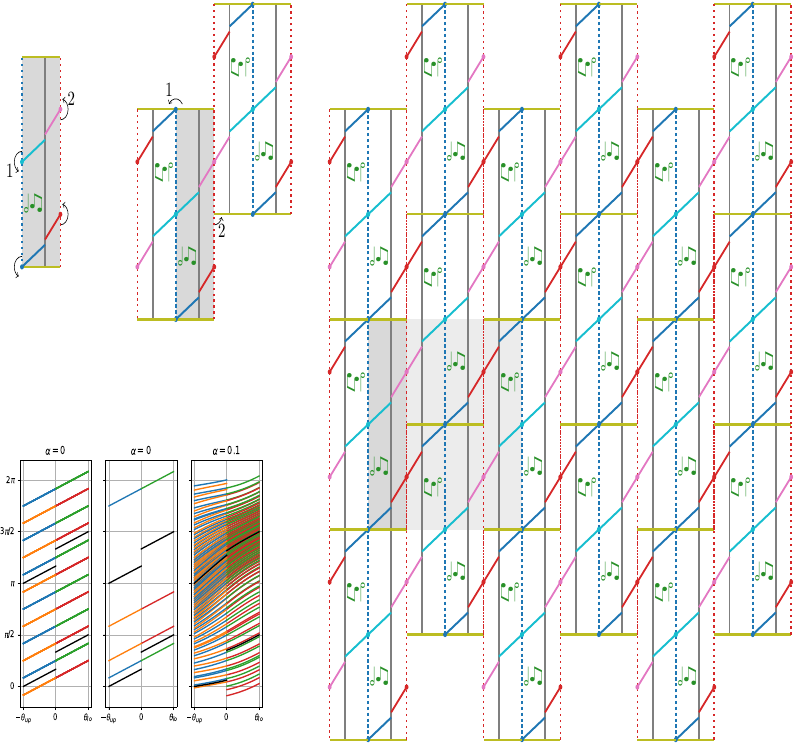}
\caption{Reflections of the generating tile.}
\label{fig:period}
\end{figure}

We illustrate the repeated reflections of the generating tile in
Figure \ref{fig:period}.  Recalling that two orthogonal reflections of
the plane yields a rotation of $\pi$ around the intersection point,
the standard reflection at $x=0$ is labelled `1', and the shifted
reflection at $x=\ell$ is labeled `2'.  The composition of these two
reflections is then translated periodically by $4\ell$ and $2\pi$ to
give the wave pattern on the plane.  A sketch of the perturbed tile,
showing the perturbed characteristics to leading order, is also shown.
The density of sketched characteristics indicates the presence of
nonlinear rarefaction and compression.  Only the forward
characteristics are shown; the backward characteristics are obtained
by applying the time reflection $\mc R$.

\subsection{The Reduced Nonlinear Problem}

We now precisely set up the nonlinear operator whose solutions provide
a tile which generates periodic solutions of compressible Euler.  To
do this we express the scalar evolution in $x$ by the nonlocal
conservation law \eqref{ycl} as the operator $\mc E^x$, defined by
\[
  \mc E^x(y^0) := y(x,\cdot), \com{for} x\in[0,\ell].
\]
By classical energy estimates, we have that $\mc E^x:H^b\to H^b$ for
any $b>3/2$ while the entropy $s(x)$ remains $H^b$, so that
$\|\frac{ds}{dx}\|_\infty<\infty$.  At any point $x$ at which the
entropy $s(x)$ has a jump, that is at any contact discontinuity, we
impose the continuity condition
\[
  p(x+,\cdot) = p(x-,\cdot), \quad 
  u(x+,\cdot) = u(x-,\cdot),
\]
which is equivalent to
\begin{equation}
  \label{cont}
  y(x+,\cdot) = y(x-,\cdot),
\end{equation}
and which implies that the $H^b$ norm is exactly preserved at any
contact.  Because our entropy profile is assumed piecewise $H^b$, the
evolution over the full interval can be regarded as a non-commuting
composition of partial evolutions, each of which is over a subinterval
in which the entropy is $H^b$.  It follows that we can write
\[
  y(\ell,\cdot) = \mc E^\ell(y^0),
\]
where by classical energy estimates, $\mc E^\ell:H^b\to H^b$ for any
$b>3/2$~\cite{Majda,Taylor}; note that the larger $b$ is, the smoother
our nonlinear solutions are known to be.  Here $\mc E^\ell$ is a
non-commuting composition of smaller evolutions for which the entropy
is $H^b$.  We provide a detailed statement of the classical existence
theorem in Corollary~\ref{cor:le} below.

If we adjust the shifted periodic boundary condition \eqref{ySR} to the
simpler even condition
\[
  \IR-\,y(\ell,\cdot) = 0,
\]
we obtain a smaller class of problems which satisfy the physically
important acoustic reflection boundary condition \eqref{abvp}, namely
$u(\ell,\cdot)=0$, global solutions of which have been an open problem
since the time of Euler.

The reflection principles outlined above reduce the problem of doubly
periodic solutions of compressible Euler to that of producing a single
tile.  To generate such a tile, we pose either of the following two
problems:
\begin{itemize}
\item The \emph{periodic tile problem}, namely
  \begin{equation}
  \label{yper}
    \mc F_P(y^0) := \IR-\,\mc T^{-T/4}\,\mc E^\ell\,y^0 = 0;
  \end{equation}
\item The \emph{acoustic boundary value problem}, namely
  \begin{equation}
  \label{yab}
    \mc F_A(y^0) := \IR-\,\mc E^\ell\,y^0 = 0.
  \end{equation}
\end{itemize}
In both of these problems, as above, $\mc E^\ell$ denotes nonlinear
evolution through the varying entropy profile from $x=0$ to $x=\ell$,
and the data $y^0=y^0(t)$ is again assumed to be even and $T$-periodic.
We use the following notation which incorporates both problems
\eqref{yper} and \eqref{yab}.  Denote the ``boundary shift operator''
by
\begin{equation}
  \label{Shat}
  \mc S := \IR-\,\mc T^{-\chi\,T/4}, \qquad \chi = 1 \com{or} 0,
\end{equation}
for the periodic or acoustic boundary, respectively, where $\chi$ is
an indicator of which boundary condition is being used.  Again, the symmetry
imposed by choosing $y^0$ even encodes the reflection at $x=0$, and
$\mc S$ imposes a self-adjoint boundary condition at $x=\ell$.  Thus
both equations \eqref{yper} and \eqref{yab} take the form
\begin{equation}
  \label{Feq}
  \mc F(y^0) := \mc S\,\mc E^\ell\,y^0 = 0.
\end{equation}
Given a fixed entropy profile, denoting $\mc L$ as the linearization
of $\mc E$ about the quiet state $p=\ol p$, the composition
$\mc S\,\mc L = D\mc F(\ol p)$ becomes a classical Sturm-Liouville
operator.

It is essential in our argument that the equation of state be
separable so that the small divisors which arise in the bifurcation
argument below be independent of the quiet state pressure for fixed
entropy profile.  On the other hand, the $k$-mode time-period
$T\equiv T(\ol p)$ of the mode we perturb does depend on the
background quiet state pressure $\ol p$.  This simplified structure
allows us to avoid a Nash-Moser argument requiring the expunging of
certain resonant states, so that we do not need problematic
diophantine estimates~\cite{BW,CW,TYperNM}.  In the bifurcation
argument, a variable count yields that changing the background
pressure $\ol p$ is necessary to solve the bifurcation equation.  This
in turn requires that the equation of state be genuinely nonlinear so
that the time period $T$ does depend non-trivially on the background
pressure.  This dependence of $T$ on $\ol p$ is interpreted as
treating $T(\ol p)$ as a Hopf bifurcation parameter.  This dependence
is treated explicitly as a separate operator in our exposition of the
argument presented in~\cite{TYperTofp}.  For clarity of exposition,
below we will treat the dependence of $T(\ol p)$ on $\ol p$
implicitly.

\begin{theorem}
  \label{thm:tile}
  Any solution of the scalar equation \eqref{yper} with even data
  $y^0$ provides a tile which generates by reflections and
  translations, a solution of the compressible Euler equations which
  is $4\ell$-periodic in space and $T$-periodic in time.  A solution
  of \eqref{yab} similarly produces a solution which is
  $2\ell$-periodic in space and $T$-periodic in time.
\end{theorem}

\begin{proof}
  Given a solution $y(x,t)$ defined on $[0,\ell]\times[0,T]$ and
  satisfying the stated conditions, we extend this to a full periodic
  solution of the system \eqref{1Dsys}, illustrated in
  Figure~\ref{fig:period}, as follows.  Recalling that $y$ satisfies a
  scalar equation, we generate the state $(p,u)$ from $y$, and then
  reflect this $2\times2$ state using the boundary conditions.

  First, we require the entropy $s(x)$ to be even and
  $2\ell$-periodic, so set
  \[
    s(x) := 
    \begin{cases}
      s(2\ell-x), &\ell\le x\le 2\ell,\\
      s(x-2\ell), & 2\ell\le x\le 3\ell,\\
      s(4\ell-x), & 3\ell\le x\le 4\ell,
    \end{cases}
  \]    
  Next, recalling \eqref{ydef}, we set
  \[
    \begin{aligned}
      p(x,t) &:= \IR+y(x,t) = \frac{y(x,t)+y(x,-t)}2,\\
      u(x,t) &:= \IR-y(x,t) = \frac{y(x,t)-y(x,-t)}2,
    \end{aligned}
  \]
  for $(x,t)\in[0,\ell]\times[0,T]$.  Recalling that $\chi=1$ or
  $\chi=0$, we extend $p$ using
  \eqref{yrefl} and \eqref{ulsymm} as
  \[
    p(x,t) :=
    \begin{cases}
      p(2\ell-x,t+\chi\,T/2), &\ell\le x\le 2\ell,\\
      p(x-2\ell,t+\chi\,T/2), & 2\ell\le x\le 3\ell,\\
      p(4\ell-x,t), & 3\ell\le x\le 4\ell,
    \end{cases}
  \]
  and similarly extend $u$ as
  \[
    u(x,t) :=
    \begin{cases}
      - u(2\ell-x,t+\chi\,T/2), &\ell\le x\le 2\ell,\\
      u(x-2\ell,t+\chi\,T/2), & 2\ell\le x\le 3\ell,\\
      - u(4\ell-x,t), & 3\ell\le x\le 4\ell.
    \end{cases}
  \]
  By construction, the solution so given is periodic and
  continuous at all boundaries $x=j\ell$, $j=0,\dots 4$.  This uses
  the fact that an even or odd that is $T$-periodic has two symmetry
  axes, $t=0$ and $t=T/2$, as in \eqref{reflper}.  In the case of the
  acoustic boundary condition, we have $\chi=0$ and the functions $p$
  and $u$ are $2\ell$-periodic.
\end{proof}

For a given entropy profile $s(x)$, each solution of \eqref{yper}
generates a space and time periodic solution, generated by a
reflection at the left boundary and a shifted reflection on the right,
so the space period of the solution is $4\ell$.  On the other hand, a
solution of \eqref{yab} generates a periodic solution by just one
reflection in $x$, so has space period $2\ell$.  We will show that the
two problems are related as follows: any $k$-mode solution of the
linearization of \eqref{yper} with $k$ even is also a solution of that
of \eqref{yab}, and conversely an even $2j$-mode linearized solution
of \eqref{yab} can also be realized as a solution of the linearization
of \eqref{yper}.  Moreover, if these linearized solutions perturb,
then by uniqueness they will coincide as solutions of the
corresponding nonlinear problems.

It appears that \eqref{yper} generates more solutions than
\eqref{yab}, but we regard \eqref{yab} as a more physically relevant
problem.  This is because the acoustic boundary condition \eqref{yab}
is simply a reflection off a wall, so the system models sound waves
bouncing between two walls which bound an unrestricted varying entropy
profile on $[0,\ell]$.  On the other hand, the velocity $u$ need not
vanish at $x=\ell$ with \eqref{yper}, and it is hard to ascribe
\eqref{yper} to a physical condition, because we expect that
controlling an entropy profile to be perfectly symmetric between two
walls (at $x=0$ and $x=2\ell$) is practically infeasible.  However, we
regard periodic solutions in which compression and rarefaction are in
perfect balance as being of fundamental importance.

In fact, the reflection symmetry imposed by the periodic boundary
condition implies that even mode solutions of the periodic boundary
value problem (BVP) automatically satisfy the acoustic BVP $u=0$ at
$x=\ell$ for both the linear and nonlinear solutions.  Thus it
suffices to consider only the periodic BVP, even mode solutions of
which are automatically the even mode solutions of the acoustic BVP.

\subsection{Factorization of the Nonlinear Operator}

Our use of symmetry replaces the \emph{periodic return problem} \eqref{F1},
namely
\begin{equation}
  \label{F2}
  \big(\mc F_1-\mc I\big)(p^0,u^0) = 0,
\end{equation}
with the \emph{tiling projection problems} \eqref{yper} or \eqref{yab},
respectively.  The major difficulty in solving \eqref{F2} by
Nash-Moser is the problem of controlling the decay rate of the small
divisors in the perturbed linearized operators at each iterate.

The authors breakthrough revelations were that, after fully
understanding the symmetries~\eqref{symm}, periodicity could be
imposed by the tiling projection operator \eqref{yper} or \eqref{yab},
respectively, which led to the discovery that the small divisors for
fixed entropy profile are uniformly independent of the constant
background pressure.  Surprisingly, we are then able to construct a
fully nonlinear operator which factors into a product of a nonlinear,
nonsingular operator $\mc N$, followed by a \emph{fixed},
\emph{constant}, \emph{diagonal} linear operator $\mc A$ which
uniformly encodes the small divisors, even as the constant state is
perturbed in the iteration.  In this sense, by exploiting symmetry and
restricting to the smaller domain of even functions, we have
``shrink-wrapped'' the problem by reducing the solution space while
retaining the fundamental nonlinear structure of periodic solutions,
namely the balance of compression and rarefaction along each
characteristic.

For fixed time period $T$, denote the nonlinear evolution operator by
\begin{equation}
  \label{Edef}
  \mc E^x:H^b\to H^b, \qquad \mc E^x(y^0) := y(x,\cdot),
\end{equation}
where $y(x,t)$ is the solution of the nonlinear equation \eqref{ycl},
and $H^b$ denotes the Sobolev space of $T$-periodic functions with
norm
\[
  \| y \| = \Big(\sum y_j^2\,j^{2b}\Big)^{1/2}, \com{where}
  y = \sum y_j \,e^{i\,j\,2\pi/T},
\]
so that $y^j$ are the Fourier series coefficients of $y$.

We restate the local existence theory for classical solutions in the
following theorem which will be established in Corollary~\ref{cor:le}
and Theorem~\ref{thm:D2Evol} below.  This holds as long as the data
$y^0$ or evolution interval $[0,X]$ is small enough that shocks do not
form, so the classical theory applies.

\begin{theorem}
  \label{thm:D2E}
  For Sobolev index $b>5/2$, the nonlinear evolution operator
  $\mc E^X:U\subset H^b\to H^b$ is twice Frechet differentiable, with
  first and second Frechet derivatives around $y^0\in H^b$ denoted by
  \[
    D\mc E^X(y^0)[\cdot]:H^b\to H^{b-1}, \com{and}
    D^2\mc E^X(y^0)[\cdot,\cdot]:H^b\times H^b\to H^{b-2},
  \]
  and these are bounded linear and bilinear operators, respectively.
\end{theorem}

Note that for us, more differentiability, that is larger values of
$b$, is better because it establishes greater regularity in the
ultimate nonlinear solutions obtained, which can be constructed to
have arbitraily high number of derivatives.  For finite
$b$, the reduction of differentiability in the ranges of
these operators reflects the principle that \emph{when
  differentiating, a derivative is lost}~\cite{Ydiff}.

Note that before shock formation, $\mc E^x$ is regular and invertible,
its inverse being backwards evolution.  The linearization of $\mc E^x$
as an operator is the Frechet derivative
\begin{equation}
  \label{DEY}
  D\mc E^x(y^0):H^b\to H^b, \qquad
  \mc E^x(y^0)\big[Y^0\big] := Y(x,\cdot),
\end{equation}
where $Y(x,t)$ is the solution of the linearized equation
\begin{equation}
  \label{Yeq}
  \del_x Y + \big(Dg(y)[Y]\big)_t = 0, \qquad Y(0,\cdot) = Y^0,
\end{equation}
where, differentiating \eqref{gy}, we get
\[
  Dg(y)[Y] = \IR-\,Y + \frac{\del v}{\del
    p}\big(\IR+y,s\big)\,\IR+\,Y.
\]
We further denote the linearized evolution operator around the quiet
base state $\ol p$, by
\begin{equation}
  \label{Ldef}
    \mc L^x := D\mc E^x(\ol p),
\end{equation}
and note that this is again invertible by backwards evolution.

It follows that the linearizations of the nonlinear operators
$\mc F_P$ and $\mc F_A$, around $y^0$ are
\[
  \begin{aligned}
    D\mc F_P(y^0) &= \IR-\,\mc T^{-T/4}\,D\mc E^\ell(y^0) \com{and}\\
    D\mc F_A(y^0) &= \IR-\,D\mc E^\ell(y^0),
  \end{aligned}
\]
respectively.  Note that using the notation \eqref{Shat}, these can be
unified as the single equation
\[
  D\mc F(y^0) = \IR-\,\mc T^{-\chi\,T/4}\,D\mc E^\ell(y^0).
\]

\begin{theorem}
  \label{thm:fact}
  For fixed period $T$, the nonlinear operators $\mc F_P$ and
  $\mc F_A$ given in \eqref{yper}, \eqref{yab}, respectively can be
  factored as
  \begin{equation}
    \label{fact}
    \mc F(y^0) = \mc A\,\mc N(y^0), \com{where}
    \mc N := (\mc L^\ell)^{-1}\,\mc E^\ell(y^0),
  \end{equation}
  and $\mc A := D\mc F(\ol p)$ is given by
  \begin{equation}
    \label{Adef}
    \mc A := \mc S\,\mc L^\ell,
  \end{equation}
  where $\mc S$ is given by \eqref{Shat}.  Here $\mc A$ is a fixed,
  constant linear operator, and $\mc N$ is a regular nonlinear
  operator satisfying
  \[
    D\mc N(\ol p) = \mc I,
  \]
  the identity map.  It follows that
  \begin{equation}
    \mc F(y^0) = 0 \com{iff}
    \mc N(y^0) \in \ker\mc A.
    \label{Neq}
  \end{equation}    
\end{theorem}

\begin{proof}
  Because $\mc L^\ell$ is invertible by backwards linearized evolution
  around $\ol p$, the operator $\mc N$ is well defined.  The
  conclusion now follows by direct substitution.
\end{proof}

The upshot of this factorization is that we will invert the nonlinear
part $\mc N$ of the operator by an application of the \emph{standard}
implicit function theorem on Hilbert spaces.  This effectively reduces
the problem of finding nonlinear periodic solutions to that of finding
those entropy profiles for which the fixed operator $\mc A$ is
invertible on the complement of a two-dimensional kernel, one
dimension being the 0-mode constant states, and the other providing
linearized modes which we will perturb.  Our purpose is to use
Sturm-Liouville theory to find appropriate time periods and obtain
conditions under which $\mc A$ has this property.  By this, we solve
\eqref{Neq} without having to resort to a more technical Nash-Moser
iteration and without the need for estimates on the decay rates of the
small divisors.

\section{Local Existence}
\label{sec:exist}

We recall the local $H^b$ existence theory, which applies for
solutions prior to the formation of shocks.  Because we consider only
periodic functions, $y(t+T)=y(t)$, we consider the Sobolev norm on
$H^b([0,T])$, which can be defined using Fourier series, namely
\begin{equation}
  \label{Hbnorm}
  \| y \|_b = \Big(\sum y_j^2\,(1+j^2)^b\Big)^{1/2}, \com{where}
  y = \sum y_j \,e^{i\,j\,2\pi/T},
\end{equation}
so that $H^b$ is a Hilbert space for all real values of
$b$~\cite{Taylor}.  Note also that for $b>1/2$, Sobolev functions are
continuous,
\[
  H^b([0,T])\subset C^0([0,T]) \subset L^\infty([0,T]).
\]

We begin by recalling the abstract local existence theorem for
(multi-dimensional) quasilinear systems of the form
\begin{equation}
  \label{ql}
  A_0\,U_t + \sum_{i=1}^d A_i\,U_{x_i} = 0,
\end{equation}
where $A_i=A_i(x,t,U)$ are symmetric.

\begin{theorem}[Kato, Majda]
  \label{thm:le}
  Let $b>d/2+1$, and assume that $U^0\in H^b$, and $U^0$ takes values
  in a neighborhood $G$ which is compactly contained in the domain of
  the $A_j$.  Then there is a unique classical solution
  $U(x,t)\in H^b$ for $t\in[0,t_*]$, with
  \[
    U\in C([0,t_*],H^b)\cap C^1([0,t_*],H^{b-1}),
  \]
  where the time of existence $t_*>0$ depends on $G$ and
  $\big\|U^0\big\|_b$.  Suppose also that the \emph{a priori}
  estimates
  \begin{equation}
    \label{gradBd}
    \big\|\nabla U\big\|_\infty \le K_1, \qquad
    \Big\|(A_0)_t + \sum_{i=1}^d ({A_j})_{x_j}\Big\|_\infty \le K_2,
  \end{equation}
  $U(x,t)\in \ol{G_1}\subset\subset G$, hold for $t\in[0,t_*]$.
  Then in addition $U$ satisfies the \emph{a priori} $H^b$ bound
  \begin{equation}
  \label{Hbbd}
    \max_{t\in[0,t_1]}\big\|U(\cdot,t)\|_b
    \le C\,e^{(K_1+K_2)\,C\,t_1}\,\big\|U^0\big\|_b,
  \end{equation}
  for $0\le t_1\le t_*$, where $C$ depends only on $b$ and $\ol{G_1}$.
\end{theorem}

This theorem is taken from Majda~\cite{Majda}, who proves it for
integer values of $b$; for a proof including non-integer values of
$b$, and for compact domains, see \cite{Taylor}.  Note that the
same result applies for $t_*<0$, because the evolution is reversible
in the absence of shocks, that is before gradient blowup.

We use this local existence theorem to set up a precise functional
framework in which we can apply the implicit function theorem.  To do
so, we start with the system \eqref{genls}, namely
\[
  p_x + u_t = 0, \qquad
  u_x - v(p,s)_t = 0,
\]
in which the known entropy field $s=s(x)$ is piecewise $H^b$, and we
are treating $x$ as the evolution variable (which simply amounts to a
relabeling of variables).  We want to show that by choosing the data
$\big(p^0(\cdot)-\ol p,u^0(\cdot)\big)$ small enough in $H^b:=H^b(0,T)$,
the solution $\big(p(x,\cdot)-\ol p,u(x,\cdot)\big)$ remains in $H^b$,
and the evolution operator
\[
  \mc E^x:(p^0-\ol p,u^0)\mapsto \big(p(x,\cdot)-\ol p,u(x,\cdot)\big),
  \com{is bounded} H^b\to H^b,
\]
uniformly for $x\in[0,\ell]$.  In order to apply Theorem~\ref{thm:le},
we thus require \emph{a priori} derivative bounds.  In order to obtain
such bounds we examine the growth of derivatives along
characteristics, generalizing an argument of Lax~\cite{Lax64}.

\begin{lemma}
  \label{lem:Linf}
  Let $b>3/2$ be given.  Given any piecewise $H^b$ entropy profile
  defined on the compact interval $[0,\ell]$, there is some $\eta>0$,
  and constants $K_i>0$, such that if
  \[
    \big(p^0-\ol p,u^0\big) \in H^b, \com{with}
    \big\|(p^0-\ol p,u^0)\big\|_b < \eta,
  \]
  then the solution $\big(p(x,\cdot)-\ol p,u(x,\cdot)\big)$ remains in
  $H^b$ for all $x \le \ell$, and we have the estimates
  \[
    \begin{aligned}
      \big\|\big(p(x,\cdot)-\ol p,u(x,\cdot)\big)\big\|_b
      &\le K_1\,\big\|(p^0-\ol p,u^0)\big\|_b,\\
      \big\|\del_t\big(p(x,\cdot),u(x,\cdot)\big)\big\|_\infty
      &\le K_2\,\big\|\del_t(p^0,u^0)\big\|_\infty,
    \end{aligned}
  \]
  uniformly for $x\in[0,\ell]$.
\end{lemma}

\begin{proof}
  According to Majda's theorem, all that needs to be proved is the
  $L^\infty$ derivative bound \eqref{gradBd}.  Initially supposing
  that the entropy $s$ is $H^b$ on all of $[0,\ell]$, differentiate
  the system \eqref{genls} in $t$.  Using the notation
  $\dot \square := \del_t \square$, we immediately get
  \[
    \Big(\frac1\sg\,\sg\,\dot p\Big)_x +\dot u_t = 0, \qquad
    \dot u_x + (\sg\,\dot p)_t = 0,
  \]
  where $\sg^2 = -v_p(p,s)$.  Using Leibniz' rule and rearranging,
  this can be written as
  \[
    \begin{aligned}
      (\sg\,\dot p)_x + \sg\,\dot u_t - \frac{\sg_x}{\sg}\,\sg \dot p
      &=0,\\
      \dot u_x + \sg\,(\sg\,\dot p)_t + \sg_t\,\sg\dot p &=0.
    \end{aligned}
  \]
  We now evaluate
  \[
    \sg_t = \frac{\del\sg}{\del p}\,\dot p, \com{and}
    \sg_x = -\frac{\del\sg}{\del p}\,\dot u +
    \frac{\del\sg}{\del s}\,d_xs,
  \]
  where we have used \eqref{genls} to write $p_x = -\dot u$.
  Thus our derivative equations are
  \begin{equation}
    \label{pudot}
    \begin{gathered}
      (\sg\,\dot p)_x + \sg\,\dot u_t + L_p \,\dot u\,\sg\dot p
      = L_s\,d_xs\,\sg\dot p,\\
      \dot u_x + \sg\,(\sg\,\dot p)_t + L_p\,(\sg\,\dot p)^2 =0,
    \end{gathered}
  \end{equation}
  where we have abbreviated $L(p,s) = \log\sg(p,s)$.

  We now evaluate \eqref{pudot} along forward and backward
  characteristics, given by $\frac{dt}{dx} = \pm\sg$, respectively.
  Using Lax's notation
  \[
    \square' := (\del_x + \sg\,\del_t)\square, \com{and}
    \square^\bpr := (\del_x - \sg\,\del_t)\square,
  \]
  for these characteristic derivatives, adding and subtracting in
  \eqref{pudot} yields
  \[
    \begin{aligned}
      \big(\sg\,\dot p + \dot u\big)'
      + L_p\,(\sg\,\dot p + \dot u)\,\sg\,\dot p
      &= L_s\,d_xs\,\sg\dot p,\\
      \big(\sg\,\dot p - \dot u\big)^\bpr
      - L_p\,(\sg\,\dot p - \dot u)\,\sg\,\dot p
      &= L_s\,d_xs\,\sg\dot p.
    \end{aligned}
  \]
  We now set
  \begin{equation}
    \label{zydef}
    z := \sg\,\dot p + \dot u, \com{and}
    y := \sg\,\dot p - \dot u,
  \end{equation}
  so that $\sg\,\dot p = (z+y)/2$, and write this as the Riccati
  system
  \begin{equation}
    \label{yzsys}
    \begin{aligned}
      z' &= - \half\,L_p\,z\,(z+y) + \half\,L_s\,d_xs\,(z+y),\\
      y^\bpr &= \half\,L_p\,y\,(z+y) + \half\,L_s\,d_xs\,(z+y),
    \end{aligned}
  \end{equation}
  where the coefficients $L_p$ and $L_s$ depend only on $p$ and $s$.

  System \eqref{yzsys} generalizes Lax's system for \x2 systems to
  Euler, the main difference here being the linear terms induced by
  the varying entropy profile.  This cannot be regarded as a \x2
  system of ODEs, because the derivatives are in different
  characteristic directions.  Nevertheless, following
  Lax~\cite{Lax64}, we compare to a scalar Riccati equation, as
  follows.

  For each $x\in[0,\ell]$ define
  \begin{equation}
    \label{LpMs}
    \begin{aligned}
    \ol{L}(x) &:= \sup_{p\in K_p}\big|L_p\big(p,s(x)\big)\big|,\\
    \ol{M}(x) &:= \sup_{p\in K_p}\big|L_s\big(p,s(x)\big)\,d_xs\big|,
    \end{aligned}
  \end{equation}
  where $K_p\subset(0,\infty)$ is a compact interval containing
  $\ol p$ and all $p(x,\cdot)$, and $s\in H^b$.  Now suppose that
  $Z^0>0$ is an initial bound for $|z^0|$ and $|y^0|$,
  \[
    Z^0 := \sup_t\big\{\big|z^0(t)\big|,\;\big|y^0(t)\big|\big\},
  \]
  and define $Z(x)$ by the equation
  \begin{equation}
    \label{Zeq}
    \frac{dZ}{dx} = \ol{L}\,Z^2 + \ol{M}\,Z, \qquad Z(0) = Z^0.
  \end{equation}
  It follows from comparing \eqref{yzsys} to \eqref{Zeq} that
  \[
    |z(t,x)|,\;|y(t,x)| \le Z(x),
  \]
  for all $t\in[0,2\pi)$, $x\in[0,\ell]$.  Note that at an extremum of
  $z$, we have $\del_tz = 0$, so also $z' = \del_xz$ there, and
  similarly for $y$.

  We now solve the comparison ODE \eqref{Zeq} explicitly, by treating
  it as a Bernoulli equation: dividing by $Z^2$ and rearranging, we get
  \[
    \frac d{dx}\Big(\frac1Z\Big) + \ol{M}\,\frac1Z = -\ol L,
  \]
  which is linear, and in turn yields
  \[
    e^{\ol I(x)}\,\frac1{Z(x)} - \frac1{Z^0} =
    - \int_0^xe^{\ol I(\xi)}\ol L(\xi)\;d\xi, \qquad
    \ol I(x):=\int_0^x\ol M(\xi)\;d\xi.
  \]
  Solving, we get the explicit expression
  \[
    Z(x) = \frac{e^{\ol I(x)}\,Z^0}{1-Z^0\,\int_0^xe^{\ol I}\ol L\,d\xi}.
  \]
  Thus as long as we choose the data small enough that
  \[
    Z^0\,\int_0^\ell e^{\ol I}\ol L\,d\xi =: \eta <1,
  \]
  then we have
  \[
    |y(x,\cdot)|,\;|z(x,\cdot)| \le Z(x) \le
    \frac{e^{\ol I(x)}}{1-\eta}\,Z^0,
    \com{for all} x\le \ell.
  \]
  Finally using \eqref{zydef}, we get the estimate
  \begin{equation}
    \label{updot}
    \big\|\big(\dot p(x,\cdot),\dot u(x,\cdot)\big)\big\|_\infty
    \le \frac{\max \sg}{\min\sg}\,\frac{e^{\ol I(x)}}{1-\eta}\,
    \big\|(\dot p^0,\dot u^0)\big\|_\infty,
  \end{equation}
  uniform for $x\in[0,\ell]$.

  For piecewise $H^b$ entropy fields, we have a finite number of
  intervals on which $s$ is bounded and $H^b$, and at entropy jumps,
  both $p$ and $u$ are continuous in $x$.  This means that $\dot p$
  and $\dot u$, being time derivatives, are unchanged across an
  entropy jump.  We thus apply estimate \eqref{updot} inductively on
  the number of subintervals in $[0,\ell]$.
\end{proof}

Note that this derivative estimate extends to entropy fields $s(x)$
that are $BV$: if $s\in BV$, then $\ol M$ in \eqref{LpMs} can be
regarded as a finite measure, and the estimates are essentially
unchanged, but we will not pursue this here.  However, the $BV$ norm
is an appropriate measure of the entropy field.

\begin{corollary}
  \label{cor:le}
  Given any $b>3/2$, tile lengths $\ell$ and $T$, compact interval
  $[p_\flat,p_\sharp]$ with $0<p_\flat<\ol p<p_\sharp$, and
  constant $S$.  There are constants $\eta>0$ and $K<\infty$,
  depending only on $(b,\ell,T,p_\flat,p_\sharp,S)$, such that, if the
  piecewise $H^b$ entropy $s(x)$ and initial values $(p^0,u^0)$
  satisfy
  \[
    TV_{[0,\ell]}\big(s(\cdot)\big) < S, \com{and}
    \big\|( p^0-\ol p, u^0)\big\|_b < \eta,
  \]
  then there is a unique bounded solution $(p,u)\in H^b$ of
  \eqref{genls}, defined for all $x\in[-\ell,\ell]$, and denoted by
  \[
    \mc E^x(p^0,u^0) := \big(p(x,\cdot),u(x,\cdot)\big),
  \]
  which is continuous in both $x$ and $t$ and satisfies \eqref{genls}
  almost everywhere.  Moreover, this solution satisfies the estimate
  \begin{equation}
    \label{puest}
    \big\|\mc E^x( p^0, u^0)-(\ol p,0)\big\|_b
    \le K\,\big\|( p^0-\ol p, u^0)\big\|_b,
  \end{equation}
  uniformly for $x\in[-\ell,\ell]$.
\end{corollary}

Despite the jumps in entropy, which are contact discontinuities, the
solutions $(p,u)$ are classical in that they are everywhere bounded
and continuous, and the PDEs \eqref{genls} are satisfied almost
everywhere.  The velocity is not differentiable (in $x$) at the
entropy jumps $\pm x_j$, $j=1,\dots,J$, say, although left and right
Dini derivatives are defined there, and the only sets on which the
equation is not satisfied pointwise are the lines
$\big\{(\pm x_j,t)\big\}$.

\begin{proof}
  This follows immediately from Theorem~\ref{thm:le} and
  Lemma~\ref{lem:Linf}.  Since $s(x)$ is $BV$ with
  $TV_{[0,\ell]}s(\cdot)\le S$, $s(x)$ is bounded, and we
  obtain fixed bounds for the coefficients $\ol L$ of \eqref{LpMs} as
  well as the coefficient of \eqref{updot} by taking suprema over
  $p\in[p_\flat,p_\sharp]$ and $s\in [s(0)-S,s(0)+S]$.  On the other
  hand, $\ol M$ appears in \eqref{updot} only through the integral 
  \[
    \int_0^{|x|} \ol M(\xi)\;d\xi \le
    \sup_{p,s}\big|L_s(p,s)\big|\,\int_0^{|x|}|d_xs(\xi)|\;d\xi
    \le \sup_{p,s}\big|L_s(p,s)\big|\,S.
  \]
  Since the bound $\|y\|_\infty<C\,\|y\|_b$ depends only on
  the same constants, the proof is complete.
\end{proof}

\section{Linearization and Sturm-Liouville}
\label{sec:lin}

We construct solutions of \eqref{genls} by perturbing solutions of the
linearized equations.  In order to carry this out we require a
detailed analysis of this linearized system.  We begin by noting that
the quiet state $(\ol p,0)$ is a solution of the Euler equations
satisfying the boundary conditions; although it is a constant solution
of \eqref{genls}, it is part of a non-constant standing wave solution
of \eqref{lagr}.  Denoting nonlinear evolution through $x$ according
to \eqref{ycl} by $\mc E^x$, we denote the linearization of this
around the constant quiet state $(\ol p,0)$ by
$\mc L^x := D\mc E^x(\ol p,0)$.  That is, $\mc L^x$ denotes evolution
according to the system
\begin{equation}
  \label{UP}
  \begin{gathered}
    P_x + U_t = 0, \qquad
    U_x + \sg^2\,P_t = 0,\\ \com{where}
    \sg = \sg(x) := \sqrt{-v_p(\ol p,s)},    
  \end{gathered}
\end{equation}
in which we use the convention that $(P,U)$ solve linearized
equations, while $(p,u)$ solve the nonlinear system.  By restricting
to our symmetry class, we equivalently obtain the nonlocal scalar
linearization of \eqref{ycl}, \eqref{gy}, namely, setting $Y := P+U$,
the linearization is
\begin{equation}
  \label{Ylin}
  Y_x + \big(\IR-Y+\sg^2\,\IR+Y\big)_t = 0.
\end{equation}
Here we regard the entropy field $s(x)$ as a given piecewise smooth
function on the interval $x\in[0,\ell]$, continuous at the endpoints.
The authors do not know of previous analyses of this linearized system
when the entropy is piecewise smooth.

\subsection{Sturm-Liouville Evolution}

We solve \eqref{UP} using separation of variables, in which we look
for $T$-periodic solutions in which $P$ and $U$ are even and odd,
respectively.  To do so, we set
\begin{equation}
  \label{UPsub}
  \begin{aligned}
    P(x,t) &:= \sum_{n\ge 0} \vp_n(x)\,\c\big(n\piot t\big), \\
    U(x,t) &:= \sum_{n> 0} \psi_n(x)\,\s\big(n\piot t\big).
  \end{aligned}
\end{equation}
Plugging in to \eqref{UP} and simplifying, we get the \emph{Sturm-Liouville}
(SL) system
\begin{equation}
  \label{SL1}
  \dot\vp_n + \w\,\psi_n = 0, \qquad
  \dot\psi_n - \sg^2\,\w\,\vp_n = 0,  \qquad
  \w := n\piot,
\end{equation}
where $\dot\square$ now denotes $\frac{d}{dx}\square$.

We use the following convenient notation: for a given (reference) time
period $T$, denote the basis $k$-mode by the $1\times2$ matrix
\begin{equation}
  \mc B_k = \mc B_k(T) := \(\c(k\frac{2\pi}Tt) & \s(k\frac{2\pi}Tt)\),
  \label{Tk}
\end{equation}
where $\c/\s$ abbreviate the trigonometric functions $\cos/\sin$.  It
then follows that any $T$-periodic function can be represented as
\[
  f(t) = \sum_{k\ge 0} \mc B_k\(a_k\\b_k\)
  = \sum a_k\,\text{c}(k\TS{\frac{2\pi}T}t)
  + b_k\,\text{s}(k\TS{\frac{2\pi}T}t),
\]
uniquely with $b_0=0$.  In this notation, \eqref{UPsub} can be written
\begin{equation}
  \label{YasT}
  Y(x,t) = P + U = \sum_{n\ge 0} \mc B_n\,\(\vp_n(x)\\\psi_n(x)\),
\end{equation}
and we obtain the following Lemma.

\begin{lemma}
  \label{lem:Lx}
  The linearized evolution operator $\mc L^x = D\mc E^x(\ol p)$ acts
  on the $k$-mode basis vector $\mc B_k = \mc B_k(T)$ given by
  \eqref{Tk}, by
  \[
    \mc L^x\,\mc B_k = \mc B_k\,\Psi\big(x;k\piot\big), \com{so}
    \mc L^x\,\mc B_k\(a\\b\) = \mc B_k\,\Psi\(a\\b\),
  \]
  where $\Psi = \Psi(x;\w)$ is the fundamental solution of
  \eqref{SL1}, so satisfies
  \begin{equation}
    \label{Psi}
    \dot\Psi(x;\w) = \w\,
    \(0&-1\\\sg^2(x)&0\)\,\Psi(x;\w),
    \qquad \Psi(0;\w) = I.
  \end{equation}
\end{lemma}

\begin{proof}
  Equation \eqref{Psi} follows by substitution of \eqref{YasT} into
  \eqref{UP} and applying standard separation of variables.
\end{proof}

Alternatively, we can express the linearized system \eqref{UP} as a
wave equation with varying speed, namely
\begin{equation}
  \label{Peq}
  P_{xx} - \sg^2\,P_{tt} = 0,
\end{equation}
and after use of \eqref{UPsub}, by separation of variables, we get the
second-order linear system
\begin{equation}
  \label{SL2}
  \ddot \vp_n + \big(n\piot\big)^2\,\sg^2\,\vp_n = 0,
\end{equation}
which is the SL system \eqref{SL1} expressed as a second-order linear
equation.

We now consider boundary values for this SL system.  In both
\eqref{yper} and \eqref{yab}, the data $y^0$ (and so also $Y^0$) posed
at $x=0$ is even, so using \eqref{YasT}, this becomes
the condition
\begin{equation}
  \label{bc0}
  U(0,\cdot) = 0, \com{equivalently} \psi_n(0) = \dot\vp_n(0) = 0.
\end{equation}
If we are solving \eqref{yab}, we get the same condition at $x=\ell$,
namely
\begin{equation}
  \label{bcab}
  U(\ell,\cdot) = 0, \com{equivalently} \psi_n(\ell)
  = \dot\vp_n(\ell) = 0.
\end{equation}
On the other hand, if we are solving \eqref{yper}, the boundary
condition is
\[
  \begin{aligned}
  0 &= \IR-\,\mc T^{-T/4}\big(P(\ell,\cdot)+U(\ell,\cdot)\big)\\
  &= \IR-\sum\Big(\vp_n(\ell)\,\c\big(n\piot t-n{\TS\frac\pi2}\big)
  +\psi_n(\ell)\,\s\big(n\piot t-n{\TS\frac\pi2}\big)\Big),
  \end{aligned}
\]
which yields the conditions
\begin{equation}
  \label{bcper}
  \begin{gathered}
    \dot\vp_n(\ell) = \psi_n(\ell) = 0, \quad n \text{ even},\\
    \vp_n(\ell) = \dot\psi_n(\ell) = 0, \quad n \text{ odd}.
  \end{gathered}
\end{equation}

Our boundary conditions for the linearized evolution can now be
expressed succinctly using the fundamental solution \eqref{Psi} as
follows.  Recall that our domain is the set of even periodic
functions.

\begin{lemma}
  \label{lem:fund}
  Suppose that constant ambient pressure $\ol p$, entropy profile
  $s(x)$ and reference period $T$ are given.  Then the linearized
  operator $\mc A$ given by \eqref{Adef} is diagonal on Fourier
  $k$-modes and is given by
  \[
      \mc A\Big[\c\big(k\piot t\big)\Big]
      = \d_{k}(T)\,\s\big(k\piot t\big),
  \]
  where the $k$-th \emph{divisor} is defined by
  \begin{equation}
    \label{dkcts}
      \d_{k}(T) := \(0&1\)\,R(k\chi\pb)\,
      \Psi\big(\ell;k\piot\big)\(1\\0\),
  \end{equation}
  where $\chi = 1$ or $\chi=0$ for the periodic or acoustic boundary
  conditions, respectively.  In particular, all divisors are bounded,
  that is, there exists a constant $C_\d$ depending only on the
  entropy profile, such that
  \begin{equation}
    \label{cdelta}
    |\d_k(T)| < C_\d, \com{uniformly for} T\in[1/C,C].
  \end{equation}
\end{lemma}

\begin{proof}
  By \eqref{Tk}, we have $\c\big(k\piot t\big) = \mc B_k
    {\scriptstyle\begin{pmatrix}1\\0\end{pmatrix}}$, and we calculate
  \[
    \IR-\,\mc B_k = \(0&1\)\,\s\big(k\piot t\big) \com{and}
    \mc T^{T/4}\,\mc B_k = \mc B_k\,R(k\pb).
  \]
  The result follows by using these and Lemma~\ref{lem:Lx} in
  \eqref{Adef}.  Since SL evolution is linear, the bound
  \eqref{cdelta} follows by continuity on compact intervals.
\end{proof}

It follows that the single mode $\c\big(k\piot t\big)$ provides a
solution given by $\mc B_k\,\Psi\big(x;k\piot\big)$ of the linearized
problem satisfying the boundary conditions if and only if
\[
  \d_{k}(T)=0,
\]
respectively.  In this case,
$\Psi\big(x;k\piot\big)\scriptstyle\begin{pmatrix}1\\0\end{pmatrix}$
solves \eqref{SL1}, and since it satisfies the boundary conditions, it
is an eigenfunction of the Sturm-Liouville problem \eqref{SL2},
corresponding to eigenvalue $\lambda_k := \big(k\piot\big)^2$.  We
will show that if this $k$-mode is nonresonant, then this perturbs to
a nonlinear pure tone solution of \eqref{genls}.

We are thus led to introduce the Sturm-Liouville (SL) operator
\begin{equation}
  \label{SLop}
  \mc L := - \frac1{\sg^2}\,\frac{d^2}{dx^2}, \com{so that}
  \mc L \phi = - \frac1{\sg^2}\,\ddot\phi,
\end{equation}
and to introduce the associated SL eigenvalue problem,
\begin{equation}
  \label{SLev}
  \mc L\,\vp = \lambda\,\vp, \com{which is}
  - \ddot\vp = \lambda\,\sg^2\,\vp,\com{on} 0<x<\ell,
\end{equation}
and subject to the boundary conditions \eqref{bc0} at $x=0$ and either
\eqref{bcper} or \eqref{bcab} at $x=\ell$.  We note that the boundary
conditions imply self-adjointness, implying that this is a
\emph{regular SL problem} as long as the weight
$\sg^2(x) = -v_p\big(\ol p,s(x)\big)$ is a positive, bounded piecewise
continuous function~\cite{Pryce,BR,CodLev}.  Moreover, because the
corresponding ODEs are linear, the solutions can be analyzed as
integral equations, and we expect our methods to extend to BV entropy
profiles with only minor technical changes.

We collect some classical results for regular SL eigenvalue problems
in the following lemma; see \cite{Pryce,CodLev} for details.  Many of
the statements in this lemma will be explicitly verified in
Section~\ref{sec:SL} below.

\begin{lemma}
  \label{lem:SLprops}
  Assume that the entropy profile is piecewise $H^b$.  The
  Sturm-Liouville system \eqref{SLev}, with initial condition
  \eqref{bc0} and target end condition \eqref{bcper} or \eqref{bcab},
  has infinitely many eigenvalues $\lambda_k$.  These are positive,
  simple, monotone increasing and isolated, and satisfy the growth
  condition $\lambda_k=O(k^2)$.  The corresponding eigenfunctions form
  an orthogonal $L^2$ basis in the Hilbert space with weight function
  $\sg^2$.
\end{lemma}

We label the eigenvalues $\lambda_k$ and corresponding eigenfunctions
$\vp_k$, scaled so that $\vp_k(0)=1$, so that \eqref{SLev} and
\eqref{SL1} hold for each $k\ge1$, and the $\lambda_k$'s are
increasing with $k$.  Because we prefer to work with the frequencies,
we define the $k$-th \emph{eigenfrequency} $\w_k$ and corresponding
\emph{reference period} $T_k$ by
\begin{equation}
  \label{Tn}
  \w_k := \sqrt{\lambda_k}, \com{and} 
  T_k := k\,\frac{2\pi}{\w_k} = k\,\frac{2\pi}{\sqrt{\lambda_k}},
\end{equation}
respectively.  It follows that the $k$-modes
\[
  \begin{aligned}
  Y_k(x,t) &:= \mc B_k(T_k)\,\Psi\big(x;k{\TS\frac{2\pi}{T_k}}\big)
  \,{\TS\begin{pmatrix}1\\0\end{pmatrix}}
  \com{or}\\
    P_k(x,t) &:= \vp_k(x)\,\c\big(k{\TS\frac{2\pi}{T_k}}t\big)
               = \vp_k(x)\,\c(\w_kt),
  \end{aligned}
\]
defined for $x\in[0,\ell]$, solve \eqref{Ylin} or \eqref{Peq},
respectively, while also satisfying the appropriate boundary
conditions.  These $k$-modes lie in the kernel of the linearized
operator $\mc A = D\mc F(\ol p)$, and we will investigate the
conditions under which they perturb to pure tone solutions of the
nonlinear problem.  We note that for given $k$, the frequency $\w_k$
is determined by the entropy profile, and this in turn yields the
reference period.  Consideration of other resonant and nonresonant
modes then refers to this fixed period $T_k$.

\subsection{Nonresonant modes}

Identifying a mode as a $k$-mode requires
us to fix the reference period $T$: clearly
$\mc B_{jk}(jT) = \mc B_k(T)$ for any $j\ge 1$.

Our goal is to perturb the linearized solution $Y_k$ (or $P_k$) to a
time periodic solution of the corresponding nonlinear equation.  Since
we want $u^0=0$, we consider perturbations of the initial data of the
form
\begin{equation}
  \label{pert}
  y^0(0) = p^0(t) := \ol p + \a\,\c\big(k{\TS\frac{2\pi}{T_k}}t\big)
   + z + \sum_{j\ge 1,j\ne k}a_j\,\c\big(j{\TS\frac{2\pi}{T_k}}t\big),
\end{equation}
in which the reference period $T_k$ is determined by the choice of SL
eigenfrequency $\w_k$ (or eigenvalue $\lambda_k=\w_k^2$) in
\eqref{Tn}, and we now regard this as fixed.  We interpret $\a$ as the
amplitude, $z$ as the quiet base state perturbation, and $a_j$ as
nonlinear adjustments, and our goal is to find $z(\a)$ and $a_j(\a)$
which produce solutions to the nonlinear problem.

We declare the $k$-mode to be \emph{nonresonant} if no other
$j$-modes satisfy the boundary condition, that is
\begin{equation}
  \label{nonres}
  \d_{j}(T_k)\ne0
  \com{for all} j\ne k,
\end{equation}
respectively for the appropriate boundary conditions \eqref{yper} or
\eqref{yab}.  Here $\d_j$ is defined by \eqref{dkcts} with $\chi=1$ or
$0$, respectively, where the reference period $T_k$ is fixed by our
choice of mode $k$.  The collection of $j$-mode basis elements
\[
  \big\{\mc B_j(T_k)\;|\;j\ge 1\big\} =
  \Big\{\(\c(j\piotk t)&\s(j\piotk t)\)\;|\;j\ge 1\Big\}
\]
span the even and odd functions of period $T_k$, respectively.

\begin{lemma}
  \label{lem:nonres}
  Fix the $k$-mode, which fixes the time period $T=T_k$.  The
  $j$-mode, defined by this time period, $\c(j\piotk t)$, is resonant
  with the fixed $k$-mode if and only if two distinct SL frequencies
  are related by the condition
  \[
    \d_j(T_k) = 0 \com{iff}
    k\,\w_l = j\,\w_k,
  \]
  for some index $l = l(j)\ne k$.
\end{lemma}

\begin{proof}
  By construction, we have
  \[
    \d_k(T_k) = 0, \com{with}
    T_k := k\,\frac{2\pi}{\w_k}.
  \]
  If $\d_j(T_k) = 0$, then $j\frac{2\pi}{T_k}$ also corresponds to
  some SL frequency $\w_l$, so that
  \[
    \w_l = j\frac{2\pi}{T_k} = \frac jk\,\w_k,
  \]
  and since $j\ne k$, we have $\w_l\ne\w_k$, which implies $l\ne k$.
\end{proof}

The resonance or nonresonance of modes corresponding to SL eigenvalues
depends on the entropy profile and constitutive law.  Moreover,
because the SL eigenvalues, and hence also frequencies, vary
continuously with the entropy profile, one should expect that the set
of profiles with resonant modes should be small in the appropriate
sense, so that generically, all modes are nonresonant.  We state the
following theorem which collects the results we need to complete the
bifurcation argument that nonresonant $k$-mode solutions of the
linearized evolution perturb to pure tone solutions of the fully
nonlinear problem.  The proof of this theorem, which relies on the
introduction of an angle variable, is postponed to
Section~\ref{sec:SL}, and proof that nonresonance is generic is
carried out in Section~\ref{sec:pwc}.

\begin{theorem}
  \label{thm:sl}
  For any piecewise $H^b$ entropy profile $s\in\mc P$ and constant
  background pressure $\ol p$, there is a unique sequence $\{\w_k\}$
  of SL eigenfrequencies, with corresponding eigenfunctions satisfying
  the SL boundary conditions \eqref{bcab} or \eqref{bcper}.  Moreover,
  $\{\w_k\}$ can be labeled as a monotone increasing sequence which
  grows like $k$, that is
  \[
    \lim_{k\to\infty}\frac{\w_k}k\to C>0, \com{so that}
    C'\le\frac{\w_k}k\le C'',
  \]
  for some constants $C$, $C'$, $C''>0$, depending only on the entropy
  profile and the quiet state $\ol p$.  In particular, the
  corresponding time periods given by \eqref{Tn} are uniformly
  bounded, $T_k\in[2\pi\,C',2\pi\,C'']$.  Each of these eigensolutions
  determines a $T_k$ time-periodic solution of the linearized
  equations on $[0,\ell]$ of the form
  \[
    P = \c(\w_k\,t)\,\vp_k(x), \qquad
    U = \s(\w_k\,t)\,\psi_k(x).
  \]
\end{theorem}

We now specialize to a \emph{separable} equation of state,
\eqref{sepEoS}.  For this class we have the following important
corollary, which states that the resonance structure depends only on
the entropy profile.  This independence of resonance structure from
average background pressure is precisely what allows us to complete
the bifurcation argument in the next section.

\begin{corollary}
  \label{cor:sep}
  Suppose that the equation of state is separable.  Then all SL
  eigenvalues factor into a quotient $\w_j = \wh\w_j/a(\ol p)$, in
  which $\wh\w_j$ are \emph{independent} of the constant background
  pressure $\ol p$.  The corresponding periods are given by
  \begin{equation}
    \label{Tn2}
    T_j := j\,\frac{2\pi}{\w_j}
    = j\,\frac{2\pi\,a(\ol p)}{\wh\w_j}
    =: a(\ol p)\,\wh T_j,
\end{equation}
  with $\wh T_j$ again independent of $\ol p$.  In particular, the
  resonance structure depends only on the entropy profile.  Moreover,
  the SL eigenfunctions $\vp_j$ and $\psi_j$ are also independent of
  $\ol p$, up to normalization or scaling by $a(\ol p)$.
\end{corollary}

\begin{proof}
  Using \eqref{sepEoS} in \eqref{UP}, we get
  \[
    \sg^2 = -v_p(\ol p,s) = \frac{-\rho_p}{\rho^2}
    = a(\ol p)^2\,\vs(s)^2,
  \]
  that is
  \begin{equation}
    \label{sas}
    \sg(x) = a(\ol p)\,\vs\big(s(x)\big),
  \end{equation}
  and plugging this into the SL operator \eqref{SLop} yields
  \begin{equation}
    \label{SLop2}
    \mc L \phi = - \frac1{a(\ol p)^2\,\vs(s(x))^2}\,\ddot\phi
    = \frac1{a(\ol p)^2}\,\wh{\mc L}\phi,
  \end{equation}
  where we have defined
  \[
    \wh{\mc L} := a(\ol p)^2\,\mc L
    =  - \frac1{\vs(s(x))^2}\,\frac{d^2}{dx^2},
  \]
  so that $\wh{\mc L}$ is independent of $\ol p$.

  Referring to \eqref{Tn}, we now solve the single eigenvalue problem
  analogous to \eqref{SLev}, namely
  \begin{equation}
    \label{SLev2}
    \wh{\mc L}\,\vp = \wh\w^2\,\vp, \com{on} 0<x<\ell.
  \end{equation}
  Since $\mc L$ and $\wh{\mc L}$ are related by \eqref{SLop2}, the
  conclusions of the corollary follow by standard linear algebra
  arguments.  Similarly substituting \eqref{sas} into \eqref{SL1}, it
  is easy to see that $\psi$ scales by $a(\ol p)$.
\end{proof}

\section{Nonlinear Perturbation of Nonresonant Modes}
\label{sec:bif}

The above development reduces the problem of existence of nonlinear
pure tone solutions of the compressible Euler equations to the problem
of finding a solution $y^0$ of the nonlinear equation
\eqref{fact}, namely
\[
  \mc F(y^0) = \mc A\,\mc N(y^0) = \mc S\,\mc L\,\mc N(y^0) = 0,
  \qquad y^0 \text{ even},
\]
where the nonlinear part $\mc N$ is bounded invertible for small
regular data, and $\mc S\,\mc L$ is the Sturm-Liouville (SL) operator
which has small divisors as described above.  Our goal is to perturb a
linearized $k$-mode solution to a pure tone solution of the nonlinear
equation, assuming that the linearized $k$-mode is nonresonant.  We
thus assume a separable equation of state \eqref{sepEoS} and that the
entropy profile $s(x)$ is fixed and known, and we pick $k$ which is
nonresonant, so that $\w_j/\w_k\notin\B Q$ for $j\ne 0,\ k$.  Recall
that these choices fix the time period $T$ by \eqref{Tn} and
\eqref{Tn2}, namely
\[
  T := T_k = \frac{2\pi\,k}{\w_k}
  = a(\ol p)\,\frac{2\pi\,k}{\wh\w_k}
  = a(\ol p)\,\wh T_k,
\]
in which $\wh T_k$ depends only on the entropy profile.  By
Lemma~\ref{lem:nonres} and Corollary \ref{cor:sep}, the nonresonant
case is characterized by the conditions $\d_j(T_k)\ne0$ for all
$j\ne k$.  Recall that the use of the indicator $\chi$ allows us to
treat both periodic and acoustic boundary value problems \eqref{yper},
\eqref{yab} simultaneously.

We now introduce the Hilbert spaces which allow us to find solutions
which are perturbations of the constant state $\ol y=0$.  To leading
order, these have the form $\a\,\c(\w_k\,t)\in\ker\{\mc S\,\mc L\}$,
and are parameterized by the amplitude $\a$.  Our program is to
construct solutions of the form
\begin{equation}
  \label{y0}
  y^0(t) = \ol p + \a\,\c(\w_k\,t) + z +
  \sum_{j\ne k}a_j\,\c(j\,t\,\piot),
\end{equation}
where the 0-mode $z$ (also in $\ker\{\mc S\,\mc L\}$), and higher mode
coefficents $a_j$ are unknowns to be found, of order $O(\a^2)$.  Here
the period $T$ implicitly depends on the quiet state by $T = a(\ol
p+z)\,\wh T_k$.

Observe first that there is one free variable and one equation
corresponding to each mode because the operator \eqref{fact} takes
even modes to odd modes.  That is, for each $\a$, after nonlinear
evolution, the projection $\IR-$ leaves only the odd components of the
evolved data, so we are left with one equation for each $j$-mode
coefficient $j\ne 0$.  On the other hand, the free parameters are the
$a_j$, $j\ne k$, and 0-mode $z$, so each equation corresponds uniquely
to an unknown, with $z$ parameterizing the $k$-mode whose amplitude is
$\a$.  Thus formally we expect to get a solution of \eqref{fact} of
the form \eqref{y0} for any $\a$ sufficiently small in the nonresonant
case.  The following development is an explicit version of the
Liapunov-Schmidt decomposition of the nonlinear operator $\mc F$ into
the auxiliary equation and corresponding bifurcation equation in the
nonresonant case, see~\cite{Golush}.

As a roadmap, we briefly describe the elements of the abstract
bifurcation problem and its solution, which we will make precise
below.  Consider the problem of solving equations of the form
\[
  F(\a,z,w) = 0, \com{with} F(0,0,0)=0,
\]
for solutions parameterized by the small amplitude parameter $\a$.  In
general the gradient $\nabla_{(z,w)}F\big|_0$ is not invertible, so
this cannot be done with a direct application of the implicit function
theorem.  For the bifurcation argument, one assumes
\begin{equation}
  \label{zw}
  z\in\ker\Big\{\nabla_{(z,w)}F\big|_0\Big\}, \qquad
  w\in\ker\Big\{\nabla_{(z,w)}F\big|_0\Big\}^\perp.
\end{equation}
and that $\frac{\del F}{\del w}\big|_0$ is invertible.  The standard
argument is then to first find
\[
  w(\a,z) \com{so that}
  F\big(\a,z,w(\a,z)\big)=0,
\]
which is known as the \emph{auxiliary equation}.  This leaves the
problem of solving for
\[
  z=z(\a) \com{such that} F\big(\a,z(\a),w(\a,z(\a))\big) = 0.
\]
In our problem, $\frac{\del F}{\del z}\big|_0$ is not invertible, so
$F$ is replaced by the equivalent function
\[
  G(\a,z) :=
  \begin{cases}
    \frac1\a\,F\big(\a,z,w(\a,z)\big), &\a\ne 0,\\[3pt]
    \frac{\del F}{\del \a}\big(0,z,w(0,z)\big), &\a = 0,
  \end{cases}
\]
and the equation $G(\a,z)=0$ is known as the \emph{bifurcation
  equation}.  The bifurcation equation can be solved by the implicit
function theorem, provided
\[
  \frac{\del G}{\del z}\Big|_{(0,0)} \equiv
  \frac{\del^2 F}{\del z\,\del\a}\Big|_{(0,0,0)} \ne 0.
\]
The decomposition of the domain of the function $F$ into a direct sum
of the kernel and its orthogonal complement as in \eqref{zw} is known
as the Liapunov-Schmidt decomposition~\cite{Golush,TYperBif}.

In our application, the infinite dimensional analog of the gradient
$\frac{\del F}{\del z}\big|_0$, required for the auxiliary equation,
is invertible with unbounded inverse because of the presence of small
divisors $\d_k$.  However, the factorization \eqref{fact} means that
these small divisors are fixed, and so we can handle them by an
appropriate adjustment of the associated Hilbert space norms.  It is
remarkable that by factoring the nonlinear operator, we are able to
avoid difficult technical issues, such as diophantine estimates, which
are common in Nash-Moser iterations, see~\cite{TYdiff2} and references
therein.  When solving the bifurcation equation, we must calculate the
second derivative of an infinite dimensional nonlinear evolution
operator.  This section presents the complete bifurcation argument
subject to the proof of Theorem~\ref{thm:D2E} that the nonlinear
evolution operator is twice Frechet differentiable.  The proof of that
theorem is carried out in Section~\ref{sec:D2E} below.

\subsection{Auxiliary Equation}

Assume $y^0$ lies in the Sobolev space $H^b$, so that for small data
$y^0$, the evolution $\mc E^x(y^0)$ stays in $H^b$ for $b>5/2$, by the
local existence theory, which applies as long as gradients remain
finite, as shown in Corollary~\ref{cor:le}.  For us, the actual value
of $b$ is secondary, and the larger $b$ is, the more regular the
periodic solutions constructed, so it suffices to take $b$ arbitrarily
large.

To precisely define our bifurcation problem we must define the
operator $\mc F$ and corresponding Hilbert spaces to carry out the
abstract argument described above.  Recall that the entropy profile
and nonresonant mode $k$ have been fixed and these determined the time
period by \eqref{Tn}, namely
\[
  T = \frac{2\,\pi\,k}{\w_k} = \frac{2\,\pi\,k}{\wh\w_k}\,a(\ol p+z),
\]
and we are perturbing the $k$-mode $\a\,\c(k\,\piot)$ for sufficiently
small amplitudes $\a$.  For given $b$ large enough, define the Hilbert
spaces for the \emph{domain} as
\begin{equation}
  \label{Hsplit}
  \begin{aligned}
    \mc H_1 &:= \big\{\ol p + z +
              \a\,\c(k\,t\,\piot)\;\big|\;z, \a\in\B R\big\}
    \com{and}\\
    \mc H_2 &:= \Big\{ \sum_{j\ne k,j>0}a_j\,\c(j\,t\,\piot)\;\Big|
    \;\sum a_j^2\,j^{2b}<\infty\Big\},    
  \end{aligned}
\end{equation}
so that $y^0$ given by \eqref{y0} is
\[
  y^0 = \ol p + \a\,\c(k\,t\,\piot) + z +
  \sum_{j\ne k,j>0}a_j\,\c(j\,t\,\piot)
  \in\mc H_1\oplus\mc H_2,
\]
and again the period $T$ depends implicitly on $z$.  In the sums
below, we will always implicitly take $j>0$ because we treat the
$0$-mode $z$ explicitly.  Now referring to \eqref{fact}, \eqref{Adef},
define the operator
\begin{equation}
  \label{Fhat}
  \begin{gathered}
    \whc F:\mc H_1\times\mc H_2\to H^b \com{by}\\
    \whc F(y_1,y_2) := \mc F(y_1+y_2) =
    \mc S\mc L\,\mc N\,(y_1+y_2),
  \end{gathered}
\end{equation}
so that $\whc F$ is a continuous nonlinear operator on
$\mc H_1\times\mc H_2$.  Recall that $\mc N = \mc L^{-1}\,\mc E^\ell$,
and Theorem~\eqref{thm:fact} states that $D\mc N(\ol p) = \mc I$.  It
follows that the partial Frechet derivative
\[
  D_{y_2}\whc F(\ol p,0):\mc H_2\to H^b \com{is}
  D_{y_2}\whc F(\ol p,0) = \mc S\,\mc L \Big|_{\mc H_2},
\]
where $\mc S$ projects onto the odd sine modes.  According to
Lemma~\ref{lem:fund}, this acts as
\[
  D_{y_2}\whc F(\ol p,0)\Big[\sum_{j\ne k} a_j\,\c(j\,t\,\piot)\Big]
  = \sum_{j\ne k} a_j\,\d_j\,\s(j\,t\,\piot),
\]
where $T=T\big|_{z=0}$.  By our choice of parameters, we have
$\d_j\ne0$ for $j\ne k$, so that $D_{y_2}\whc F(\ol p,0)$ is
injective, but the inverse is not bounded as a map $H^b\to H^b$ for
any $b$, because of the presence of the small divisors $\d_k$.
However, because the small divisors are fixed, we can define a new
norm on the target space $H^b$ so that $D_{y_2}\whc F(\ol p,0)$
becomes an \emph{isometry}, and in particular is bounded invertible on
its range, as in \cite{TYdiff2}.  Thus we define the scaled $H^b$ norm
on the \emph{range} of $\whc F$ by
\begin{equation}
  \begin{aligned}
  \|y\|^2 &:= \beta^2 + \sum_{j\ne k} c_j^2\,\d_j^{-2}\,j^{2b},
            \com{and set}\\
  \mc H_+ &:= \Big\{ y = \sum_{j\ne k}c_j\,\s(j\,t\,\piot)\;\Big|
  \;\|y\|<\infty\Big\}, \com{and}\\
    \mc H &:= \big\{ \beta\,\s(k\,t\,\piot) \big\}
            \oplus \mc H_+.
  \end{aligned}
  \label{Hplus}
\end{equation}
Here for convenience we have set $\beta = c_k\,k^{2b}$, which isolates
the $k$-mode kernel.  According to \eqref{cdelta}, each divisor $\d_j$
is bounded, $0<|\d_j|\le C_\d$, so that
\[
  \|y\|^2_{H^b} \le C_\d\,\|y\|^2,
  \com{and} \mc H\subset H^b,
\]
and it is clear that $\mc H$ is a Hilbert space.  Finally, let $\Pi$
denote the projection 
\[
  \begin{gathered}
  \Pi:\mc H\to \mc H_+, \com{by}\\
  \Pi\Big[\beta\,\s(k\,t\,\piot)+\sum_{j\ne k}a_j\,\s(j\,t\,\piot)\Big]
  := \sum_{j\ne k}a_j\,\s(j\,t\,\piot).
  \end{gathered}
\]
Note that we have constructed these spaces so that
\[
  \ker\{D_{y_2}\whc F(\ol p,0)\} = \mc H_1, \qquad
  \text{ran}\{D_{y_2}\whc F(\ol p,0)\} = \mc H_+,
\]
so that $\Pi D_{y_2}\whc F(\ol p,0) = D_{y_2}\whc F(\ol p,0)$, and moreover
$D_{y_2}\whc F(\ol p,0):\mc H_2\to\mc H_+$ is an isometry, and thus bounded
invertible on its range $\mc H_+$.  The invertibility of
$D_{y_2}\whc F(\ol p,0)$ allows us to solve the auxiliary equation
  by a regular application of the classical
implicit function theorem on Hilbert spaces.

\begin{lemma}
  \label{lem:himodes}
  There is a neighborhood $\mc U\subset\mc H_1$ of the origin on which
  there is a unique solution $W(\a,z) \in\mc H_2$ of the auxiliary
  equation.  That is there exists a unique $C^2$ map
  \[
  \begin{gathered}
    W:\mc U\to\mc H_2, \com{written}\\
    W\big(z+\a\,\c(k\,t\,\piot)\big) =: W(\a,z) \in\mc H_2,
  \end{gathered}
  \]
  such that, for all $z+\a\,\c(k\,t\,\piot)\in\mc U$, we have
  \begin{equation}
    \begin{aligned}
      \Pi\,&\whc F\big(\ol p+z+\a\,\c(k\,t\,\piot),W(\a,z)\big)  =\\
      &\Pi\,\mc F\big(\ol p+z+\a\,\c(k\,t\,\piot)+W(\a,z)\big) = 0.
    \end{aligned}
    \label{Weq}
  \end{equation}
  Moreover, $W(\a,z)$ satisfies the estimates
  \begin{equation}
    \label{West}
    \lim_{\a\to0}\frac{W(\a,z)}{|\a|} = 0, \com{and}
    \lim_{\a\to0}\frac{\del W(\a,z)}{\del\a} = 0,
  \end{equation}
  uniformly for $z$ in a neighborhood of 0.
\end{lemma}

\begin{proof}
  This is a direct application of the implicit function theorem,
  assuming Theorem~\ref{thm:D2E}, so that $\mc E$ is twice Frechet
  differentiable.  This states that if $\mc H_1$, $\mc H_2$ and $\mc
  H_+$ are Hilbert spaces, and
  \[
    \mc G := \Pi\whc F:\Omega\subset\mc H_1\times\mc H_2\to\mc H_+
  \]
  is a continuously differentiable map defined on an open neighborhood
  $\Omega$ of $(0,0)$, and satisfying $\mc G(0,0)=0$, and if the
  linear (partial derivative) map
  $D_{y_2}\mc G(0,0):\mc H_2\to\mc H_+$ is bounded invertible, then
  there is an open neighborhood $\mc U_1\subset \mc H_1$ of 0 and a
  unique differentiable map $W:\mc U_1\to\mc H_2$, such that
  $\mc G\big(x,W(x)\big)=0$ for all $x\in\mc U_1$, see
  e.g.~\cite{Jost}.  Because we have built our Hilbert spaces so that
  $D_{y_2}\mc G$ is an isometry, the existence of the $C^2$ solution
  $W(\a,z)$ follows immediately.

  The fact that the period $T$ depends implicitly on the quiet state
  perturbation $z$ plays no significant part here because the value of
  $z$ does not materially affect the calculation of the derivative
  $D_{y_2}\mc G$.

  Theorem \ref{thm:D2E} together with the implicit function theorem
  imply that $W(\a,z)$ is twice differentiable.  Since $\mc F(\ol p+z)=0$,
  we have
  \[
    W(\a,z) \in \mc H_2, \com{with} W(0,z) = 0,
  \]
  and, setting $y=z+\a\,\c(k\,t\,\piot) + W(\a,z)$, we have by \eqref{Weq}
  \[
    0 = \Pi\, \mc F(y) =
    \Pi\,\IR-\,\mc L_0\,\mc N\big(\ol p+z+\a\,\c(k\,t\,\piot) + W(\a,z)\big).
  \]
  Differentiating with respect to $\a$ and setting $\a=0$,
  we get
  \[
    \begin{aligned}
      0 = \frac{\del\Pi\mc F(y)}{\del\a}\Big|_{\a=0}
      &= \mc S\,\mc L\,D\mc N(\ol p+z)\Big[\c(k\,t\,\piot) +
      \frac{\del W}{\del\a}\Big|_{\a=0}\Big] \\
      &= \mc S\,\mc L\,D\mc N(\ol p+z)\Big[\frac{\del W}{\del\a}\Big|_{\a=0}\Big],
    \end{aligned}
  \]
  because, in addition to $\mc L$, $D\mc N(\ol p+z)$ respects modes,
  although for $z\ne0$ the linear wavespeed is changed.  Since $\mc
  S\,\mc L:\mc H_2\to\mc H_+$ is invertible, and $D\mc N(\ol p+z) = \mc I +
  O(z)$, we have both
  \[
  W(0,z) = 0 \com{and} \frac{\del W(0,z)}{\del\a} = 0.
  \]
  Since $W(\a,z)$ is twice differentiable in a neighborhood of $0$,
  the estimates \eqref{West} follow.
\end{proof}

It is remarkable that we can solve the infinite dimensional auxiliary
equation \emph{without requiring} any estimates on the decay rate of
the small divisors, because we have reframed the problem as the
vanishing of a composition of operators, \eqref{Neq}.  Indeed, the
faster the small divisors decay, the smoother the corresponding
periodic solution must be, as seen by the norm \eqref{Hplus}.  By
splitting the Hilbert spaces into orthogonal complements in
\eqref{Hsplit} and \eqref{Hplus}, we have explicitly carried out the
Liapunov-Schmidt decomposition of the nonlinear operator $\mc F$
around 0, as anticipated in \cite{TYperBif}.  It remains only to solve
the \emph{scalar} bifurcation equation.

\subsection{Solution of the bifurcation equation}

To complete the solution of equation \eqref{Neq}, it
remains to show that we can ensure, after use of \eqref{Weq}, that the
remaining component of $\mc F(y^0)$, which is the component orthogonal
to the range of $D_{y_2}\whc F(\ol p,0)$, also vanishes.  From the
decomposition \eqref{Hplus}, this is the (scalar) \emph{bifurcation
  equation},
\begin{equation}
  \label{bif}
  \begin{aligned}
  f(\a,z) &:=
  \Big\langle \s(k\,t\,\piot),
  \whc F\big(\ol p+z+\a\,\c(k\,t\,\piot),W(\a,z)\big)\Big\rangle\\
  &= \Big\langle \s(k\,t\,\piot),
  \mc S\,\mc L\,\mc N\big(\ol p+z+\a\,\c(k\,t\,\piot)
  + W(\a,z)\big)\Big\rangle = 0,
  \end{aligned}
\end{equation}
by \eqref{Fhat}.  Here $\a$ is the amplitude of the linearized
solution, and $z$ is the 0-mode adjustment which can be regarded as
bringing compression and rarefaction back into balance, as described
in \cite{TYdiff2}.  Recall that the period $T$ depends also on this
adjustment, and our choice of $T$ in \eqref{Tn2} amounts to the
introduction of a Hopf bifurcation parameter.  Although \eqref{bif} is
a scalar equation in two variables $\a$ and $z$, the difficulty arises
because the scalar function $f(\a,z)$ factors through the infinite
dimensional solution $W(\a,z)$ of the auxiliary equation.  We will see
presently that the existence of a solution of \eqref{bif} is a
consequence of the genuine nonlinearity of the system, which states
that the nonlinear wavespeed depends monotonically on $z$, which in
turn perturbs the quiet state pressure.

Since $\mc S$ and $\mc L$ are independent of $\a$, differentiation by
$\a$ commutes with both $\mc S$ and $\mc L$.  Similarly, because
$\mc N = \mc L^{-1}\,\mc E$, differentiation of $\mc N$ by $\a$
amounts to differentiation of the nonlinear evolution $\mc E$.  When
applied to $\whc F\big(\ol p+z+\a\,\c(k\,t\,\piot),W(\a,z)\big)$, this
yields
\begin{equation}
  \label{FSLN}
  \frac{\del\whc F}{\del\a} =
  \frac{\del\mc F}{\del\a} =
  \mc S\,\mc L\,\frac{\del\mc N}{\del\a}
  = \mc S\,\frac{\del\mc E}{\del\a}.
\end{equation}

The scalar function $f(\a,z)$ given by \eqref{bif} is defined on the
neighborhood $\mc U\subset\mc H_1$, which with a slight abuse of
notation can be regarded as $\mc U\subset\B R^2$, so we write
\[
  f:\mc U\subset\B R^2\to\B R, \com{with} f(0,0)=0.
\]
As in our description of the bifurcation argument above, we would like
to apply the implicit function theorem to $f$, to get a curve
$z=z(\a)$ on which $f\big(\a,z(\a)\big)=0$.  We cannot apply this
directly, because $\frac{\del f}{\del z}\big|_{(0,0)} = 0$, since all
0-modes are killed by the projection $\mc S$.  Thus we consider the
second derivative $\frac{\del^2f}{\del z\del\a}\big|_{(0,0)}$, and if
this is nonzero, we can conclude the existence of a solution.

One way of effectively calculating the second derivative is to replace
$f$ by the function
\begin{equation}
  \label{gdef}
  \begin{aligned}
    g(\a,z) &:= \frac1\a\, f(\a,z), \quad \a\ne0,\\[2pt]
    g(0,z) &:=  \frac{\del f}{\del\a}(0,z), 
  \end{aligned}
\end{equation}
which is consistently defined because $W(0,z)=0$, and so also
$f(0,z) = 0$, for all $z$ near 0.  It is then clear that
\[
  f(\a,z) = 0 \com{iff} g(\a,z) = 0 \com{for} \a \ne 0,
\]
and it suffices to apply the implicit function theorem to $g$.

To calculate $\del g/\del z$ at $(0,0)$, we first calculate
$\del\mc F/\del\a$ and set $\a=0$.  From \eqref{FSLN}, and using
\eqref{West}, we have
\begin{align}
  \frac{\del\mc F(y)}{\del\a}\Big|_{\a=0}
  &= \mc S\,D\mc E\big(\ol p + z+\a\,\c(k\,t\,\piot)+W(\a,z)\big)
    \Big[\c(k\,t\,\piot) +
    \frac{\del W}{\del\a}\Big]\Big|_{\a=0}\nonumber\\
  &= \mc S\,D\mc E\big(\ol p + z\big)\big[\c(k\,t\,\piot)\big],
    \label{dFda}
\end{align}
and recall that $D\mc N(\ol p) = \mc L$ and $D\mc E(\ol p+z)$ are diagonal,
because $z$ parameterizes the perturbation of the background quiet
state, and all linearized operators around quiet states respect modes,
notwithstanding the dependence of the period $T$ on $z$.

Using \eqref{dFda} in \eqref{bif}, \eqref{gdef}, it follows that
\[
  g(0,z) = \Big\langle \s(k\,t\,\piot),
  \mc S\,D\mc E(\ol p+z)\big[\c(k\,t\,\piot)\big]\Big\rangle,
\]
in which $T$ depends on $z$ according to \eqref{Tn2}, namely
$T = a(\ol p+z)\,\wh T_k$.  Now
differentiating in $z$, we have
\[
  \begin{aligned}
  \frac{\del g}{\del z}\Big|_{(0,0)}
  &= \Big\langle \s(k\,t\,\piot),
  \mc S\,D^2\mc E(\ol p)
    \big[1,\c(k\,t\,\piot)\big]\Big\rangle\\
  &\quad{} + \frac{\del}{\del T}
    \Big\langle \s(k\,t\,\piot),
    \mc S\,D\mc E(\ol p)\big[\c(k\,t\,\piot)\big]\Big\rangle
    \,\frac{\del T}{\del z}\Big|_{z=0},
  \end{aligned}
\]
where the inputs to the bilinear operator $D^2\mc E$ correspond to the
$0$-mode and $k$-mode given by $z$ and $\a\,c(k\,t\,\piot)$, in which
we are differentiating, respectively.  However, since $D\mc E(\ol p)$
is diagonal and the inner product does not depend on the length of the
period but only the Fourier decomposition, it follows that
$\frac{\del}{\del T}\big\langle\cdot,\cdot\big\rangle = 0$, and the
derivative simplifies to
\begin{equation}
  \label{dgdz1}
  \frac{\del g}{\del z}\Big|_{(0,0)}
  = \Big\langle \s(k\,t\,\piot),
  \mc S\,D^2\mc E(\ol p)
    \big[1,\c(k\,t\,\piot)\big]\Big\rangle.
\end{equation}

The final step to complete the bifurcation argument and proof of the
existence of a space and time periodic solution is to show that this
derivative \eqref{dgdz1} is nonzero,
\[
\frac{\del g}{\del z}\Big|_{(0,0)}  \ne 0.
\]
According to \eqref{Shat}, we have
\[
\mc S = \IR-\,\mc T^{-\chi\,T/4},
\]
where as usual, $\chi$ indicates the boundary condition used.  For any
function $Z(t)$, using the adjoints $\IR-^\dag=\IR-$ and ${\mc
  T^\tau}^\dag=\mc T^{-\tau}$, we calculate
\[
\begin{aligned}
  \big\langle \s(k\,t\,\piot),\mc S\,Z(t)\big\rangle
  &= \Big\langle \mc T^{\chi\,T/4}\,\IR-\,\s(k\,t\,\piot),
  Z(t)\Big\rangle\\
  &= \Big\langle \s\big(k\,(t-\chi\,\TS{\frac{T}4})\,\piot\big)
  ,Z(t)\Big\rangle\\
  &= \Big\langle \s\big(k\,t\,\piot-\chi\,k\,\TS{\frac{\pi}2}\big)
  ,Z(t)\Big\rangle,
\end{aligned}
\]
since $\s(\cdot)$ is odd, and so \eqref{dgdz1} reduces to
\begin{equation}
  \label{dgdz}
\frac{\del g}{\del z}\Big|_{(0,0)} =
\frac{\del^2f}{\del z\,\del\a}\Big|_{(0,0)} =
\Big\langle \s\big(k\,t\,\piot-\chi\,k\,{\TS\frac\pi2}\big),
D^2\mc E(\ol p)
\big[1,\c(k\,t\,\piot)\big]\Big\rangle.
\end{equation}

The final estimate we need is contained in the following technical
lemma, which is also proved in Section~\ref{sec:D2E} below.  The lemma
is essentially the statement that genuine nonlinearity, which implies
that the average wavespeed varies with background quiet state, implies
non-degeneracy of the bifurcation equation.

\begin{lemma}
  \label{lem:D2Enz}
  For a genuinely nonlinear constitutive equation, the second
  derivative $D^2\mc E(\ol p)$ satisfies
  \[
    \Big\langle \s\big(k\,t\,\piot-\chi\,k\,{\TS\frac\pi2}\big),
    D^2\mc E(\ol p)
    \big[1,\c(k\,t\,\piot)\big]\Big\rangle\ne 0
  \]
  so that the derivative
  \begin{equation}
  \label{dgzdfza}
  \frac{\del g}{\del z}\Big|_{(0,0)} = 
  \frac{\del^2f}{\del z\,\del\a}\Big|_{(0,0)} \ne 0.
  \end{equation}
\end{lemma}

The proof of Lemma \ref{lem:D2Enz} relies on a detailed analysis of
the nonlinear evolution of the system and is given in
Lemma~\ref{lem:duh} below.  Assuming it for now, we complete the
analysis of the bifurcation equation, which implies that nonresonant
$k$-mode solutions of the linearized equation perturb to nonlinear
pure tone solutions.

\begin{theorem}
  \label{thm:bifurc}
  Assume that the $k$-mode is nonresonant.  Then there exists a
  one-parameter family of solutions of the form \eqref{pert} of
  equation \eqref{yper} or \eqref{yab}, respectively, parameterized by
  the amplitude $\a$ in a neighborhood of $0$.  This in turn generates
  a periodic pure tone solution of the compressible Euler equations.
\end{theorem}

\begin{proof}
  Lemma~\ref{lem:himodes} shows that the auxiliary equation has a
  unique solution in a neighborhood of the origin.  It remains only to
  show that the bifurcation equation \eqref{bif}, namely
  $f\big(\a,z\big) = 0$, for the $k$-mode can always be solved
  uniquely in a neighborhood of the origin.  Using \eqref{dgzdfza} the
  implicit function theorem implies the existence of a unique $z(\a)$
  such that \eqref{y0} gives the data $y^0$ which solves the
  equation \eqref{bif} for $\a$ in a neighborhood of the origin.
\end{proof}

The proof of Theorem \ref{thm:bifurc} will be complete once we prove
Theorems~\ref{thm:D2E} and \ref{thm:sl}, as well as
Lemma~\ref{lem:D2Enz}.  In the next section we carry out a detailed
analysis of the SL system by introducing an angle variable, and prove
Theorem~\ref{thm:sl}.  In the following section we differentiate the
evolution operator and prove Theorem~\ref{thm:D2E} and
Lemma~\ref{lem:D2Enz}.  After that we further study the resonance
structure and give a sense in which the fully nonresonant entropy
profiles are generic.

\section{Analysis of Sturm-Liouville Operator}
\label{sec:SL}

To proceed we need a detailed analysis of the SL operator, which
includes a complete description of both the fundamental solution and
the full set of eigenfrequencies.  Once again we treat the general
entropy profile as given, and solve for the base frequencies $\w_k$
that yield $k$-mode solutions of the linearized equation.  Recall that
starting with the linearization \eqref{UP}, we separated variables by
setting
\[
  P(x,t) := \vp(x)\,\c(\w\,t), \qquad
  U(x,t) := \psi(x)\,\s(\w\,t),
\]
which yields the SL system \eqref{SL1}, namely
\begin{equation}
  \label{SL3}
  \dot\vp + \w\,\psi = 0, \qquad
  \dot\psi - \w\,\sg^2\,\vp = 0,
\end{equation}
with $\sg^2 = - v_p(\ol p,s) = \sg^2(x)$.  According to \eqref{bc0},
we first consider initial values
\[
  \vp(0) = c_0 \com{and} \psi(0) = 0,
\]
for appropriate $c_0\ne 0$.  Those values of $\w_k$ that meet the
boundary conditions \eqref{bcab} or \eqref{bcper} will then determine
the eigenfrequencies, which in turn define the reference period
$T_k:=k\,\frac{2\pi}{\w_k}$.

\subsection{Angle Variable}

For fixed entropy profile $s(x)$, the SL evolution is efficiently
analyzed using the angle $\t=\t(x)=\t(x,\w)$ of the vector
$(\vp,\psi)$ which satisfies \eqref{SL3}.  This in turn can be
effectively captured with the use of \emph{modified Pr\"ufer
variables}, which are
\begin{equation}
  \label{pruf}
  \vp(x) := r(x)\,\frac1{\varepsilon(x)}\,\c\big(\t(x)\big), \qquad
  \psi(x) := r(x)\,\varepsilon(x)\,\s\big(\t(x)\big),
\end{equation}
see \cite{Pryce,CodLev,BR}.  Here we interpret $r(x)$ as the radial
length or amplitude, $\varepsilon(x)$ as the eccentricity, and $\t(x)$ as the
angle variable.  This is a degenerate description in which we are free
to choose the eccentricity $\varepsilon(x)>0$, and having done so, both $r(x)$
and $\t(x)$ will be determined by the equations.

We solve equation \eqref{SL3} on an interval $[x_I,x_E]$ with the
assumption that the entropy $s$ (and so also $\sg$) is $H^b$ on this
interval.  Plugging in \eqref{pruf} into \eqref{SL3} and simplifying,
we get the system
\[
  \begin{aligned}
    \frac{\dot r}r\,\c(\t) - \frac{\dot \varepsilon}{\varepsilon}\,\c(\t)
    - \s(\t)\,\dot\t + \w\,\varepsilon^2\,\s(\t) &= 0, \\
    \frac{\dot r}r\,\s(\t) + \frac{\dot \varepsilon}{\varepsilon}\,\s(\t)
    + \c(\t)\,\dot\t - \w\,\frac{\sg^2}{\varepsilon^2}\,\c(\t) &= 0,
  \end{aligned}
\]
with initial conditions
\[
  \t(x_I) = \t_I \com{and} r(x_I) = r_I,
\]
to be determined by the initial values $\vp(x_I)$ and $\psi(x_I)$.
After use of elementary trig identities, this in turn becomes
\[
  \begin{aligned}
    \frac{\dot r}r - \frac{\dot \varepsilon}{\varepsilon}\,\c(2\t)
    + \w\,\Big(\varepsilon^2-\frac{\sg^2}{\varepsilon^2}\Big)\s(\t)\,\c(\t) &=0, \\
    \dot\t + \frac{\dot \varepsilon}{\varepsilon}\,\s(2\t)
    - \w\,\Big(\frac{\sg^2}{\varepsilon^2}\,\c^2(\t) + \varepsilon^2\,\s^2(\t)\Big) &=0.
  \end{aligned}
\]
Thus, once $\varepsilon(x)$ is chosen, the system is reduced to a nonlinear
scalar ODE for the angle $\t(x)$, coupled with an integration for the
amplitude $r(x)$.  Moreover, it is now clear that we should choose
$\varepsilon$ such that
\begin{equation}
  \label{adef}
  \varepsilon^2 = \frac{\sg^2}{\varepsilon^2}, \com{that is}
  \varepsilon (x) := \sqrt{\sg(x)}.
\end{equation}
With this choice the equations simplify further, and we get the
nonlinear scalar equation for $\t(x)$,
\begin{equation}
  \label{theta}
  \dot\t = \w\,\sg - \frac{\dot\sg}{2\sg}\,\s(2\t),
  \qquad \t(x_I) = \t_I,
\end{equation}
coupled with a linear homogeneous equation for $r(x)$, namely
\[
  \frac{\dot r}r = \frac{\dot\sg}{2\sg}\,\c(2\t),
  \qquad r(x_I) = r_I,
\]
which immediately yields the quadrature
\begin{equation}
  \label{rint}
  r(x) = r_I\,\exp\Big\{ \int_{x_I}^x
  \c\big(2\t(y)\big)\;d\log\sqrt\sg(y)\Big\}.
\end{equation}

We can now describe the fundamental solution $\Psi(x-x_I;\w)$ of the
linear system on the interval $x\in[x_I,x_E]$ by meeting the
appropriate initial conditions.  The idea is that the evolution is
represented by a rotation in the angle variable, so can be described
by first changing to that variable, rotating the appropriate amount,
and then reverting to the original variables.

\begin{lemma}
  \label{lem:SLfund}
  A fundamental matrix of the system \eqref{Psi} for the evolution of
  $\vp(x)$ and $\psi(x)$, for $x\in[x_I,x_E]$, for which the entropy
  is $H^b$, is
  \begin{equation}
    \label{SLfund}
    \Psi(x-x_I;\w) = \frac{r(x)}{r_I}\,M\big(\sg(x_-)\big)
    R\big(\t(x)-\t_I\big)\,M^{-1}\big(\sg({x_I}_+)\big),
  \end{equation}
  where $M(\cdot)$ is defined by
  \begin{equation}
    \label{MJmat}
    M(q) := \( 1/\sqrt{q}&0\\0&\sqrt{q}\), \qquad M^{-1}(q)=M(1/q),
  \end{equation}
  $R(\t)$ is the usual \x2 rotation matrix, and $\t(x)$ and $r(x)$ satisfy
  \eqref{theta} and \eqref{rint}, respectively.
\end{lemma}

In this representation of the fundamental solution, the first column
describes the linearized evolution of the even modes $\c(\w\,t)$, and
the second column describes evolution of the odd modes $\s(\w\,t)$ in
the linearized PDE system \eqref{UP}.

\begin{proof}
  Writing the substitution \eqref{pruf} in matrix notation and using
  \eqref{adef}, we have
  \[
    \(\vp(x)\\\psi(x)\) =
    r(x)\,M\big(\sg(x)\big)\,R\big(\t(x)\big)\,\(1\\0\),
  \]
  for any $x\in[x_I,x_E]$.  It follows that
  \[
  \(1\\0\) = R(-\t_I)\,M\big(1/\sg(x_I)\big)\,\frac1{r_I}\,
  \(\vp(x_I)\\\psi(x_I)\),
  \]
  and substituting this back in and comparing to the defining relation
  \[
  \(\vp(x)\\\psi(x)\) = \Psi(x-x_I;\w)\,\(\vp(x_I)\\\psi(x_I)\)
  \]
  completes the proof.
\end{proof}

In SL theory, the fundamental solution is also referred to as a
\emph{transfer matrix}~\cite{Pryce}.

\subsection{Jump Condition}

Recall that we restrict to piecewise $H^b$ entropy profiles.  The
representation \eqref{pruf} of solutions implicitly assumes that $\sg$
is absolutely continuous, so does not apply at jump discontinuities.
At an entropy jump the ODE \eqref{theta} formally degenerates to
$\dot\t = c\,\d$, a Dirac measure.  Instead of looking for weak
solutions, at entropy jumps we refer back to the original SL system
\eqref{SL3} and impose continuity of both $\vp$ and $\psi$.  This is
consistent with taking $p$ and $u$ continuous at contact
discontinuities in the Euler equations \eqref{genls}.

At an entropy jump, we thus update $\t$ and $r$ by directly using
\eqref{pruf} and imposing continuity on both $\vp$ and $\psi$.  Thus,
assume that the entropy is discontinuous at $x$ and has left and right
limits $s(x_-)$ and $s(x_+)$, repectively.  Using \eqref{pruf} and
imposing continuity, we get the equations
\begin{equation}
  \label{jumpmp}
  \begin{aligned}
    r_-\,\frac1{\sqrt{\sg_-}}\,\c(\t_-)
    &= r_+\,\frac1{\sqrt{\sg_+}}\,\c(\t_+),\\
    r_-\,\sqrt{\sg_-}\,\s(\t_-)
    &= r_+\,\sqrt{\sg_+}\,\s(\t_+),
  \end{aligned}
\end{equation}
where again the subscripts $\mp$ denote left and right limits at $x$,
respectively.  Dividing the equations, we get
\begin{equation}
  \label{thjump}
  \sg_-\,\tan(\t_-) = \sg_+\,\tan(\t_+),
\end{equation}
while we can also rewrite them as
\begin{equation}
  \label{rjump}
  \frac{r_+}{r_-}\,\c(\t_+) =
  \frac{\sqrt{\sg_+}}{\sqrt{\sg_-}}\,\c(\t_-),\qquad
  \frac{r_+}{r_-}\,\s(\t_+) =
  \frac{\sqrt{\sg_-}}{\sqrt{\sg_+}}\,\s(\t_-).
\end{equation}

To proceed we now define the \emph{jump} $J$ at the discontinuity at
$x$ by
\begin{equation}
  \label{Jdef}
  J := \frac{\sg_-}{\sg_+} = \frac{\sg(s(x_-))}{\sg(s(x_+))}.
\end{equation}
We define the function
\begin{equation}
  \label{hdef}
  h(J,z) :=
  \begin{cases}
    \text{Arctan}(J\,\tan z\big) + m\,\pi, &-\pi/2<z-m\,\pi<\pi/2,\\
    z &z = m\,\pi \pm\pi/2,
  \end{cases}
\end{equation}
where $\text{Arctan}(\square)\in(-\pi/2,\pi/2)$ is the principal
branch.  Geometrically, $h$ gives the angle that results from a vector
with angle $z$ scaling according to \eqref{thjump}; in particular,
$h(J,z)$ is always in the same quadrant as $z$, and $h(J,z)=z$ for $z$
on the coordinate axes.

\begin{lemma}
  \label{lem:h}
  For any fixed $J$, $h(J,\cdot)$ is a monotone real analytic function
  of $z$, and for any fixed $z$, $h(\cdot,z)$ is a real analytic,
  bounded function of $J$.  The partial derivatives of $h$,
  \begin{equation}
    \label{hder}
    \frac{\del h}{\del z} =
    \frac{J\,(1+\tf^2(z))}{1 + J^2\,\tf^2(z)}, \qquad
    \frac{\del h}{\del J} = \frac{\tf(z)}{1 + J^2\,\tf^2(z)},
  \end{equation}
  where we have abbreviated $\tan(\cdot)=:\tf(\cdot)$, are
  uniformly bounded,
  \[
    \min\Big\{\frac1J,J\Big\}\le \frac{\del h}{\del z} \le
    \max\Big\{\frac1J,J\Big\}, \qquad
    \Big|\frac{\del h}{\del J}\Big| \le\frac1{2\,J^2}.
  \]
\end{lemma}

\begin{proof}
  The derivatives \eqref{hder} are directly calculated for
  $z\ne(m\pm\frac12)\pi$, at which points $\tf(z)=\pm\infty$.  The
  function is defined in \eqref{hdef} to be continuous at those
  points, and the derivatives are bounded there,
  \[
    \lim_{\tf(z)\to\pm\infty}
    \frac{J\,(1+\tf^2(z))}{1 + J^2\,\tf^2(z)} = \frac1J, \qquad
    \lim_{\tf(z)\to\pm\infty}
    \frac{\tf(z)}{1 + J^2\,\tf^2(z)} = 0.
  \]
  Analyticity and boundedness now follow because the denominators in
  \eqref{hder} never vanish.  The derivative bounds follow by
  elementary Calculus.
\end{proof}

We are now in a position to fully describe the effect of a finite
entropy jump.

\begin{lemma}
  \label{lem:jump}
  In Pr\"ufer coordinates, the angle
  changes across an entropy jump as
  \begin{equation}
    \label{thJ}
    \t_+ = h(J,\t_-),
  \end{equation}
  while the radius changes as
  \begin{equation}
    \label{rJ}
    \begin{aligned}
      r_+ &= r_- \, \rho(J,\t_-), \com{where}\\
      \rho(J,\t_-) &:= \sqrt{\frac1J\,\c^2(\t_-) + J\,\s^2(\t_-)},
    \end{aligned}
  \end{equation}
  this being uniformly bounded,
  \[
    \min\big\{\sqrt J, 1/\sqrt J\big\}
    \le \rho(J,\cdot) \le \max\big\{\sqrt J, 1/\sqrt J\big\}.
  \]
\end{lemma}

\begin{proof}
  Using the definitions \eqref{adef} and \eqref{MJmat}, we write
  \eqref{pruf} as
  \begin{equation}
    \label{prufpp}
    \(\vp(x)\\\psi(x)\) = r(x)\,M\big(\sg(x)\big)
    \(\c\big(\t(x)\big)\\\s\big(\t(x)\big)\),
  \end{equation}
  so continuity of $\vp$ and $\psi$ gives \eqref{jumpmp},
  which is in turn equivalent to \eqref{thjump}, \eqref{rjump}.  Using
  \eqref{Jdef}, \eqref{hdef} in \eqref{thjump} yields \eqref{thJ},
  while eliminating $\t_+$ in \eqref{rjump} yields \eqref{rJ}.  The
  bound on $\rho(J,\cdot)$ is immediate.
\end{proof}

For piecewise $H^b$ entropy profiles, we solve the ODE \eqref{theta}
on intervals of smooth entropy (with appropriate initial condition),
and couple these evolutions with \eqref{thJ} at entropy jumps, to get
a uniquely defined piecewise $H^b$ angle $\t$ throughout the interval
$[0,\ell]$.  Use of \eqref{rint} on subintervals, together with
\eqref{rjump} at jumps, explicitly determines $r(x)$, and use of
\eqref{pruf} finally gives $\vp$ and $\psi$ throughout the interval
$[0,\ell]$.  Thus the angle variable, which solves the scalar ODE
\eqref{theta}, determines the entire SL evolution.

To be specific, let the entropy field $s\in\mc P$ be piecewise $H^b$
on $[0,\ell]$ with jumps at the points
\[
0<x_1<x_2<\dots<x_N<\ell,
\]
and set $x_0:=0$ and $x_{N+1}:=\ell$.  Define the jump $J_j$ at $x_j$ by
\eqref{Jdef}, namely
\begin{equation}
  \label{Jj}
  J_j = \frac{\sg(s(x_{j-})}{\sg(s(x_{j+})},
\end{equation}
and for given $\w>0$ define the angle $\t(x)$ inductively as follows:
first set $\t_{0+}=\t_0$, which we assume given.  On the interval
$(x_j,x_{j+1})$, let $\t(x)$ be the solution of \eqref{theta} with
initial data $x_I:=x_j$ and $\t_I:=\t_{j+}$: this determines $\t$ up
to $x_{j+1-}$.  At the $j$-th jump $x_j$, we use \eqref{thJ}, that is
we set
\[
\t_{j+} := h(J_j,\t_{j-}).
\]
Inductively, this determines the angle $\t(x)$ on the entire interval
$[0,\ell]$.

\begin{lemma}
  \label{lem:angle}
  This angle function $\t(x)$ as defined above, together with the
  entropy profile $s(x)$, fully determine the solution to the SL
  initial value problem \eqref{SL3}, with arbitrary initial condition
  $(\vp_0,\psi_0)$.  Moreover, the radius variable $r(x)$ is uniformly
  bounded: that is, there exist finite positive constants $\ul r$ and
  $\ol r$ depending only on $\sg(x)$, such that
  \[
    r(x)\in [\ul r,\ol r] \com{for all} x\in[0,\ell].
  \]
\end{lemma}

\begin{proof}
  First, use \eqref{pruf} to determine the initial values $\t_0$ and
  $r_0$, namely
  \[
    \tan\t_0 := \frac{\psi_0}{\sg_0\,\vp_0}, \qquad r_0 :=
    \sqrt{\sg_0\,\vp_0^2 + \psi_0^2/\sg_0},
  \]
  where $\sg_0 := \sg(s(0))$.  The algorithm described above then
  determines $\t(x)$ for $x\in[0,\ell]$.  Now define $r(x)$
  inductively using \eqref{rint} and \eqref{rjump}: this becomes
  explicit once $\t(x)$ is regarded as known.  Now use \eqref{pruf} to
  determine $\vp(x)$ and $\psi(x)$: by construction, these satisfy
  \eqref{SL3} on each subinterval $(x_j,x_{j+1})$, and are continuous
  at each $x_j$.
  
  It remains to find bounds for $r$: on an interval in which $s$ is
  continuous, this is provided by \eqref{rint}.  On the other hand,
  at an entropy jump, by \eqref{rJ} we have
  \begin{equation}
    \label{rJ1}
    \min\big\{1/\sqrt J,\sqrt J\big\} \le \frac{r_+}{r_-}
    \le \max\big\{1/\sqrt J,\sqrt J\big\},
  \end{equation}
  and combining finitely many such changes in $r$ yields the uniform
  lower and upper bounds $\ul r$ and $\ol r$, respectively.  
\end{proof}

\subsection{Eigenfrequencies}

We are interested in characterizing the $k$-mode eigenfrequencies and
corresponding periods, that yield periodic solutions of the linearized
equations.  That is, we seek to characterize periods
$T_k=k\,\frac{2\pi}{\w_k}$, such that the $k$-mode basis vector $\mc
B_k(T_k)$, defined in \eqref{Tk}, satisfies the boundary conditions
\eqref{bcab} or \eqref{bcper}, respectively, or equivalently
\[
  \d_{k}(T_k)=0,
\]
respectively, these being given by \eqref{dkcts}.  Using \eqref{pruf},
we express the boundary conditions in terms of the angle variable
$\t(x)$, as follows.  For the periodic boundary condition
\eqref{bcper}, with $\chi=1$, we require
\begin{equation}
  \label{kodd}
  \vp_k(\ell) = 0, \com{so} \t(\ell) = \Big(n+\frac12\Big)\pi,
\end{equation}
for $k$ odd, and for $k$ even, or the acoustic boundary condition
$\chi=0$, we need
\begin{equation}
  \label{keven}
  \psi_k(\ell) = 0, \com{so} \t(\ell) = n\,\pi,
\end{equation}
for some $n$.  Recall that it suffices to consider only the periodic
boundary condition because the even modes of that problem will
automatically solve the acoustic boundary condition $u=0$ at $x=\ell$
by our reflection principle.

We can methodically enumerate both SL boundary conditions \eqref{kodd}
and \eqref{keven} in terms of the frequency $\w$ and angle
$\t(x) = \t(x,\w)$, using the \emph{angle boundary condition} for the
$k$-th mode,
\begin{equation}
  \label{abc}
  \t(\ell,\w) = k\,\frac\pi2,
\end{equation}
which we interpret as an implicit condition for the coefficient
$\w=\w_k$ of \eqref{theta}.

We express \eqref{abc} more explicity as follows.  Recall that a
piecewise continuous function on a closed interval is continuous at
the boundaries, and is bounded on $[0,\ell]$.  For definiteness,
suppose that the entropy profile has jumps $J_j$ defined by \eqref{Jj}
at points
\[
0<x_1<x_2<\dots<x_N<\ell,
\]
and is $H^b$ on each subinterval.  Denote the left and right limits of
$\t$ at the entropy jump $x^j$ by $\t_{j-}$ and $\t_{j+}$,
respectively.  For each $\w$, denote the nonlinear evolution of the
angle across the subinterval $(x_j,x_{j+1})$ by
\[
\Phi(\t_{j+},\w,x_j,x_{j+1}) := \t(x_{j+1-}),
\]
where $\t(x)$ solves \eqref{theta} in $(x_j,x_{j+1})$ with $x_I=x_j$
and $\t_I = \t_{j_+}$.  According to our description above, we have,
inductively, $\t_{0+}=0$,
\begin{equation}
  \label{Phind}
  \t_{j-} = \Phi(\t_{j-1+},\w,x_{j-1},x_j), \qquad
  \t_{j_+} = h(J_j,\t_{j-}).
\end{equation}
Our angle boundary condition \eqref{abc} then becomes
\begin{equation}
  \label{Phinl}
  \t(\ell,\w) = \t_{N+1-} =  \Phi(\t_{N+},\w,x_N,\ell) = k\,\frac\pi2,
\end{equation}
which is an implicit equation for $\w_k$.

According to \eqref{theta}, the angle evolution across a subinterval
is given implicitly by
\begin{equation}
  \label{Phi}
  \begin{aligned}
  \t_{j-} &= \Phi(\t_{j-1+},\w,x_{j-1},x_j)\\
  &= \t_{j-1+} + \int_{x_{j-1}}^{x_j} \w\,\sg\;dx
  + \int_{x_{j-1}}^{x_j} s\big(2\t(x)\big)\;d\log\sqrt\sg(x).
  \end{aligned}
\end{equation}

Thus, given an entropy profile, we find the base reference period
$T_k$ of the $k$-mode by solving \eqref{Phinl} for $\w_k$, where
$\t(x,\w)$ solves \eqref{theta}, and then using \eqref{Tn}.  Each such
$k$-mode with reference period $T_k$ then determines a periodic
solution of the linearized equation.

Note that \eqref{abc} is equivalent to the \emph{periodic} boundary
condition \eqref{bcper}, but permits more frequencies than the
\emph{acoustic} reflection boundary condition \eqref{bcab} allows,
these latter two being equivalent only for even values of $k$.
However, the even modes form a closed subspace, so if we start from
$a_m=0$ for all odd $m$, this persists when we perturb to the
nonlinear problem.  Thus for notational convenience, we use
\eqref{Phinl} to identify all frequencies, with the understanding that
if we are using \eqref{bcab}, all linearizations and perturbations are
restricted to even modes.

\begin{lemma}
  \label{lem:wk}
  For each integer $k\ge 1$, and any piecewise $H^b$ entropy profile
  $s = s(x)\in\mc P$, there is a unique $\w_k$ such that $\t(x,\w_k)$
  as defined above satisfies \eqref{abc}.  Moreover, $\{\w_k\}$ is a
  strictly monotone increasing sequence,
  \[
  \w_{k_1} < \w_{k_2} \com{for} k_1 < k_2,
  \]
  and the eigenfrequencies $\w_k$ grow like $k$, that is
  \begin{equation}
    \label{wkok}
    \lim_{k\to\infty}\frac{\w_k}k = \Lambda, \com{where}
    \Lambda := \frac{\pi/2}{\int_0^\ell \sg(y)\;dy}.
  \end{equation}
\end{lemma}

\begin{proof}
  We show that for any fixed $x>0$, the function $\t(x,\w)$ is a
  continuous, strictly monotone increasing function of $\w$.  In
  particular, taking $x=\ell$ in \eqref{abc} implies that $\w_k$
  exists for each $k\ge1$ and is increasing.

  We first consider an interval $(x_I,x_E)$ on which $s(x)$ is $H^b$.
  Note that $\sg$ is determined by $s(x)$, and differentiate the ODE
  \eqref{theta} with respect to $\w$.  Denoting $\frac{\del\t}{\del
    \w}$ by $\zeta$, this yields
  \begin{equation}
    \label{zeta}
    \dot\zeta = \sg - \frac{\dot\sg}{\sg}\,\c(2\t)\,\zeta,
  \end{equation}
  which by \eqref{rint} can also be written
  \[
    \dot\zeta = \sg - 2\,\frac{\dot r}{r}\,\zeta.
  \]
  It follows that $r^2$ is an integrating factor, and we integrate to
  get
  \begin{equation}
    \label{zetasol}
  r^2(x)\,\zeta(x) = r^2(x_I)\,\zeta(x_I)
  + \int_{x_I}^xr(y)^2\,\sg(y)\;dy.
  \end{equation}
  Note that our left boundary condition $\t(0)=0$ implies
  $\zeta(0)=0$.  Thus $\zeta(x)>0$ for $x>x_I$ provided
  $\zeta(x_I)\ge 0$, and in particular $\zeta(x_E)>0$.

  On the other hand, at the entropy jump $x_j$, assume that
  \[
  \t(x_{j-},\w_\flat) \le \t(x_{j-},\w_\sharp)
  \com{for any} 0\le \w_\flat \le \w_\sharp.
  \]
  Since $\t(x_{j+},\w)$ is given by \eqref{thJ}, continuity and
  monotonicity of $h(J,\cdot)$ for any fixed $J>0$ implies that
  \[
  \t(x_{j+},\w_\flat) \le \t(x_{j+},\w_\sharp).
  \]
  It now follows by induction on subintervals that $\t(\ell)$ is a
  strictly monotone function of $\w$ for every piecewise $H^b$ entropy
  profile.  This in turn implies the existence, uniqueness and
  monotonicity of the eigenfrequencies $\w_k$.

  To get the growth rate of $\w_k$ with respect to $k$, it suffices to
  obtain global upper and lower bounds for $\zeta=\frac{\del\t}{\del
    \w}$, uniform for $x\in[0,\ell]$.  On subintervals where $s$ is
  $H^b$, \eqref{zetasol} immediately gives the bound as long as $r(x)$
  stays positive and finite; this in turn follows immediately from
  \eqref{rint}.  At entropy jumps, differentiating \eqref{thJ}, we
  have
  \[
    \zeta(x_{j+},\w) = \frac{\del\t_{j_+}}{\del \w} = \frac{\del
      h(J,z)}{\del z}\Big|_{z=\t_{j-}}\,\zeta(x_{j-},\w),
  \]
  so we need a uniform bound for $\frac{\del h(J,z)}{\del z}$.
  Referring to \eqref{hder} and rearranging, we have
  \[
    \frac{\del h}{\del z} = \frac1{1+J^2\,\tf^2}\,J
    + \frac{J^2\,\tf^2}{1+J^2\,\tf^2}\,\frac1J,
  \]
  which is a convex combination of $J$ and $1/J$, so that
  \[
    \min\big\{J,1/J\big\} \le
    \frac{\zeta(x_{j+},\w)}{\zeta(x_{j-},\w)} \le
    \max\big\{J,1/J\big\},
  \]
  which is the desired bound.
  
  Because there are finitely many jumps, we thus obtain uniform lower
  and upper bounds for $\zeta(x,\w)$ for any $x$, uniform in
  $\w\ge\epsilon>0$, and in particular for $x=\ell$.  This shows that
  $\t$ grows positively and linearly with $\w$, and since $\t(\ell)$
  grows as $k$, so does $\w_k$, so that $\w_k/k\to \Lambda$.

  We can evaluate the limit $\Lambda$ by formally taking the limit
  $k\to\infty$: as $k\to\infty$, so does $\t(x)$, and the integral in
  \eqref{rint}, which is oscillatory, tends to zero.  Thus the
  limiting value of $r(x)$ on the interval $(x_{j-1},x_j)$ is a
  constant,
  \[
  \lim_{k\to\infty}r(x) = r^\infty_{j-} = r^\infty_{j-1+},
  \]
  where the superscript denotes the limiting value, and 
  \eqref{zetasol} gives
  \[
  \zeta^\infty(x_{j-}) = \zeta^\infty(x_{j-1+}) +
  \int_{x_{j-1}}^{x_j} \sg(y)\;dy.
  \]
  For large values of $\w$ (or $k$), we have
  \[
  \t(x,\w) = \w\,\Big(\zeta^\infty(x) +
  \frac1\w\,\int_0^\w\big(\zeta(x,w)-\zeta^\infty(x)\big)\;dw\Big),
  \]
  which yields
  \[
  \lim_{\w\to\infty}\frac{\t(x,\w)}{\w} = \zeta^\infty(x).
  \]
  On the other hand, \eqref{hdef}, \eqref{thJ} imply that at any jump,
  \[
  |\t_{j+} - \t_{j-}|<\frac\pi2, \com{so that}
  \lim_{\w\to\infty}\Big|\frac{\t_{j+}}{\w}-\frac{\t_{j-}}{\w}\Big|=0,
  \]
  and we can take $\zeta^\infty$ to be continuous at the jumps.
  Combining evolutions and jumps by induction, we conclude that
  \[
  \lim_{\w\to\infty}\frac{\t(\ell,\w)}{\w} = \zeta^\infty(\ell)
  = \int_0^\ell \sg(y)\;dy.
  \]
  Finally, since $\t(\ell,\w_k)=k\,\frac\pi2$, we conclude that
  \[
  \lim_{k\to\infty}\frac{k}{\w_k} =
  \frac2\pi\,\int_0^\ell \sg(y)\;dy,
  \]
  from which our formula for $\Lambda$ follows.
\end{proof}

\section{Differentiation of the Evolution Operator}
\label{sec:D2E}

In order to complete the bifurcation argument, we must differentiate
the evolution operator twice.  We do so for the full \x2 system
\eqref{genls}, with the assumption that $s=s(x)$ is known and
piecewise $H^b$, and that the derivatives $p_t$ and $u_t$ remain
bounded, $(p_t,u_t)\in L^\infty$, so that $p$ and $u$ remain in $H^b$.
We work in the context of the classical local existence theory as
described in Corollary~\ref{cor:le}, so that the uniform estimates
stated there may be assumed.  In this section we do not assume the
symmetries \eqref{tsymm} or \eqref{xsymm}, although they are preserved
under differentiation.  Because $p$ and $u$ are taken to be continuous
at entropy jumps, the case of piecewise $H^b$ entropy profiles is
handled as a finite composition of non-commuting evolutions over
sub-intervals on which the entropy is $H^b$, as above.

We begin with the system \eqref{genls} in conservative form, namely
\begin{equation}
    \label{nleq}
  p_x + u_t = 0, \qquad
  u_x - v(p,s)_t = 0,
\end{equation}
with $s=s(x) \in H^b([0,L])$.  Regarding $L$ as fixed, we consider the
evolution operator from $x=0$ to $x=L$, defined in \eqref{Edef}, namely
\begin{equation}
    \label{EL}
  \mc E^L:U\subset H^b\to H^b, \qquad \mc E^L(y^0) := y(L,\cdot),
\end{equation}
where $y(x,t)=(p,u)$ is the solution of \eqref{nleq}, defined
throughout the interval $[0,L]$.  Here, we are not assuming the
symmetries \eqref{xsymm}, so \eqref{EL} cannot be written as a scalar
equation in $y$.  Instead, with a slight abuse of notation used in
this section, we use $y(x,t)$ to denote the \emph{vector} $(p,u)$, so
the data is $y^0=(p^0,u^0)$, and $H^b$ is the space of vector
functions $(p,u)\in H^b\times H^b$.  The classical local existence
theory, Theorem~\ref{thm:le} and Corollary~\ref{cor:le}, then states
that for $b>3/2$, the operator $\mc E^L$ indeed defines a solution in
$H^b$, provided the data $y^0=(p^0,u^0)$ or interval $[0,L]$ is chosen
small enough that the derivatives $p_t$ and $u_t$ remain bounded
throughout the interval $x\in[0,L]$.

The local existence theorem, Theorem~\ref{thm:le}, is based on $H^b$
energy estimates~\cite{Majda,Taylor}.  In general the inhomogeneous
linearized evolution equation can be written in the abstract form
\begin{equation}
  \label{LXF}
  \mc L\,X = F,\qquad X(0,\cdot) = X^0,
\end{equation}
where $\mc L$ is the linearized evolution operator and $F$ represents
linear and/or nonlinear corrections.  The corresponding energy
estimate takes the form
\begin{equation}
  \label{Hbenergy}
  \|X(x,\cdot)\|_b \le K \,\|X^0\|_b + \int_0^x
  \|F(\chi,\cdot)\|_b\;d\chi,
\end{equation}
where we are evolving in $x$, and this estimate holds as long as the
nonlinear solution remains $H^b$.  This can be extended to establish
the Frechet differentiability of the nonlinear evolution operator
$\mc E^L$, which is worked out in detail in~\cite{Ydiff}, as follows.

Evaluating $\mc E^L(y^0)$ implicitly requires knowledge of the entire
evolution, that is $\mc E^x(y^0)$ for $x\in[0,L]$, and the
linearization will similarly implicitly require knowledge of the
entire linear and nonlinear solution across $[0,L]$, even though only
the initial data $(p^0,u^0)$ is perturbed.  Recall that $\mc E^L$ is
Frechet differentiable at $y^0$, with derivative
\[
  D\mc E^L(y^0)[\cdot]: H^{b} \to H^{b-1},
\]
if there exists a linear and bounded operator, denoted
$D\mc E^L(y^0)$, such that
\begin{equation}
  \label{DExdef}
  \mc E^L(y^0+Y^0) - \mc E^L(y^0) - D\mc E^L(y^0)[Y^0] =
  o(\|Y^0\|_{b-1}),
\end{equation}
that is such that
\[
  \frac1{\|Y^0\|_{b-1}}\,
  \Big\|\mc E^L(y^0+Y^0) - \mc E^L(y^0) - D\mc E^L(y^0)[Y^0]\Big\|_{b-1}
  \to 0,
\]
as $\|Y^0\|_{b-1}\to 0$.  We will use the convention that square
brackets and upper case $[Y^0]$ denote inputs to linear or
multi-linear operators, while parentheses and lower case $(y^0)$ are
inputs to nonlinear functionals.

We briefly summarize how to linearize the evolution operator; again,
see~\cite{Ydiff} for details.  Consider the nonlinear system for a
(small) perturbation $(p+\delta p,u+\delta u)$ of $(p,u)$, namely
\[
  \big(p+\delta p\big)_x + \big(u+\delta u\big)_t = 0, \qquad
  \big(u+\delta u\big)_x - v\big(p+\delta p,s\big)_t = 0,
\]
so that $(\delta p,\delta u)$ satisfies
\begin{equation}
  \label{delpu}
  \begin{gathered}
    \delta p_x + \delta u_t = 0, \qquad
    \delta u_x - \delta v_t = 0,\\
    \delta v := v\big(p+\delta p,s\big) - v(p,s).
  \end{gathered}
\end{equation}
The linearization (or Frechet derivative) $D\mc E^L$ at
$y^0=(p^0,u^0)$, evaluated on data $Y^0=(P^0,U^0)$, is now obtained by
expanding $\delta v$ and retaining the linear terms in \eqref{delpu},
while replacing $(\delta p,\delta u)$ by $(P,U)$.  This yields the
linear system
\begin{equation}
  \label{lineq}
  P_x + U_t = 0, \qquad
  U_x - \Big(\frac{\del v(p,s)}{\del p}\,P\Big)_t = 0,
\end{equation}
with data $Y^0=(P^0,U^0)$, and we have
\begin{equation}
  \label{DE}
  D\mc E^L(y^0)[Y^0] = Y(L,\cdot).
\end{equation}
Here the coefficient $\frac{\del v(p,s)}{\del p}$ is determined by the
nonlinear solution $y=(p,u)$, but the equation is linear in $Y=(P,U)$.

The error between the perturbed and linearized evolutions,
\[
  \big(\ol{\delta P},\ol{\delta U}\big) := \big(\delta p,\delta u\big) - (P,U),
\]
evolves according to the system
\[
  \begin{gathered}
    \ol{\delta P}_x + \ol{\delta U}_t = 0, \qquad
    \ol{\delta U}_x -
    \Big(\frac{\del v(p,s)}{\del p}\,\ol{\delta P}\Big)_t = E_t, \\
    E := v\big(p+\delta p,s\big) - v(p,s)-
    \frac{\del v(p,s)}{\del p}\,P.
  \end{gathered}
\]
Applying the energy estimate \eqref{Hbenergy}, which is (5.1.33) in
\cite{Taylor} and based on Gronwall's inequality, to this system for
the error, gives an estimate sufficient to establish \eqref{DExdef},
which implies that $\mc E^x$ is Frechet differentiable at $y^0$ with
derivative $D\mc E^x(y^0)$ as constructed here.  Because the nonlinear
term $v=v(p+\delta p,s)$ is developed in a Taylor expansion, the error
term is controlled by derivatives of $p$, which implies an essential
loss of derivative in the linearization.  It follows that provided the
evolution $\mc E^L:H^b\to H^b$ is bounded, then for $p^0\in H^b$, the
linearization $D\mc E^L$ is a bounded linear map from $H^{b-1}$ to
$H^{b-1}$~\cite{Ydiff}.

We can similarly define higher Frechet derivatives by an inductive
process.  To differentiate a second time, we evaluate the \emph{linear
  map} $D\mc E^L(y^0)\big[\cdot\big]$, which is nonlinear in $y^0$.
That is, we perturb $y^0$ in an \emph{independent direction}
$Y_2^0$, so we are evaluating
\[
  D\mc E^L(y^0+Y_2^0)[Y_1^0] - D\mc E^L(y^0)[Y_1^0],
\]
and taking the linear part to get $D^2\mc E(y^0)[Y_1^0,Y_2^0]$.  This
is achieved by differentiating \eqref{lineq} in the independent
direction $Y_2^0$.  Denoting the solution of this second
linearization by
\[
  Z = (Q,V), \quad Z(x,\cdot) = D^2\mc E^x(y^0)[Y_1^0,Y_2^0],
\]
the equation for $Z=(Q,V)$ becomes
\begin{equation}
  \label{D2eq}
  Q_x + V_t = 0, \qquad
  V_x - \Big(\frac{\del v(p,s)}{\del p}\,Q\Big)_t =
  \Big(\frac{\del^2 v(p,s)}{\del p^2}\,P_1\,P_2\Big)_t.
\end{equation}
Here $y=(p,u)$ is the result of the nonlinear evolution and
$Y_1=(P_1,U_1)$ and $Y_2=(P_2,U_2)$ are the results of linearized
evolution by \eqref{lineq}, with data $Y_1^0$ and $Y_2^0$,
respectively.  There is no explicit dependence on the velocity
component $U$ of $Y$, because the original system \eqref{nleq} is
already linear in $u$.  Since the data $Y_1^0$ and $Y_2^0$ are
independent, the data for $Z$ vanishes: that is,
$Z^0 = (Q^0,V^0) = (0,0)$.  This vanishing of the data, together with
Gronwall, implies that the error in the expansion is appropriately
bounded with a cubic error estimate.  Note that equation \eqref{D2eq}
for $Z$ is linear inhomogeneous, with the same linear evolution as
$D\mc E(y^0)$, and with inhomogeneous term that can be regarded as
explicit once $y$, $Y_1$ and $Y_2$ have been calculated through the
nonlinear and linear equations \eqref{nleq} and \eqref{lineq},
respectively.

Repeated use of the energy estimate \eqref{Hbenergy} yields the
following theorem; see~\cite{Ydiff} for details.

\begin{theorem}
  \label{thm:D2Evol}
  Theorem \ref{thm:D2E} holds.  Moreover, the nonlinear evolution
  operator $\mc E^L:U\subset H^b\to H^b$ is twice Frechet
  differentiable, with first and second Frechet derivatives around
  $y^0\in H^b$ denoted by
  \[
    \begin{aligned}
      D\mc E^L(y^0)&[\cdot]:H^b\to H^{b-1}, \com{and}\\
      D^2\mc E^L(y^0)&[\cdot,\cdot]:H^b\times H^b\to H^{b-2},
    \end{aligned}
  \]
  respectively.  Both operators are given by evolutions of the form
  \eqref{LXF}, where the linearization $D\mc E^L(y^0)[Y^0]$ is
  homogeneous, while $D^2\mc E^L(y^0)$ is inhomogeneous.
\end{theorem}

Recall that we are working with piecewise $H^b$ entropy profiles.
Since Theorem~\ref{thm:D2Evol} applies only to classical solutions,
and since $x$ is our evolution variable, we treat entropy jumps
explicitly and regard the full evolution as a non-commuting
composition of $H^b$ evolutions, separated by entropy jumps across
which $p$ and $u$ and their time derivatives are unchanged.  Since our
variables $p$ and $u$ are both continuous across jump discontinuities
in the entropy, there is no jump operator between the $H^b$
evolutions.  In other variables, such as Riemann invariants or angle
variables, a jump operator is necessary to connect the $H^b$
evolutions, as in \eqref{thjump}, \eqref{rjump}.  For our further
analysis, we differentiate the composition of evolution operators
twice.

\begin{corollary}
  \label{cor:D2Ecomp}
  Suppose that our evolution operator is a composition of $N$ $H^b$
  evolutions,
  \begin{equation}
    \label{comp}
    \mc E^L = \mc E^{L_N}\,\dots\, \mc E^{L_2}\,\mc E^{L_1}.
  \end{equation}
  Then the linearization $D\mc E^L$ is the composition of individual
  linearizations,
  \[
    D\mc E^{L} = D\mc E^{L_N}\,\dots\, D\mc E^{L_2}\,D\mc E^{L_1},
  \]
  and the second Frechet derivative is
  \[
    D^2\mc E^L = \sum_{k=1}^N D\mc E^{L_N}\,\dots\, D\mc E^{L_{k+1}}
    D^2\mc E^{{L_k}}\circ D\mc E^{L_{k-1}}\,\dots\,D\mc E^{x_1},
  \]
  where $\circ$ denotes
  \[
    D^2\mc E^{{L_k}}\circ \mc A[Y_1^0,Y_2^0] :=
    D^2\mc E^{{L_k}}\big[\mc A\,Y_1^0,\mc A\,Y_2^0\big],
  \]
  and each linear and bilinear map $D^i\mc E^{L_k}$ is evaluated at
  the appropriately evolved nonlinear state,
  \[
    D^i\mc E^{L_k}\big[\cdot,\cdot\big] =
    D^i\mc E^{L_k}(\mc E^{L_{k-1}}\,\dots\,
    \mc E^{L_1}y^0)\big[\cdot,\cdot\big].
  \]
\end{corollary}

\begin{proof}
  The first statement is the chain rule, which for two operators is
  \[
    D\big(\mc E^{L_2}\,\mc E^{L_1}\big)(y^0)[Y_1^0] = 
    D\mc E^{L_2}\big(\mc E^{L_1}(y^0)\big)
    \big[D\mc E^{L_1}(y^0)[Y_1^0]\big].
  \]
  Since the right hand side is a product, differentiating again, in
  direction $Y_2^0$, follows from Leibniz' rule,
  \[
    \begin{aligned}
      D^2\big(\mc E^{L_2}\,\mc E^{L_1}\big)&(y^0)[Y_1^0,Y_2^0]
      = D\mc E^{L_2}\big(\mc E^{L_1}(y^0)\big)
        \big[D^2\mc E^{L_1}(y^0)[Y_1^0,Y_2^0]\big]\\
      &{} + D^2\mc E^{L_2}\big(\mc E^{L_1}(y^0)\big)
        \big[D\mc E^{L_1}(y^0)[Y_1^0],D\mc E^{L_1}(y^0)[Y_2^0]\big].
    \end{aligned}
  \]
  The result now follows by induction.
\end{proof}

Although Theorem~\ref{thm:D2Evol} and Corollary~\ref{cor:D2Ecomp} are
abstract theorems which apply for arbitrary vector inputs $y=(p,u)$
and $Y=(P,U)$, we now revert to the symmetry assumptions \eqref{symm}
made throughout the rest of the paper, which restrict consideration to
data with vanishing velocity component, $u^0=U^0=0$.  This in turn
allows us to regard $p$ as even and $u$ as odd throughout, and to
again treat $y=p+u$, $Y=P+U$ as scalar functions, as in \eqref{ydef}.

Lemma~\ref{lem:himodes} requires the abstract use of
Theorem~\ref{thm:D2Evol} in that it uses twice differentiablity of the
evolution operator $\mc E^L$.  On the other hand,
Lemma~\ref{lem:D2Enz} requires explicit calculation of this second
derivative $D^2\mc E^L$ at the constant state $y^0=0$.  This is
accomplished in the next subsection.

\subsection{Differentiation at the Quiet State}

To complete the proof that nonresonant linear modes perturb to
nonlinear solutions, we must prove Lemma~\ref{lem:D2Enz}, which allows
us to solve the bifurcation equation, and requires calculation of the
second derivative \eqref{dgdz}, namely
\[
  \Big\langle \s\big(k\,t\,\piot-\chi\,k\,{\TS\frac\pi2}\big),
  D^2\mc E(\ol p)
  \big[1,\c(k\,t\,\piot)\big]\Big\rangle.
\]
This in turn requires calculation of the second Frechet derivative in
which one input is the zero-mode,
\begin{equation}
  \label{coeff}
  G(\w) := D^2\mc E(\ol p)\big[1,\c(\w\,\cdot)\big], \com{where}
  \w := k\,\piot,
\end{equation}
of the evolution operator at the constant quiet state.  Here $G(\w)$
is a function of $t$ which depends on the frequency $\w$, and
$\c(\w\,\cdot)$ denotes the input function which is a pure cosine mode
of frequency $\w$, namely $\c(\w\,\cdot):t\mapsto \c(\w\,t)$.

\begin{lemma}
  \label{lem:d2ects}
  The second derivative \eqref{coeff} is given by the solution of the
  linear inhomogeneous SL system
  \begin{equation}
    \label{ODE}
    \begin{aligned}
      \dot{\wh\vp} + \w\,\wh\psi &= 0,\\
      \dot{\wh\psi} - \sg^2\,\w\,\wh\vp &= - v_{pp}\,\w\,\vp_k, 
    \end{aligned}
  \end{equation}
  with vanishing initial data $\wh\vp(0) = \wh\psi(0) = 0$,
  with coefficient $\sg^2 = - v_p(\ol p,s)$, and where
  $v_{pp} = v_{pp}(\ol p,s)\ne 0$.  More precisely, we have
  \begin{equation}
    \label{dEdadz}
    G(\w) = G(\w):t\mapsto
    \wh\vp(\ell)\,\c(\w t)+\wh\psi(\ell)\,\s(\w t),
  \end{equation}
  where $\w = \w_k$ is given by \eqref{Tn}, namely $\w = k\piot$.
\end{lemma}

\begin{proof}
  The quantity $G(\w)$ is directly obtained by evolving system
  \eqref{D2eq} with the appropriate substitutions of the data.
  Because we are linearizing around the quiet state $y=\ol p$, the
  data for the nonlinear equation \eqref{nleq} is this constant state
  $(\ol p,0)$, which is an exact solution of the nonlinear equation,
  and we have $y=\ol p$ throughout the evolution.

  Next, we evolve the data $Y_1^0 = 1$ and $Y_2^0 = \c(\w \cdot)$,
  which are the inputs to $D^2\mc E(\ol p)[\cdot,\cdot]$ in
  \eqref{coeff} by the linearized equation \eqref{lineq}, evaluated on
  the solution we are linearizing around, namely $y=p = \ol p$.  This
  becomes the linear equation
  \begin{equation}
    \label{linpbar}
    P_x + U_t = 0, \qquad
    U_x + \sg^2\,P_t = 0, \qquad
    \sg^2 := -\frac{\del v}{\del p}(\ol p,s),
  \end{equation}
  with initial data $(P_1^0,U_1^0) = (1,0)$ and
  $(P_2^0,U_2^0) = \big(\c(\w\cdot),0\big)$, respectively.  The
  first of these is solved by inspection, to give $(P_1,U_1)=(1,0)$
  throughout the evolution, and is the observation that quiet states
  are preserved by the evolution.  The second is solved by separating
  variables on the corresponding linearized mode,
  \[
    P_2 = \vp_k(x)\,\c(\w\,t), \qquad
    U_2 = \psi_k(x)\,\s(\w\,t).
  \]
  Plugging this in to \eqref{linpbar} gives the SL system
  \[
    \dot\vp_k + \w\,\psi_k = 0, \qquad
    \dot\psi_k - \w\,\sg^2\,\vp_k = 0,
  \]
  with initial data $\big(\vp_k(0),\psi_k(0)\big)=(1,0)$.  Note that
  this is exactly the previously studied SL system \eqref{SL1},
  \eqref{SL3}.

  To calculate the second derivative $G(\w)$ of \eqref{coeff}, we now
  solve system \eqref{D2eq} at the nonlinear quiet state solution
  $\ol p$, and with $P_1$ and $P_2$ as above, and vanishing initial
  data.  System \eqref{D2eq} then becomes
  \begin{equation}
    \label{D2eqpbar}
    \begin{gathered}
      Q_x + V_t = 0, \qquad
      V_x + \sg^2\,Q_t =
      \big(v_{pp}\,\vp_k\,\c(\w t)\Big)_t,\\
      Q(0,t)=V(0,t) = 0,
    \end{gathered}
  \end{equation}
  where
  \[
    \sg^2=-\frac{\del v}{\del p}(\ol p,s), \com{and}
    v_{pp} = \frac{\del ^2v}{\del p^2}(\ol p,s),
  \]
  and these coefficients depend continuously on the $H^b$ entropy
  $s=s(x)$.  We again solve by separating variables, so set
  \[
    Q = \wh\vp_k(x)\,\c(\w\,t), \qquad
    V = \wh\psi_k(x)\,\s(\w\,t).
  \]
  Plugging this ansatz into \eqref{D2eqpbar} gives the inhomogeneous
  SL system \eqref{ODE}, and this completes the proof of the lemma.
\end{proof}

\subsection{Proof of Lemma~\ref{lem:D2Enz}}

We solve \eqref{ODE} using Duhamel's principle with the fundamental
solution $\Psi(x;\w)$ given in \eqref{Psi}: it is straight-forward to
check that the solution of \eqref{ODE} is
\begin{equation}
  \label{duhamel}
  \(\wh\vp(x)\\\wh\psi(x)\) = -\w\,\int_0^x
  \Psi(x-x';\w)\(0\\1\)\vp_k(x')\,v_{pp}\;dx',
\end{equation}
and where $v_{pp}=v_{pp}\big(\ol p,s(x')\big)$.  Recall that genuine
nonlinearity is the condition that $v_{pp}(p,s)>0$, and we have
assumed this holds uniformly.  It turns out that this condition alone
is sufficiently strong to ensure that we can solve the bifurcation
equation \eqref{bif}.

\begin{lemma}
  \label{lem:duh}
  The solution of \eqref{ODE} can be written as
  \begin{equation}
    \label{hatpsi}
    \(\wh\vp(x)\\[2pt]\wh\psi(x)\) = \Psi(x;\w)\,\(a(x)\\b(x)\),
  \end{equation}
  for appropriate functions $a$ and $b$, and we have $b(x) < 0$ for
  all $x > 0$.
\end{lemma}

\begin{proof}
  Because $\Psi$ is a fundamental solution, using the ansatz
  \eqref{hatpsi} in \eqref{ODE} yields the simplified system
  \begin{equation}
    \label{abODE}
    \(a\\b\)^{\dot{}} = - v_{pp}\,\w\,\vp\,\Psi^{-1}\,\(0\\1\),
    \qquad a(0) = b(0) = 0.
  \end{equation}
  Next, since $\det\Psi \equiv 1$, writing
  \[
    \Psi = \(\vp&\wt\vp\\\psi&\wt\psi\),
    \com{we have}
    \Psi^{-1} = \(\wt\psi&-\wt\vp\\-\psi&\vp\),
  \]
  and the second component of \eqref{abODE} simplifies to
  \[
    \dot b = -v_{pp}\,\w\,\vp^2 < 0.
  \]
  It follows that $b(x) < 0$ for all $x > 0$, as required.
\end{proof}

In summary, Lemma~\ref{lem:d2ects} gives the second derivative of the
evolution in terms of the fundamental solution, which is calculated in
Lemma~\ref{lem:SLfund}.  Lemma~\ref{lem:duh} explicitly calculates a
nonvanishing term due to genuine nonlinearity, as seen by the $v_{pp}$
term in \eqref{abODE}.  This then allows us to solve the bifurcation
equation as stated in Lemma~\ref{lem:D2Enz} above.  For convenience we
restate the Lemma.

\begin{lemma}[Lemma ~\ref{lem:D2Enz}]
  For a genuinely nonlinear constitutive equation, the second
  derivative $D^2\mc E(\ol p)$ satisfies
  \begin{equation}
  \label{D2E01c}
  \Big\langle \s\big(k\,t\,\piot-\chi\,k\,{\TS\frac\pi2}\big),
  D^2\mc E(\ol p)
  \big[1,\c(k\,t\,\piot)\big]\Big\rangle\ne 0
  \end{equation}
  so that the derivative
  \[
    \frac{\del g}{\del z}\Big|_{(0,0)} = 
    \frac{\del^2f}{\del z\,\del\a}\Big|_{(0,0)} \ne 0.
  \]
\end{lemma}

\begin{proof}[Proof of Lemma~\ref{lem:D2Enz}]
  Recall that we must show that the quantity \eqref{D2E01c} is
  nonzero, which according to \eqref{dgdz} is
  \[
    \frac{\del g}{\del z}\Big|_{(0,0)} =
    \frac{\del^2f}{\del z\,\del\a}\Big|_{(0,0)} =
    \Big\langle \s\big(k\,t\,\piot-\chi\,k\,{\TS\frac\pi2}\big),
    D^2\mc E(\ol p)
    \big[1,\c(k\,t\,\piot)\big]\Big\rangle.
  \]
  Using \eqref{coeff}, we substitute  \eqref{dEdadz} to get
  \[
    \frac{\del g}{\del z}\Big|_{(0,0)} = \Big\langle \s\big(\w\,t
    -\chi\,k\,{\TS\frac\pi2}\big), G\big(k\,\piot\big)\Big\rangle =
    \begin{cases}
      \wh\vp(\ell), & k\text{ odd and } \chi=1,\\
      \wh\psi(\ell), & \text{otherwise},
    \end{cases}
  \]
  respectively, where $(\wh\vp,\wh\psi)$ solve \eqref{ODE}.  Recall
  that $\chi$ is an indicator of which boundary condition is used:
  $\chi=1$ corresponds to the periodic condition, $\chi=0$ to the
  acoustic condition.
  
  Evaluating \eqref{hatpsi} at $x=\ell$, we get
  \[
    \(\wh\vp(\ell)\\[2pt]\wh\psi(\ell)\) =
    \(\vp(\ell)&\wt\vp(\ell)\\\psi(\ell)&\wt\psi(\ell)\)
    \,\(a(\ell)\\b(\ell)\),
    \com{with} b(\ell)<0.
  \]
  By our choice \eqref{abc} of $\w$, for $k$ odd and $\chi=1$, we
  have \eqref{kodd}, so that $\vp(\ell)=0$, which implies that
  \[
    \wh\vp(\ell) = \wt\vp(\ell)\,b(\ell) \ne 0,
  \]
  and similarly for $k$ even or $\chi=0$, we have \eqref{keven}, which
  is $\psi(\ell)=0$, so we must have
  \[
    \wh\psi(\ell) = \wt\psi(\ell)\,b(\ell) \ne 0,
  \]
  and hence \eqref{D2E01c} holds as claimed and the proof is complete.
\end{proof}

The proof of Lemma~\ref{lem:D2Enz} completes the the proof of
Theorem~\ref{thm:bifurc}, which asserts the global existence of pure
tone time-periodic solutions, which are perturbations of a nonresonant
linearized mode.  In the next section we show that non-resonant modes
are generic.

\section{Piecewise Constant Entropy}
\label{sec:pwc}

In order to give a precise sense in which nonresonant modes are
generic, even though resonant modes are dense, we restrict to the
finite dimensional case of piecewise constant entropy profiles with a
fixed number of distinct jumps.  In this case we can derive explicit
formulas which give the dependence of the SL eigenfrequencies on the
entropy profile.  We note that piecewise constant profiles are
foundational because they are dense in $L^1$ and they form the basis
of many numerical approaches to SL problems~\cite{Pryce}.

In Theorem \ref{thm:gen} below, we establish that nonresonant profiles
form a set of full measure in this space of piecewise constant
profiles with a finite number of jumps, even though in Lemma
\ref{lem:res}, we show that resonant entropy profiles are dense in
$L^1$.  This density of resonant profiles provides some explanation of
the difficulty of such problems and why Nash-Moser methods are a
technically difficult as they are.

For piecewise constant entropy profiles, we exactly solve the
linearized evolution and therefore explicitly calculate the divisors
$\d_j(T)$ and eigenfrequencies $\w_k$ defined in \eqref{Tn} and
\eqref{dkcts}, respectively.  In this case, the evolution is a
non-commuting composition of \emph{isentropic} evolution operators,
and all nonlinear effects of entropy changes are encoded in the
Rankine-Hugoniot jump conditions at the jumps, which are imposed by
continuity of $p$ and $u$.  The corresponding linearization around
each constant quiet state $\ol p$ is again the SL system \eqref{SL3},
but now with constant coefficient $\sg$ on each subinterval.
This changing constant $\sg$ is reflected as a change of sound speed
in the PDE.

\subsection{Parameterization of the Profile}

We introduce parameters for the piecewise constant entropy profile, as
follows.  Suppose there are $N$ distinct entropy levels, and the
jumps are located at points
\[
  0 =: x_0 < x_1 < x_2 < \dots < x_{N-1} < x_N := \ell,
\]
and denote the successive values of the entropy by $s_i$, so that
\begin{equation}
  \label{pwcs}
  s(x) = s_i, \qquad x_{i-1} < x < x_i.
\end{equation}
We denote this class of piecewise constant entropy fields by
$\mc P_N$.

The corresponding linearized wavespeeds are similarly denoted
\begin{equation}
  \label{pwcsg}
  \sg_i := \sg(\ol p,s_i), \qquad x_{i-1} < x < x_i,
\end{equation}
this function being given in \eqref{sigma}.  As in \eqref{Jdef},
\eqref{Jj}, we define the $i$-th jump by
\begin{equation}
  \label{Ji}
  J_i := \frac{\sg_i}{\sg_{i+1}}, \qquad i = 1,\dots, N-1.
\end{equation}
We denote the width of the $i$-th subinterval on which the entropy is
constant, by
\begin{equation}
  \label{Lidef}
    L_i := x_i - x_{i-1}, \com{so that}
    \ell = \sum_{i=1}^N L_i.
\end{equation}
In \cite{TYperStr, TYperLin}, we found that the linearized isentropic
evolution can be interpreted as a rotation in Riemann invariants, and
in this paper we identify this rotation with the angle variable
introduced in \eqref{pruf}, \eqref{theta}.  In this context we define
the \emph{evolution angle}
\begin{equation}
  \label{thdef}
\t_i := \sg_i\,L_i.
\end{equation}

The SL evolution across the entire interval $[0,\ell]$ is the
composition of the distinct \emph{isentropic} SL evolutions,
\begin{equation}
  \label{pwccomp}
  \Psi(\ell;\w) = \Psi_N\,\circ\dots\circ\Psi_1, \qquad
  \Psi_i:=\Psi(x_i-x_{i-1};\w),
\end{equation}
where each $\Psi_i$ is the fundamental matrix of the evolution over the
entire subinterval on which $\sg=\sg_i$ is constant.  For each of these,
the SL system \eqref{SL1} reduces to a constant coefficient system,
and the corresponding angle representation \eqref{theta},
\eqref{rint}, namely
\[
  \vp(x) := r(x)\,\frac1{\sqrt{\sg_i}}\,\c\big(\t(x)\big), \qquad
  \psi(x) := r(x)\,\sqrt{\sg_i}\,\s\big(\t(x)\big),
\]
reduces to the simpler system
\begin{equation}
  \label{cpwc}
  \dot\t = \w\,\sg_i, \qquad \dot r = 0, \qquad x_{i-1} < x < x_i.
\end{equation}
Therefore $\Psi_i$ is no more than the evolution in $(\vp,\psi)$
determined by \eqref{cpwc}.  In the following lemma, we show that
linearized isentropic evolution can be realized as rotation in Riemann
invariants, which is \eqref{cpwc}, and transformation between physical
variables $(\vp,\psi)$ and Riemann invariants $(r\,\c\t,r\,\s\t)$ is
conjugation by the jump matrix \eqref{MJmat}.  This is a restatement
of the linearized structures developed by the authors
in~\cite{TYperStr, TYperLin,TYperBif}.

\begin{lemma}
  \label{lem:pws}
  When the entropy is piecewise constant, the fundamental matrix of
  the SL system \eqref{SL3}, satisfying
  \[
    \(\vp(\ell)\\\psi(\ell)\) = \Psi(\ell;\w)
    \(\vp(0)\\\psi(0)\),
  \]
  can be written as the $2\times2$ matrix product
  \begin{equation}
    \label{pwcprod}
    \begin{aligned}
      \Psi(\ell;\w) = M(\sg_N)\,&R(\w\,\t_N)\,M(\sg_N)^{-1}\,M(\sg_{N-1})
                                  \dots\\
    &\dots\,R(\w\,\t_2)\,M(\sg_2)^{-1}\,M(\sg_{1})\,R(\w\,\t_1)
    \,M^{-1}(\sg_1) \\ = M(\sg_N)\,&R(\w\,\t_N)\,M(J_{N-1})
                                  \dots\\
    &\dots\,R(\w\,\t_2)\,M(J_{1})\,R(\w\,\t_1)
    \,M^{-1}(\sg_1).
    \end{aligned}
  \end{equation}
  Here $R(\cdot)$ is the usual rotation matrix and $M(\cdot)$ is the
  jump matrix defined in \eqref{MJmat}, and $J_i$ is given by
  \eqref{Ji}.  In particular the linearized evolution of data with
  given angle $\t^0$ can be evaluated inductively as a series of
  rotations and jumps in Riemann invariants, as
  \[
    \Psi(\ell;\w)\,r^0\(\c(\t^0)\\\s(\t^0)\) =
    r_N\,\(\c(\g_N)\\\s(\g_N)\),
  \]
  where the intermediate angles $\g_m$ are defined by
  \begin{equation}
    \label{gamdef}
    \begin{aligned}
      \g_0&:=h(1/\sg_1,\w\t^0)\\
      \g_{m} &:= h\big(J_{m},\g_{m-1} + \w\t_{m}),
                   \qquad m = 1,\dots,N-1,\\
      \g_N &:= h(\sg_{N},\g_{N-1} + \w\t_{N}),
    \end{aligned}
  \end{equation}
  and the intermediate radii are 
  \begin{equation}
    \label{irdef}
    \begin{aligned}
      r_0&:= \rho(1/\sg_1,r^0)\\
      r_m &:= r_{m-1}\,\rho\big(J_{m},\g_{m-1} + \w\t_{m}),
      \qquad m = 1,\dots,N-1,\\
      r_N &:= r_{N-1}\,\rho(\sg_{N},\g_{N-1} + \w\t_{N}),
    \end{aligned}
  \end{equation}
  and where the functions $h$ and $\rho$ are defined in \eqref{hdef}
  and \eqref{rJ}, respectively.  Moreover, the change in radius
  $r_N/r^0$ is uniformly bounded above and below by a constant
  depending only on the entropy profile.
\end{lemma}

\begin{proof}
  We constructed the general fundamental solution in
  Lemma~\ref{lem:SLfund}.  Using the simplified evolution \eqref{cpwc}
  in \eqref{SLfund}, we get, for each $i$,
  \[
    \Psi(L_i;\w) =
    M(\sg_i)\,R(\w\,\sg_i\,L_i)\,M^{-1}(\sg_i),
  \]
  where $M(\cdot)$ is given by \eqref{MJmat}, with $\t_i = \sg_i\,L_i$
  by \eqref{Lidef}.  Since the matrix $M$ is multiplicative,
  $M(q_1\,q_2) = M(q_1)\,M(q_2)$, assertion \eqref{pwcprod} follows
  from \eqref{Ji}.

  The inductive formulas \eqref{gamdef} and \eqref{irdef} follow
  directly from the observation that the angle changes linearly in
  linearized isentropic evolution, together with repeated application
  of Lemma~\ref{lem:jump} at each jump $J_m$.
\end{proof}

Recall that the eigenfrequencies $\w_k$ are defined by the implicit
condition \eqref{abc}, namely
\[
  \t(\ell,\w_k) = k\,\frac\pi2, \com{or}
  \t(\ell,\w_k) = k\,\pi,
\]
with initial angle $\t^0=0$, for the indicator $\chi=1$ or $\chi=0$,
corresponding to periodic and acoustic boundary conditions,
respectively.  As above, for convenience, we treat the more general
periodic condition $\chi=1$; restricting to even $k$ then recovers the
other case.

\begin{corollary}
  \label{cor:efreq}
  For a piecewise constant entropy profile, the $k$-th eigenfrequency
  is the root $\w = \w_k$ of the implicit equation
  \begin{equation}
    \label{weq}
    \w\,\t_N + \g_{N-1} = k\,\frac\pi2,
  \end{equation}
  where $\g_0=0$ and each $\g_m$ is given inductively by
  \eqref{gamdef}, namely
  \begin{equation}
    \label{wmdef}
    \g_{m} := h\big(J_{m},\g_{m-1} + \w\,\t_{m}),
    \qquad m = 1,\dots,N-1.    
  \end{equation}
\end{corollary}

\begin{proof}
  First note that for any $J>0$, and any integer $n$,
  \[
    h\Big(J,n\,\frac\pi2\Big) = n\,\frac\pi2
  \]
  so that, because $\t^0=0$, the initial and final entropy adjustments
  $M(1/\sg_1)$ and $M(\sg_N)$ do not affect the angle, and can be
  ignored.  Similarly changes to the radii $r_m$ play no part.  Thus
  the closed system \eqref{wmdef}, \eqref{weq} determines the
  eigenfrequency $\w_k$ and the proof is complete.
\end{proof}

\subsection{Time Period and Frequencies}

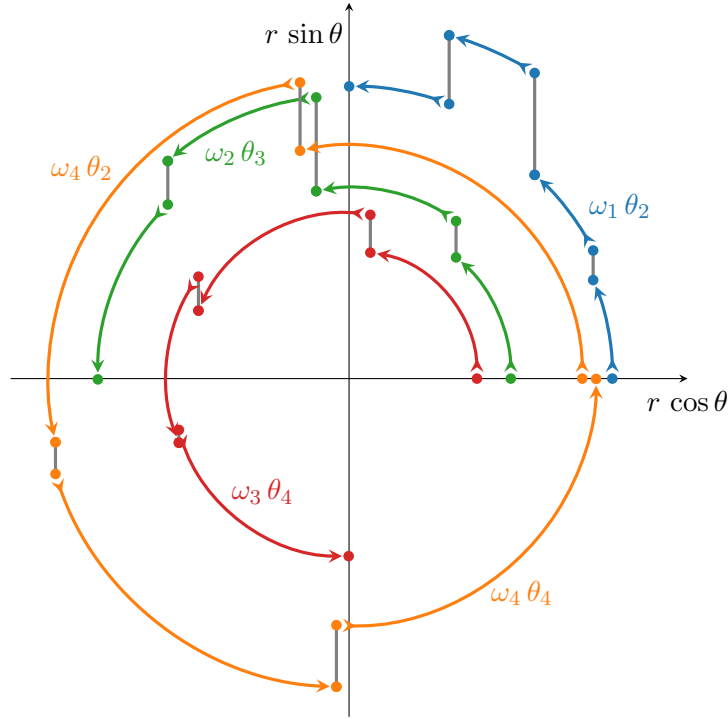
\begin{figure}[thb]
  \begin{tikzpicture}[>=stealth]
  \draw[->] (-4.5,0) -- (4.5,0);
  \draw[->] (0,-4.5) -- (0,5);
  \node at (4.5,-0.3) {$r\,\cos\t$};
  \node at (-0.6,4.6) {$r\,\sin\t$};

  \def\th0{22}

\foreach \rad/\mult/\C in
  {3.5/1/C0,2.15/2.208/C2,3.1/4.6425/C1,1.7/3.655/C3}{%
\fill[color=\C] (\rad, 0) circle (2 pt);
\draw[very thick,color=\C,>->,shorten <=2.4pt,shorten >=2.5pt]
(\rad0, 0) arc (0:\mult*\th0:\rad) coordinate (ee);
  \foreach \J/\th in {1.3/20,1.5/15,0.8/20}
  {%
    \draw[very thick,color=C7] 
    let \p2 = (ee) in (ee) -- (\x2,\J*\y2) coordinate (ff);
    \fill[color=\C] (ee) circle (2 pt);
    \fill[color=\C] (ff) circle (2 pt);
    \draw[very thick,color=\C,>->,shorten <=2.4pt,shorten >=2.5pt]
    let \p1 = (ff), \n1 = {veclen(\x1,\y1)},
    \n2 = {atan(\y1/\x1)+180*(\x1<0)}
    in (ff) arc [radius=\n1, start angle=\n2, delta angle=\mult*\th]
    coordinate (ee);
    \fill[color=\C] (ee) circle (2 pt);
  }
}

\node[C0] at (3.6,2.3) {$\w_1\,\theta_2$};
\node[C2] at (-1.5,3) {$\w_2\,\theta_3$};
\node[C1] at (-3.55,2.8) {$\w_4\,\theta_2$};
\node[C3] at (-1.15,-1.5) {$\w_3\,\theta_4$};
\node[C1] at (2.3,-2.85) {$\w_4\,\theta_4$};
\end{tikzpicture}
\caption{Rotations and jumps generating frequencies}
\label{fig:freqs}
\end{figure}

To illustrate the calculation of the eigenfrequencies, in
Figure~\ref{fig:freqs}, we show the corresponding rotations and jumps
(scaled for visibility) for the first four nonzero eigenfrequencies,
in the plane of Riemann invariants $(r\,\c\t,r\,\s\t)$.  Here the
entropy field consists of four constant states separated by three
jumps.  The circular arcs represent linear evolution of the $k$-mode
through the entropy level $\t_m$ so are rotations by $\w_k\,\t_m$,
and the vertical segments represent the action of the jumps, as
encoded by the function $h(J,\t)$, in which the second component is
scaled but the first remains constant.  The constraints on each
eigenfrequency are that the $k$-th curve should start at the $x$-axis
$\t=0$ and should end in the multiple $k\frac\pi2$, which is the $k$-th
axis, counting anti-clockwise.  The four curves are color coded, and
some arcs are labeled: for example, $\w_1\,\t_2$ (blue) is the
evolution of the 1-mode between jumps $J_1$ and $J_2$.  Recall that
according to Lemma \ref{lem:wk}, the eigenfrequencies $\w_k$ are
unique and grow linearly with $k$, as can be seen in this
diagram, because the frequency, or rate of rotation, must increase in
order to reach the $k$-th axis over the same interval $[0,\ell]$.

Recall that the time period $T=T_k$ defined in \eqref{Tn} is
determined by the eigenfrequency of the linearized mode that we are
perturbing.  Thus, for a given entropy profile, we \emph{first} pick
the mode $k$, which then determines the time period by \eqref{Tn},
namely
\[
  T = T_k = k\,\frac{2\pi}{\w_k},
\]
where $\w_k$ is the $k$-th eigenfrequency of that entropy profile.
Having fixed the time period, we can \emph{then} calculate the small
divisors corresponding to the $j$-th modes, $\d_j(T)$, according to
\eqref{dkcts}, namely
\[
  \d_{j}(T) := \(0&1\)\,R(j\chi\pb)\,
  \Psi\big(\ell;j\piot\big)\(1\\0\).
\]
Here the vector $\(1\\0\)$ represents the coefficient of the even
$j$-mode $\c(j\piot)$ at $x=0$, $\Psi(\ell;j\piot)$ evolves this mode
from $x=0$ to $x=\ell$, $R(j\chi\pb)$ shifts the result by a quarter
period for the periodic boundary condition, and $\(0&1\)$ projects
onto the odd mode which vanishes when the boundary condition is met.
The divisor $\d_j(T)$ then measures the degree to which this $j$ mode
fails to satisfy the boundary condition.  Note that these small
divisors are now just numbers, and in general are not directly related
to the other eigenfrequencies of the given entropy profile.  The
$j$-th mode is \emph{resonant}, if the corresponding $j$-th divisor
vanishes.  According to Lemma~\ref{lem:nonres}, this occurs if and
only if $\w_j$ is rationally related to $\w_k$, that is
\[
  \d_j(T) = 0 \com{iff} k\,\w_l = j\,\w_k,
\]
which is the statement that $q\,\w_k$, $q:=j/k$, is another
eigenfrequency of the same entropy profile.

\subsection{Resonance Conditions}

In order to better understand the resonance structure, we give an
explicit geometric procedure for determining the $k$-th eigenfrequency
$\w_k$ for a given piecewise constant entropy profile.  This geometric
construction is based on the observation that all $\t_m$'s in
\eqref{wmdef}, \eqref{weq} are multiplied by $\w$, so can be replaced
by another variable $Z_m$.  It is convenient to introduce the vector
notation
\[
  \begin{aligned}
    \Theta &:= (\t_1,\dots,\t_N)\in\B R_+^N,\\
    J &:= (J_1,\dots,J_{N-1})\in\B R_+^{N-1},\\
    Z &:= (Z_1,\dots,Z_N)\in\B R_+^N,
  \end{aligned}
\]
and we also write $(J,\T)\in \mc P_N=\B R_+^{2N-1}$.

Referring to Corollary \ref{cor:efreq}, we define a function
$\Gamma^J = \g_{N-1}$ on $\B R_+^{N-1}$, as follows.  First set
$\g_0=0$, and inductively define
\emph{intermediate angle functions}
\begin{equation}
  \label{gmdef}
  \begin{gathered}
    \g_{m}:\B R_+^m\to \B R_+ \com{by}\\
    \g_m(Z_1,\dots,Z_m)
    := h\big(J_{m}, Z_{m} + \g_{m-1}(Z_1,\dots,Z_{m-1})\big),
  \end{gathered}
\end{equation}
for $m = 1,\dots,N-1$, and set
\begin{equation}
  \label{Gamma}
  \Gamma^J(Z_1,\dots,Z_{N-1}) := \g_{N-1}(Z_1,\dots,Z_{N-1}).
\end{equation}
Note that each $\g_m$ also depends on the jumps $(J_1,\dots,J_m)$,
which we treat as parameters.  We calculate the derivatives of the
$\g_m$ recursively.

\begin{lemma}
  \label{lem:gam}
  Each function $\g_m$ depends only on the values of $J_l$ and $Z_l$
  for $l\le m$, and the derivatives satisfy
  \begin{equation}
    \label{gmder}
    \frac{\del\g_m}{\del J_m} =
    \frac{\del h}{\del J}\Big|_{(J_m,Z_m+ \g_{m-1})}, \qquad
    \frac{\del\g_m}{\del Z_m} =
    \frac{\del h}{\del Z}\Big|_{(J_m,Z_m+ \g_{m-1})},
  \end{equation}
  derivatives of $h$ being given in \eqref{hder}, together with the
  recursive relations
  \begin{equation}
    \label{gamderj}
    \begin{aligned}
      \frac{\del\g_{m}}{\del J_l}
      &= \frac{\del h}{\del Z}\Big|_{(J_{m},Z_m+\g_{m-1})}\,
        \frac{\del\g_{m-1}}{\del J_l},\\
      \frac{\del\g_{m}}{\del Z_l}
      &= \frac{\del h}{\del Z}\Big|_{(J_{m},Z_m+\g_{m-1})}\,
        \frac{\del\g_{m-1}}{\del Z_l},
    \end{aligned}
  \end{equation}
  for $l<m$.
\end{lemma}

\begin{proof}
  This follows immediately from the definition \eqref{gmdef} together
  with repeated use of the chain rule.
\end{proof}

Using these functions $\g_m$, we now rewrite the equations
\eqref{wmdef}, \eqref{weq}, as the system
\begin{equation}
  \label{Zteq}
  Z_N = k\,\frac\pi2 - \Gamma^J\big(Z_1,\dots,Z_{N-1}\big), \qquad
  Z = \w\,\Theta,
\end{equation}
where the second is a vector equation, $Z_m=\w\,\t_m$, $m=1,\dots,N$.
We intepret this system as follows: $Z_N$ is the graph of the function
$k\frac\pi2-\Gamma^J$, and $Z=\w\,\Theta$ describes the ray in
$\B R_+^N$ determined by the vector $\Theta$, and parameterized by $\w$.
Both equations are satisfied at the point $\wh Z$ at which the graph
and ray intersect, and the eigenfrequency $\w=\w_k$ is then the scaling
parameter such that
\[
  \wh Z = \w_k\,\Theta, \com{that is}
  \w_k = \wh Z_m/\t_m,
\]
for any $m$.  Note that both the ray and function $\Gamma^J$ are
independent of $k$, and as $k$ varies, the only difference is in the
height of the graph.  Also note that this interpretation does
\emph{not} require that $k$ be an integer.

For piecewise constant entropy profiles, we write the $k$-th
eigenfrequency as the function
\begin{equation}
  \label{wk}
  \w_k = \w_k(J,\Theta),\qquad \w_k:\mc P_N=\B R_+^{2N-1}\to \B R_+.
\end{equation}
We showed in Lemma~\ref{lem:wk} that the $k$-th eigenfrequency is
unique, so that $\w_k$ is well-defined, and in the case of piecewise
constant entropies, this is in fact a $C^\infty$ function.

\begin{lemma}
  \label{lem:wkf}
  For piecewise constant entropy profiles with $N-1$ jumps, the
  eigenfrequency $\w_k(J,\Theta)$ defined in \eqref{wk} is a
  $C^\infty$ function of the parameters $(J,\T)$.
\end{lemma}

\begin{proof}
  Recall that $\w=\w_k$ is given by \eqref{Zteq}, which
  we rewrite as
  \[
    \w_k\,\t_N + \G^J\big(\w_k\,\t_1,\dots,\w_k\,\t_{N-1}\big) =
    k\,\frac\pi2.
  \]
  We now change perspective, allowing $k$ to be a continuous variable,
  and instead regard $\w$ as the parameter.  That is, we define the
  function
  \[
    \k:\B R_+\times\mc P_N\to\B R_+, \qquad
    \k := \k(\w;J,\T),
  \]
  by
  \begin{equation}
    \label{kappa}
    \k(\w;J,\T) := \frac2\pi\,\big[\w\,\t_N +
    \G^J(\w\,\t_1,\dots,\w\,\t_{N-1})\big],
  \end{equation}
  so that eigenfrequencies occur at integral values of $\k$,
  \[
    \w_k(J,\T) = \w \com{if and only if}
    \k(\w;J,\T) = k \in \B N.
  \]
  In other words, for each fixed $(J,\T)\in\mc P_N$, $\w=\w_k$ and
  $k = \k(\w)$ are inverse functions.

  It follows that $\k$ is a $C^\infty$ function, and we calculate
  \begin{equation}
    \label{dkdw}
    \frac\pi2\,\frac{\del \k}{\del\w} = \t_N + \sum_{i=1}^{N-1}
    \frac{\del\Gamma^J}{\del Z_i}\,\t_i > 0,
  \end{equation}
  since each term, given recursively by \eqref{gamderj},
  \eqref{gmder}, is positive.  Moreover, by Lemmas \ref{lem:h} and
  \ref{lem:gam}, it follows that this derivative is uniformly bounded,
  that is
  \begin{equation}
    \label{kder}
    m(J,\T) \le \frac{\del \k}{\del\w} \le M(J,\T),
  \end{equation}
  these bounds depending only on the entropy profile and independent of
  $k$, $\w$.  The result now follows by the inverse function theorem.
\end{proof}

\subsection{Density of Resonant Modes}

According to Lemma~\ref{lem:nonres}, the $j$-mode resonates with the
$k$-mode if and only if there is some $l\ne k$ such that
\begin{equation}
  \label{wres}
  k\,\w_l(J,\Theta) = j\,\w_k(J,\Theta).
\end{equation}
For each fixed triple $(k,j,l)$, we regard this as a constraint on the
entropy profile, and taking the union as the triple varies generates
\emph{all possible} resonant entropy profiles within the class of
piecewise constant profiles with $N-1$ entropy jumps.  We proceed by
showing that each such restriction generates a codimension one
$C^\infty$ manifold in $\mc P_N=\B R_+^{2N-1}$, and so has Lebesgue
measure zero in $\B R_+^{2N-1}$.  Since there are countably many such
restrictions, the set of all resonant entropy profiles is also of
Lebesgue measure zero.  This means that generically, piecewise
constant entropy fields with $N-1$ jumps, are \emph{fully
  nonresonant}, in that no two modes resonate with each other.  This
in turn implies that, within the space of piecewise constant entropy
profiles, generically, \emph{all} linear modes perturb to pure tone
periodic solutions of the nonlinear compressible Euler equations.

Recall that piecewise constant functions are dense in $L^1$.  Because
resonance requires only the existence of a single pair of
eigenfrequencies satisfying \eqref{wres}, it follows that resonant
profiles are dense in $L^1$.

\begin{lemma}
  \label{lem:res}
  The set $\mc Z$ of all piecewise constant resonant entropy profiles
  is dense in $L^1$.
\end{lemma}

\begin{proof}
  We first reduce to a finite dimensional entropy profile.  Given an
  arbitrary $\wh s\in L^1$ and $\epsilon>0$, choose a piecewise
  constant entropy profile with $N-1\ge 1$ jumps, denoted
  $s(J_0,\T_0)\in\mc P_N$, so that
  \[
    \|\wh s - s(J_0,\T_0)\|_{L^1} < \frac\epsilon2.
  \]
  Since the norm is continuous, we can find some $\d>0$ such that
  whenever
  \[
    (J,\T)\in \ol{B_\d} := \Big\{( J,\T)\subset\mc P_{N}\;\Big|\;
    \big|( J,\T)-(J_0,\T_0)\big|_{\B R^{2N-1}}\le\d\Big\},
  \]
  then also
  \[
    \|s(J,\T) - s(J_0,\T_0)\|_{L^1} < \frac\epsilon2,
  \]
  so it suffices to show that the resonant set $\mc Z_N$ is relatively
  dense in $\ol{B_\d}$.

  Recall that the profile $(J,\T)\in\mc Z_N$ is resonant if
  \eqref{wres} holds for \emph{any} triple $(k,j,l)$, namely
  \begin{equation}
    \label{reseq}
    k\,\w_l = j\,\w_k, \com{or}
    \frac jl\, \frac{\w_k(J,\T)}k = \frac{\w_l(J,\T)}l.
  \end{equation}
  From Lemma~\ref{lem:wk}, we have that $\w_l$ grows linearly with $l$:
  that is, by \eqref{wkok}, if we set
  \[
    \Lambda(J,\T) :=
    \frac{\pi/2}{\int_0^\ell \sg(y)\;dy}, \com{then}
    \lim_{l\to\infty} \frac{\w_l(J,\Theta)}l 
    = \Lambda(J,\T).
  \]
  For $k$ fixed, we exhibit a curve 
  \[
    (J_\nu,\T_\nu) := \big(J(\nu),\T(\nu)\big) \in \ol{B_\d}, \qquad
    \nu\in(-\ol\nu,\ol\nu),
  \]
  with $J(0)=J_0$ and $\T(0)=\T_0$, along which $\Lambda(J,\T)$ varies
  but $\w_k$ is constant, as follows.  Set $\ol\w_k = \w_k(J_0,\T_0)$,
  and let jumps $J_{i,\nu}$ and widths $\t_{i,\nu}$ vary smoothly with
  $\nu$, for $i=1,\dots,N-1$.  Then interpret \eqref{kappa}, which
  defines $\w_k$, as a constraint on $\t_{N,\nu}$, that is, take
  \[
    \t_{N,\nu} := \frac k{\ol\w_k}\,\frac\pi2 -
    \frac1{\ol\w_k}\,\G^{J_\nu}
    \big(\ol\w_k\,\t_{1,\nu},\dots,\ol\w_k\,\t_{N-1,\nu}\big).
  \]
  
  Now pick $\nu_1$, $\nu_2$ such that
  \[
    \lambda :=
    \Lambda(J_{\nu_2},\T_{\nu_2})-\Lambda(J_{\nu_1},\T_{\nu_1}) > 0,
  \]
  say.  Since $\w_l/l\to \Lambda$, there is some $L$ large enough so that
  \[
    \Big|\frac{\w_l(J_{\nu_i},\T_{\nu_i})}l -
    \Lambda(J_{\nu_i},\T_{\nu_i})\Big|< \frac{\lambda}3,
    \qquad i = 1,2,
  \]
  whenever $l>L$, and in particular,
  \[
    \frac{\w_l(J_{\nu_2},\T_{\nu_2})}l -
    \frac{\w_l(J_{\nu_1},\T_{\nu_1})}l > \frac\lambda3,
  \]
  whenever $l>L$.  Now referring to \eqref{reseq}, by density of the
  rationals, we can find integers $j$ and $l$, with $l>L$, such that
  \[
    \frac{\w_l(J_{\nu_1},\T_{\nu_1})}l <
    \frac jl\,\frac{\ol\w_k}k <
    \frac{\w_l(J_{\nu_2},\T_{\nu_2})}l.
  \]
  Finally, by continuity, it follows that there is some
  $\nu\in(\nu_1,\nu_2)$ such that
  \[
    \frac{\w_l(J_{\nu},\T_{\nu})}l = \frac jl\,\frac{\ol\w_k}k,
  \]
  which is exactly the resonance condition \eqref{reseq}.  Since this
  profile $(J_{\nu},\T_{\nu})\in\ol{B_\d}$, the proof is complete.
\end{proof}

\subsection{Resonant Sets as a Union of Smooth Manifolds}

Piecewise constant entropy profiles consist of $N$ entropy levels
$s_i$ separated by $N-1$ jumps, which for fixed $\ol p$ and $s_1$, are
parameterized by the vector $(J,\Theta)\in\mc P_N$.  Referring to
\eqref{wres}, for each triple $(k,j,l)\in\B N^3$, define the set
\begin{equation}
  \label{Zkjl}
  \mc Z^N_{k,j,l} := \Big\{(J,\Theta)\in\mc P_N\;\Big|\;
  k\,\w_l(J,\Theta) = j\,\w_k(J,\Theta)\Big\},
\end{equation}
where $\w_k$ is the function \eqref{wk}.  The \emph{resonant sets} are
defined to be the union of all of these,
\begin{equation}
  \label{Zres}
  \mc Z_N = \bigcup_{k,j,l}\mc Z^N_{k,j,l}, \com{and}
  \mc Z = \bigcup_N\mc Z_N.
\end{equation}
Here each $\mc Z_N\subset\mc P_N$ is the set of resonant profiles with
$N-1$ jumps, and $\mc Z$ is the set of all piecewise constant resonant
profiles.

\begin{theorem}
  \label{thm:Z}
  For fixed $N$, and for each triple $(k,j,l)\in\B N^3$, the set
  \[
    \mc Z^N_{k,j,l}\subset\mc P_N = \B R_+^{2N-1},
  \]
  of piecewise constant profiles satisfying the resonance
  condition \eqref{wres}, forms a $C^\infty$ submanifold of dimension
  $2N-2$, that is of codimension one, in the parameter space $\mc P_N$
  of entropy profiles with $N-1$ jumps.
\end{theorem}

\begin{proof}
  The set $\mc Z^N_{k,j,l}$ is defined by \eqref{Zkjl}, that is by the
  single equation
  \[
    k\,\w_l(J,\Theta) = j\,\w_k(J,\Theta).
  \]
  We rewrite this condition in terms of the function $\k$ defined in
  \eqref{kappa}, as follows.  We set
  \[
    \w_k =: \w,\quad \w_l =: \ol\w, \com{and} j/k =: q,
  \]
  so the condition becomes
  \[
    \k(\w;J,\T) = k, \quad
    \k(\ol\w;J,\T) = l, \com{and}
    \ol\w = q\,\w.
  \]
  Eliminating $\ol w$, we write this as the system
  \[
    q\,\k(\w;J,\T) = q\,k = j, \qquad
    \k(q\,\w;J,\T) = l,
  \]
  and the profile $(J,\T)\in\mc Z^N_{k,j,l}$ if this system is
  satisfied for the given $k$ and $l$ and for $q = j/k$.

  We thus regard $k$, $l$ and $q$ as fixed and having been given, and
  considering \eqref{kappa} define the map
  \begin{equation}
    \label{Fdef}
    G_q:\B R_+\times\mc P_N\to\B R_+^3, \qquad
    G_q(\w;J,\T) :=
    \(\k(\w;J,\T)\\q\,\k(\w;J,\T)\\\k(q\,\w;J,\T)\).
  \end{equation}
  Then the particular resonant set we are looking at is the inverse
  image
  \[
    \mc Z^N_{k,j,l} = G_q^{-1}\(k\\j\\l\),
  \]
  and the result will follow if we prove that the derivative $DG_q$
  has rank two as a linear map $\B R^{2N}\to\B R^3$, the first two
  rows being dependent.  Assuming this is so, the implicit function
  theorem implies $G_q^{-1}(j,l)$ is a submanifold of the stated
  dimension.

  It remains to show that $DG_q$ is rank two.  Dropping the first row
  and writing $\wh G$, we have
  \begin{equation}
    \label{gradG}
    D\wh G_q = \(     q\,\frac{\del\k}{\del\w}\big|_\w
    & q\,\nabla_{(J,\T)}\k\big|_\w\\[2pt]
    q\,\frac{\del\k}{\del\w}\big|_{q\,\w}
    & \nabla_{(J,\T)}\k\big|_{q\,\w}\),
  \end{equation}
  and we must show that this is full rank; this will certainly
  follow if the gradients $q\,\nabla_{(J,\T)}\k\big|_\w$ and
  $\nabla_{(J,\T)}\k\big|_{q\,\w}$ are independent.  Referring to the
  definition \eqref{kappa}, it is immediate that
  \[
    q\frac{\del\k}{\del\t_N}\Big|_\w = \frac2\pi\,q\,\w
    = \frac{\del\k}{\del\t_N}\Big|_{q\,\w},
  \]
  and so the result will follow if we can prove that the other
  components of the gradient aren't equal, that is provided we can
  show
  \begin{equation}
    \label{gradF}
    \nabla_{(J,\T)}\big(q\,\k(\w;J,\T) - \k(q\,\w;J,\T)\big) \ne 0.
  \end{equation}

  To this end, and again referring to \eqref{kappa}, we define the
  function
  \begin{equation}
    \label{FJdef}
    F:\B R^{2N-2}\to\B R, \qquad
    F(J,Y) := \Gamma^J(q\,Y) - q\,\Gamma^J(Y),
  \end{equation}
  where
  \[
    Y\in\B R^{N-1} \com{is}
    Y := \big(\w\,\t_1,\dots,\w\,\t_{N-1}\big).
  \]
  We use the (backwards) recursive formulas for the derivatives from
  Lemma \ref{lem:gam} to calculate the derivatives of $F$.  First,
  using \eqref{FJdef}, \eqref{Gamma}, and \eqref{gmder}, we get
  \[
    \frac{\del F}{\del Y_{N-1}}
    = \frac{\del h}{\del Z}
    \bigg|_{(J_{N-1},\ol{qY}_{N-1})}\,q
    - q\,\frac{\del h}{\del Z}
    \bigg|_{(J_{N-1},\ol{Y}_{N-1})},
  \]
  where we use the notation
  \[
    \begin{aligned}
      \ol Y_m &:= Y_m+\g_{m-1}(Y_1,\dots,Y_{m-1}),\\
      \ol{qY}_m &:= qY_m+\g_{m-1}(qY_1,\dots,qY_{m-1}),      
    \end{aligned}
  \]
  in evaluating functions at intermediate stages.  Next, using
  \eqref{hder}, we write this as
  \[
    \begin{aligned}
      \frac{\del F}{\del Y_{N-1}}
      &= \frac{q\,J\,(1+\tf^2(\ol{qY}_{N-1}))}
        {1 + J^2\,\tf^2(\ol{qY}_{N-1})}
        - \frac{q\,J\,(1+\tf^2(\ol Y_{N-1}))}
        {1 + J^2\,\tf^2(\ol Y_{N-1})}\\[3pt]
      &= \frac{q\,J\,(J^2-1)\,\big(\tf^2(\ol Y_{N-1})-
        \tf^2(\ol{qY}_{N-1})\big)}
        {\big(1 + J^2\,\tf^2(\ol{qY}_{N-1})\big)\,
        \big(1 + J^2\,\tf^2(\ol Y_{N-1})\big)},
    \end{aligned}
  \]
  with $J$ abbreviating $J_{N-1}$, and in which $\tf:=\tan$.  It
  follows that
  \[
    \frac{\del F}{\del Y_{N-1}} = 0
    \com{iff}
    \tf(\ol{qY}_{N-1}) = \pm\tf(\ol Y_{N-1}),
  \]
  unless $J=J_{N-1}=1$.  Similarly, from \eqref{hder} we get
  \[
    \begin{aligned}
      \frac{\del F}{\del J_{N-1}}
      &= \frac{\del h}{\del J}
      \bigg|_{(J_{N-1},\ol{qY}_{N-1})}
      - q\,\frac{\del h}{\del J}
      \bigg|_{(J_{N-1},\ol{Y}_{N-1})}\\[2pt]
      &= \frac{\tf(\ol{qY}_{N-1})}{1 + J^2\,\tf^2(\ol{qY}_{N-1})}
        - q\,\frac{\tf(\ol Y_{N-1})}{1 + J^2\,\tf^2(\ol Y_{N-1})}.
    \end{aligned}
  \]
  Thus, if we assume $J_{N-1}\ne 1$, and assume both
  \[
    \frac{\del F}{\del Y_{N-1}}\bigg|_{(J,\ol Y_{N-1})} = 0 \com{and}
    \frac{\del F}{\del Y_{N-1}}\bigg|_{(J,\ol Y_{N-1})} = 0,
  \]
  then we must have both
  \[
    \tf(\ol{qY}_{N-1}) = \pm\tf(\ol Y_{N-1}) \com{and}
    \tf(\ol{qY}_{N-1}) - q\,\tf(\ol Y_{N-1}) = 0,
  \]
  which in turn implies
  \[
    \tf(\ol{qY}_{N-1}) = \tf(\ol Y_{N-1}) = 0.
  \]
  It then follows from the definitions \eqref{FJdef}, \eqref{Gamma}
  and \eqref{hdef}, that
  \[
    F(J,\ol Y_{N-1}) = \arctan\big(J_N\,\tf(\ol{qY}_{N-1})\big) -
    q\,\arctan\big(J_N\,\tf(\ol Y_{N-1})\big) = 0.
  \]

  On the other hand, if $J_{N-1}=1$, which is the degenerate case in
  which there is no entropy jump, we carry out the same calculation
  for $\frac{\del F}{\del Y_{N-2}}$ and $\frac{\del F}{\del J_{N-2}}$,
  and continue by backward induction as necessary.  If all $J_m=1$,
  which is the most degenerate isentropic case, then since
  $h(1,\ol Y_{N-1})=\ol Y_{N-1}$, we get $F = 0$ identically.

  We have thus shown that, as long as there is at least one nontrivial
  jump $J_m\ne 1$, so that the flow is not everywhere isentropic, it
  follows that
  \[
    \nabla_{(J,\ol Y_{N-1})} F = 0 \com{implies} F(J,\ol Y_{N-1}) = 0,
  \]
  which implies \eqref{gradF}.  This in turn implies that $DG_q$ has
  full rank, and completes the proof.
\end{proof}

\subsection{Genericity of Nonresonant Profiles}

Even though the space of resonant profiles is dense, these should be
non-generic, because each satisfies at least one, and at most
countably many, constraints \eqref{wres}.  However, proving this in
the full set of possible $L^1$ entropy profiles is problematic,
because we have not fully anayzed the linearized operator for general
entropy profiles in $L^1$.  Since every $L^1$ function can be
approximated in $L^1$ by piecewise constant step functions, Theorem
\ref{thm:Z} does however provide a rigorous sense in which the set of
resonant entropy profiles is non-generic.

\begin{corollary}
  \label{cor:pwc0}
  For each $N>1$, the subset of fully nonresonant piecewise constant
  entropy profiles has full measure in the finite dimensional set
  $\mc P_N$.  Specifically, the Lebesgue measure of the resonant set
  $\mc Z_N$ given in \eqref{Zres} is zero.  In this sense, the
  \emph{fully nonresonant} entropy profiles are generic in the space
  of piecewise constant profiles.
\end{corollary}

\begin{proof}
  Restricting to $\mc P_N$, any entropy profile which has a resonant
  pair of linear modes must satisfy the constraint \eqref{wres} for
  some triple of indices.  We therefore define the resonant set to be 
  \[
    \mc Z_N = \bigcup_{j,k,l}\mc Z^N_{j,k,l},
  \]
  which is a countable union.  Since by Theorem \ref{thm:Z}, each of
  these $\mc Z^N_{j,k,l}$ is a $C^\infty$ manifold of codimension one,
  it has Lebesgue measure zero.  Thus the countable union $\mc Z^N$,
  which is the set of all resonant profiles with $N-1$ jumps, has
  measure zero, and the set of fully nonresonant profiles has full
  Lebesgue measure in $\mc P_N$.
\end{proof}

This in turn provides a sense in which the nonresonant profiles are
generic among all entropy profiles in $BV$ or $L^1$.

\begin{theorem}
  \label{thm:gen}
  The probability that any one piecewise constant approximation of an
  $L^1$ entropy profile has resonant linearized modes is zero, in the
  following sense.  Given any piecewise constant approximation to an
  $L^1$ entropy profile, the probability of that approximation having
  any resonant linearized modes is zero, within the class of all
  piecewise constant profiles with the same number of jumps.
\end{theorem}

\section{Context and Conclusions}
\label{sec:con}

The modern theory of continuum mechanics arguably began in Book I of
Newton's \emph{Principia}, in which set out his famous Laws of Motion.
In an effort to understand resistive forces, in Book II Newton
attempted to describe both the dynamics of continuous media and the
propagation of sound waves.  In the 1750's, Euler developed the
correct extension of Newton's laws of motion to the continuum, and
then linearized the equations to produce the wave equation which
D'Alembert had earlier derived to describe infinitesimal displacements
of a vibrating string.  By this Euler provided a \emph{mechanical
  explanation for music}: vibrations of an instrument produce
sinusoidal oscillations in air pressure, frequencies of which
correspond to the pure tones of sound we hear when, say, a violin is
played.

In the mid-19th century, mathematicians including Stokes and Riemann
discovered a problem with this theory: solutions containing
compressions could not be sustained, and shock formation appeared to
be unavoidable, and inconsistent with the musical tones of the linear
theory.  Challis identified the issue in 1848, after which Stokes
argued that oscillations break down in finite time~\cite{Stokes}.
Riemann's proof that compression \emph{always} produces shock waves in
\emph{isentropic} flows was made definitive in the celebrated
Glimm-Lax decay result of 1970, that periodic solutions of \x2 systems
necessarily form shocks and decay to average at rate $1/t$.  These
results reinforced the belief that the nonlinear and linearized
solutions apparently produce qualitatively different phenomena for
oscillatory sound waves -- the nonlinear system produces shocks while
the linear system does not.

The results in this paper introduce a new nonlinear physical
phenomenon which has the effect of of bringing the linear and
nonlinear theories into alignment when the entropy profile is
non-constant.  Namely, characteristics move ergodically through the
periods, bringing compression and rarefaction into balance and
avoiding shock formation, by sampling regions of both equally, on
average.  Alternatively, the pattern of waves reflected by the entropy
gradients, or ``echoes'' has the effect of attenuating compression and
rarefaction in the wave family opposite the transmitted wave, on
average.  Instead of Riemann invariant coordinates propagating as
constant along characteristics (as in \x2 isentropic and barotropic
systems \cite{S}), in this new \x3 regime, every characteristic cycles
through, and hence samples, a dense set of values of each Riemann
invariant.  This is a new point of view on shock free wave propagation
in compressible Euler.

Our results raise the interesting question as to whether solutions
with shocks fall in the actual regime of ordinary sounds of speech and
musical tones heard in nature and everyday life.  Glimm-Lax theory,
and indeed most existence theories, are based on approximating smooth
solutions by weak shock waves, but as far as we can tell, only strong
shocks are typically observed in nature.  Our results and the success
of the field of Acoustics strongly indicate that this regime of
nonlinear shock-free wave propagation could well be more fundamental
to ordinary sounds and musical tones than formation and propagation of
``weak shock waves’’.  Interestingly, we note that equi-temperament
tuning of the piano makes frequencies irrationally related, which is
precisely our \emph{non-resonance} condition, sufficient to imply
perturbation of linear pure tones to nonlinear pure tones.

The essential physical ideas in this paper were first understood by
the authors within the context of the theory of nonlinear wave
interactions introduced by Glimm and Lax in~\cite{G,GL}.  The
interaction of a nonlinear (acoustic) wave of strength $\gamma$ with
an entropy jump $[s]$ produces an ``echo'', which is a reflected
acoustic wave, whose strength is $O([s]\,\gamma)$, on the order of the
incident acoustic wave~\cite{TY}.  On the other hand, the interaction
of any two (weak) acoustic waves is linear, with an error which is
\emph{cubic} in wave strength.  We began this project with the insight
that, therefore, \emph{the echoes produced by finite entropy jumps are
  at the critical order sufficient to balance rarefaction and
  compression}.  By this we might also expect that a theory of
nonlinear superposition of the ``pure tone'' nonlinear sound waves
constructed here could also produce perturbations which provide
general shock-free solutions of the nonlinear equations, although
these would no longer be periodic in time.  Mathematically, this
raises the question as to whether, by the same mechanism,
quasi-periodic mixed modes of the linearized theory also perturb to
nonlinear.

We remark that the nonlinear effect of compression forming into shock
waves, is strongest in one spatial dimension, because geometrical
effects mitigate steepening of waves in higher dimensions.  Thus the
study of effects which attenuate shock formation get more to the heart
of the issue in one dimension.  Since the time of Stokes and Riemann,
there has been a long history of study of the mechanisms which lead to
shock formation.  The authors suggest that the greater scientific
mystery now might lie in the nonlinear mechanisms which \emph{prevent}
shock formation.

Finally, regarding mathematical methods, our results demonstrate that
the problem of expunging resonances inherent in many applications, can
be overcome when enough symmetries are present to impose periodicity
by \emph{projection}, rather than by \emph{periodic return}.  Compared
to Nash-Moser techniques required to prove existence of periodic
solutions in simpler settings, we find here that taking into account
the essential symmetries has led to the discovery of a simpler method
for solving a more complicated problem.

 \providecommand{\url}[1]{{\tt #1}}

\end{document}